%% file: main.tex
\documentclass[reqno,twoside]{article} 

\input{settings.tex}  

\def\essayauthor{Mats Wahlberg}

\def\essaytitle{The Banach-Tarski Paradox}

\def\subtitle{How I Learned to Stop Worrying and Love the Axiom of Choice}

\def\semester{Spring}

\def\essayyear{2022}

\def\credits{15}

\def\degree{Bachelor in Mathematics}

\def\degreecredits{180}

\def\department{Department of Mathematics and Mathematical Statistics}

\begin{document}

\newtitleGM 	

\selectlanguage{english}
\begin{abstract}
	\noindent This thesis presents the strong and weak forms of the Banach-Tarski paradox
	based on the Hausdorff paradox. It provides modernized proofs of the paradoxes and
	necessary properties of equidecomposable and paradoxical sets. The
	historical significance of the paradox for measure theory is covered, along with its
	incorrect attribution to Banach and Tarski.
	Finally, the necessity of the axiom of choice is discussed and contrasted with
	other axiomatic and topological assumptions that enable the paradoxes.
	\\[4\baselineskip]
\end{abstract}

\selectlanguage{swedish}

\begin{abstract}
	\noindent Den här uppsatsen presenterar den starka och svaga formen av Banach-Tarskis
	paradox baserade på Hausdorffs paradox. Den tillhandahåller moderniserade bevis av
	paradoxerna och nödvändiga egenskaper av likuppdelningsbara och paradoxala mängder.
	Den historiska betydelsen av paradoxen på måtteori tas upp samt dess felaktiga
	tillskrivning till Banach och Tarski.
	Till sist diskuteras behovet av urvalsaxiomet som ställs i kontrast mot andra
	axiomatiska och topologiska antaganden som möjliggör paradoxerna.
\end{abstract}

\selectlanguage{english}

\cleardoublepage
\thispagestyle{empty}
\tableofcontents
\vspace{0.7in}

\cleardoublepage
\setcounter{page}{1}
\pagestyle{plain}
\pagenumbering{arabic}
		
\include{introduction}

\include{preliminaries}

\include{banach-tarski}

\include{hausdorff}

\include{proof}

\include{discussion}

\section{References}

	\begingroup
		\renewcommand{\section}[2]{}%
		\bibliography{refs}
	\endgroup



\end{document}

%% file: settings.tex
\newcommand{\SPH}{\mathbb{S}^2}
\newcommand{\CIRC}{\mathbb{S}^1}
\newcommand{\SPHMD}{\mathbb{S}^2 \setminus D}
\newcommand{\BLL}{\mathbb{B}^3}
\newcommand{\BLLMC}{\mathbb{B}^3 \setminus\{\vect{0}\}}

\newcommand{\ACC}{\nameref{ax_choice}\xspace}
\newcommand{\WRD}[1]{\texttt{#1}}
\newcommand{\PWR}[1]{\textsuperscript{#1}}
\newcommand{\INV}{\PWR{-1}}

\usepackage[shortcuts]{extdash}
\usepackage{enumerate}

\usepackage{cancel}
\usepackage{braket}
\usepackage{xspace}

\usepackage{amsmath}
\usepackage[linktoc=all]{hyperref}
\usepackage[nameinlink]{cleveref}
\usepackage[alphabetic]{amsrefs}
\usepackage{siunitx}
\usepackage{tkz-euclide}
\usepackage{amsthm}

\crefname{prop}{proposition}{propositions}
\crefname{cor}{corollary}{corollaries}
\crefname{claim}{claim}{claims}
\crefname{enumi}{quote}{quotes}
\crefname{axiom}{axiom}{axioms}
\crefname{equation}{}{}

\newtheorem{claim}{Claim}

\usepackage{etoolbox}
\AtEndEnvironment{proof}{\setcounter{claim}{0}}

\crefname{examplex}{example}{examples}
\crefname{remarkx}{remark}{remarks}

\usepackage[utf8]{luainputenc}
\usepackage[LGR,T1]{fontenc}
\DeclareMathAlphabet{\mathgtt}{LGR}{cmtt}{m}{n}
\newcommand{\GWRD}[1]{$\mathgtt{#1}$}

\usepackage[english,swedish]{babel}
\usepackage[margin=0.5in,twoside]{geometry}               
\usepackage{graphicx,
            epstopdf,
            float,
            subcaption}             
\usepackage[section]{placeins}      
\usepackage{
            amssymb,
            amsthm,
            amsfonts,
            mathrsfs,
            dsfont,
            esint,
            mathtools,
            mathrsfs}               
\usepackage{fancyhdr}               
\usepackage{booktabs,
            multirow}               
\usepackage{appendix}               
\usepackage{csquotes}               
\usepackage{color}                  
\usepackage{pdflscape}              
\usepackage{multicol}               
\usepackage{pgfplots}               
\usepackage{lmodern}                
\usepackage{newfloat}               
\usepackage{bm}
\usepackage{physics}
\usepackage{wasysym}
\usepackage{thmtools}

\let\oldsection\section
\def\section{\cleardoublepage\oldsection}

\makeatletter
    \newtoggle{biblatexloaded}
    \newtoggle{mintedloaded}
    \@ifpackageloaded{biblatex}{\toggletrue{biblatexloaded}}{\togglefalse{biblatexloaded}}
    \@ifpackageloaded{minted}{\toggletrue{mintedloaded}}{\togglefalse{mintedloaded}}
\makeatother

\pgfplotsset{compat=newest} 
\usetikzlibrary{patterns,
				arrows,
				positioning,
				automata,
				calc}

\iftoggle{biblatexloaded}{
    \DeclareNameAlias{sortname}{last-first}
    \DeclareNameAlias{default}{last-first}

    \urlstyle{sf}

}

\iftoggle{mintedloaded}{
    \newmintedfile{c}{fontfamily=tt,
                        fontsize=\footnotesize,
                        tabsize=4,
                        numberblanklines=true,
                        numbers=left,
                        numbersep=5pt,
                        breakautoindent=false,
                        xleftmargin=0.7cm,} 

    \newminted[ccodecap]{c}{fontfamily=tt,
                            fontsize=\normalsize,
                            tabsize=4,
                            frame=lines,
                            breaklines=true,
                            breaksymbolleft=\tiny\ensuremath{\hookrightarrow},
                            breakautoindent=true,} 
    \newmintinline{c}{fontfamily=tt,
                        breaklines=true,
                        fontsize=\normalsize} 

    \newminted[shellcode]{shell-session}{fontfamily=tt,
                                            fontsize=\normalsize,
                                            tabsize=4,} 

    \newminted[shellcodecap]{shell-session}{fontfamily=tt,
                                            fontsize=\normalsize,
                                            tabsize=4,
                                            frame=lines,
                                            breaklines=true,
                                            breaksymbolleft=\tiny\ensuremath{\hookrightarrow},
                                            breakautoindent=true} 
    \newmint[shell]{shell-session}{fontfamily=tt,
                                    fontsize=\normalsize,
                                    frame=lines,
                                    breaklines=true,
                                    breaksymbolleft=\tiny\ensuremath{\hookrightarrow},
                                    breakautoindent=true} 
    \newmintedfile[shellfile]{shell-session}{fontfamily=tt,
                                                fontsize=\normalsize,
                                                tabsize=4,
                                                breakautoindent=false,} 
    \newmintinline[shellinline]{shell-session}{fontfamily=tt,
                                                fontsize=\normalsize,
                                                breaklines=true,
                                                breakautoindent=false} 
    \newmintedfile{matlab}{fontfamily=tt,
                            fontsize=\footnotesize,
                            tabsize=4,
                            numberblanklines=true,
                            numbers=left,
                            numbersep=5pt,
                            breakautoindent=false,
                            xleftmargin=0.7cm,} 
    \newmintedfile[matlabfileframe]{matlab}{fontfamily=tt,
                                            fontsize=\footnotesize,
                                            tabsize=4,
                                            frame=lines,
                                            numberblanklines=true,
                                            numbers=left,
                                            numbersep=5pt,
                                            breakautoindent=false,
                                            xleftmargin=0.7cm,} 

    \newminted[mcodecap]{matlab}{fontfamily=tt,
                                    fontsize=\normalsize,
                                    tabsize=4,
                                    frame=lines,
                                    breaklines=true,
                                    breaksymbolleft=\tiny\ensuremath{\hookrightarrow},
                                    breakautoindent=true,} 
    \newmintinline[minline]{matlab}{fontfamily=tt,
                                    breaklines=true,
                                    fontsize=\normalsize} 

    \SetupFloatingEnvironment{listing}{name=Code} 
}

\usepackage[margin=3ex,
			font=small,
			labelfont=bf,
			labelsep=endash]{caption}

\hypersetup{colorlinks=true,
			hidelinks}

\newcommand{\vect}[1]{\boldsymbol{\mathbf{#1}}} 

\let\originalleft\left
\let\originalright\right
\renewcommand{\left}{\mathopen{}\mathclose\bgroup\originalleft}
\renewcommand{\right}{\aftergroup\egroup\originalright}

\graphicspath{{./figs/}}

\numberwithin{equation}{section}
\numberwithin{figure}{section}

\renewcommand{\emptyset}{\varnothing}

\newtheorem{theorem}{Theorem}[section]
\newtheorem{lemma}[theorem]{Lemma}
\newtheorem{prop}[theorem]{Proposition}
\newtheorem{cor}[theorem]{Corollary}
\theoremstyle{definition}
\newtheorem{definition}[theorem]{Definition}
\newtheorem{examplex}[theorem]{Example}

 \newenvironment{example}
  {\pushQED{\qed}\renewcommand{\qedsymbol}{$\diamondsuit$}\examplex}
  {\popQED\endexamplex}

\theoremstyle{remark}
\newtheorem{remarkx}[theorem]{Remark}

 \newenvironment{remark}
  {\pushQED{\qed}\renewcommand{\qedsymbol}{$\diamondsuit$}\remarkx}
  {\popQED\endremarkx}

\newcommand{\RE}{\mbox{$\mathbb{R}$}}
\newcommand{\Z}{\mbox{$\mathbb{Z}$}}

\newcommand{\Q}{\mbox{$\mathbb{Q}$}}

\newcommand{\N}{\mbox{$\mathbb{N}$}}
\newcommand{\A}{\mbox{$\mathscr{A}$}}

\newcommand{\B}{\mbox{$\mathscr{B}$}}

\newenvironment{bottompar}{\par\vspace*{\fill}}{\vspace{3\baselineskip}\clearpage}
\author{\essayauthor}

\newcommand*{\newtitleGM}{
	\pagenumbering{gobble}
	\pagestyle{empty}
	\begingroup\begin{center}
  		\includegraphics[width=.25\textwidth]{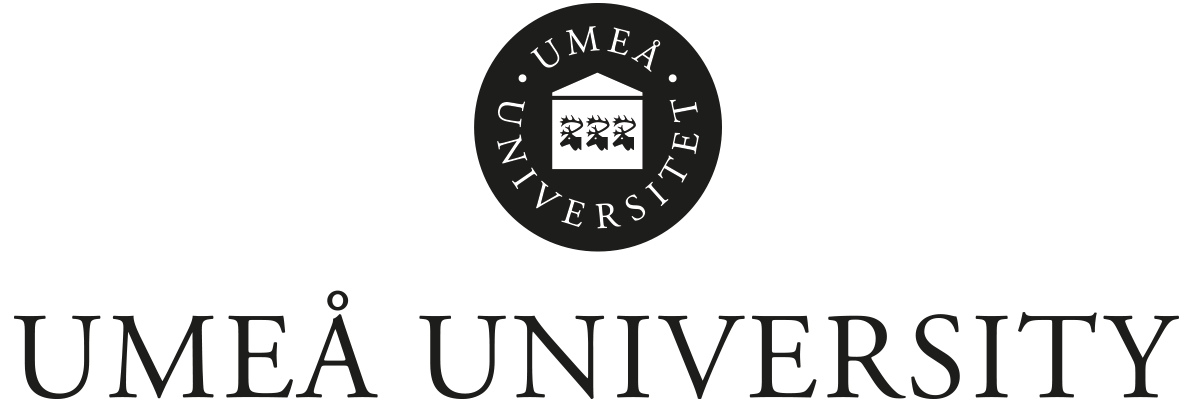} 
  		\\[6\baselineskip]
		\vfill	  	
	  	{\Huge\bfseries\essaytitle}\\[\baselineskip]
		{\large\slshape\bfseries\subtitle}\\[\baselineskip]
  	  	{\large\essayauthor}
  	  	\vfill
		\begin{bottompar}
			\noindent
			Bachelor Thesis, \credits hp\\
			\degree, \degreecredits hp\\
    		\semester\ \essayyear\\
			\department
 		\end{bottompar}
	\end{center}\endgroup
	\newgeometry{top=1.75in, bottom=1.75in, right=1.6in, left=1.6in}
	\cleardoublepage
} 

%% file: introduction.tex
\section{Introduction}

\renewenvironment{quote}{%
  \list{}{%
    \leftmargin1.22cm   
    \rightmargin\leftmargin
  }
  \item\relax
}
{\endlist}

Humanity has always been fascinated by paradoxes: seemingly valid reasonings on acceptable premises
that give unreasonable results. In mathematics, paradoxes have been the cause for particular
concern as they seem to question the trust we put in our results and the assumptions on which they
are built. Sometimes, these results turn out to be either antinomous: contradictions such as
Russell’s paradox, or falsidical: fallacies in the reasoning like Zeno’s infamous paradoxes.

However, there are also veridical paradoxes: valid results that challenge our expectations.
Generally, only these paradoxes are of interest within the mathematical community, and for a good
reason. Veridical paradoxes have been explored since the early days of mathematics, from famous
results by Galileo to significantly older. However, few have gained such recognizability as the
\textit{Banach-Tarski paradox}, reaching well beyond the mathematical community and given such
colorful interpretations as the \enquote{pea and Sun paradox} by the general public.

\begin{figure}[h]
	\centering
	\def\svgscale{0.5}
	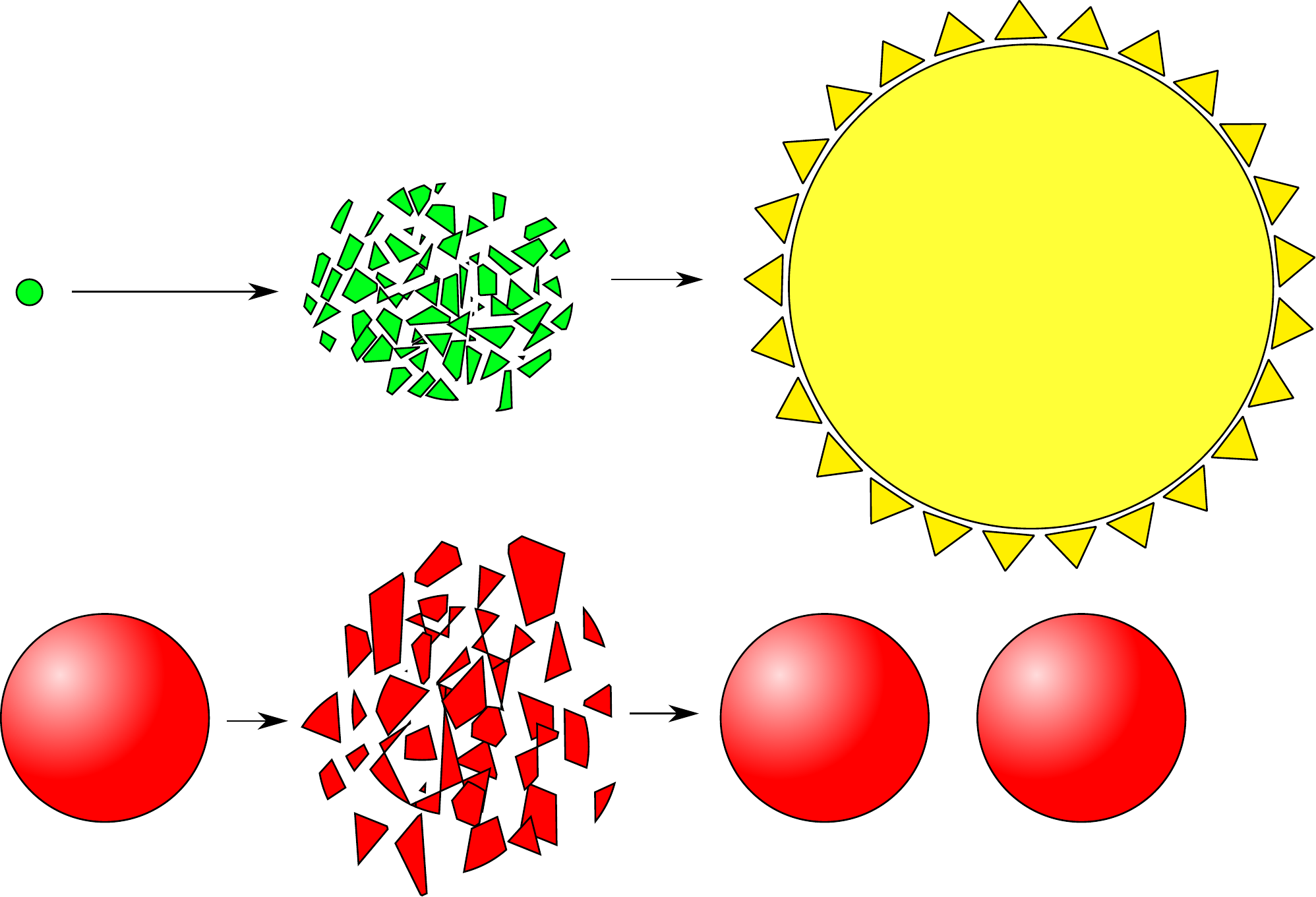
	\caption{Popular illustrations of the Banach-Tarski paradox, involving duplication of
	objects and constructing vastly different objects from each other.}
	\label{fig:illustration}
\end{figure}

Some popular descriptions of the Banach-Tarski paradox, illustrated in \Cref{fig:illustration}, are
as follows:
\begin{quote}
	\begin{enumerate}[Quote 1.]
		\item \enquote{A pea can be chopped up and reassembled into the Sun.\label{quote1}}
		\item \enquote{A ball can be divided into a finite number of pieces and reassembled
			into two identical copies of the original ball.\label{quote2}}
	\end{enumerate}
\end{quote}

Although expressed in layman's terms, these two descriptions also illustrate the two different
mathematical forms of the paradox, sometimes referred to as its strong and weak forms by
\Cref{quote1,quote2}, respectively. The weak form is usually the one referred to as
the Banach-Tarski paradox. Mathematical definitions of both forms are provided in \Cref{sec:bt},
where we also will see that they are equivalent. The weak form offers us an easier avenue for
proving both forms. However, it also has a curious history by originally
being presented by Felix Hausdorff ten years prior to the work of Stefan Banach and Alfred Tarski.

Hausdorff's result has remained relatively unknown to this day and was not credited by Banach and
Tarski even
though they credit their use of an intermediary step from Hausdorff's work, on which they apply
a series of constructions identical to what Hausdorff originally presented. Ironically, this intermediary
result is today known as the \textit{Hausdorff paradox}. Of everything presented in this
thesis, the Hausdorff paradox is still by far the single most significant result with the most
sophisticated proof. Its significance has inspired this thesis to dedicate a section,
\Cref{sec:hausdorff}, to it. How the Hausdorff paradox leads to the
Banach-Tarski paradox is covered in \Cref{sec:proof}. Further historical discussion on its
attribution is in \Cref{sec:attribution}, and \Cref{sec:non-measurable} covers the historical significance of the result.

A formal mathematical statement of the strong form of the Banach-Tarski paradox is as follows:
\begin{restatable*}[The strong Banach-Tarski paradox]{theorem}{btprdx}
	\label{b-t1}
	Let $A$ and $B$ be bounded subsets of $\RE^3$ with nonempty interiors, then $A$ and $B$ are
	$G_3$\-/equidecomposable.
\end{restatable*}

It is this result to which this thesis is dedicated. From the basic definitions needed for its
formulation and interpretation to its construction, history and to what extent the result and
assumptions it
is based on can be questioned. It is worth noting that the Banach-Tarski paradox is veridical, but
only if we accept the \nameref{ax_choice}. This axiom has a long history of controversies among
mathematicians but is these days generally accepted. However, it is still common to see the
\nameref{ax_choice} blamed as the sole cause of the Banach-Tarski paradox. \Cref{disc_ac} provides some arguments as to why
this reasoning overlooks far more critical flaws in our core assumptions, enabling
many other paradoxes independently of the \nameref{ax_choice}. It also covers how these assumptions can be
changed to prevent paradoxes without discarding the \nameref{ax_choice}.

Although the result is often discussed from its analytical implication, its proof relies on the
time-honored tradition of unifying seemingly different fields of mathematics to achieve results
beyond what each field could do on its own. In some cases, it has also brought something back to the
individual fields, from the immediate result for measure theory to sparking more research into
topics like free groups and pointless topologies. The presented proofs combine fields such as
analysis, linear algebra, group theory, set theory and even discrete mathematics.

One of the most critical proofs, that of \Cref{G_6}, uses modular arithmetics, which is a modern
approach not part of the original proofs. For the same proof, this thesis utilizes a Pythagorean
triple in an initial construction and observes the potential of a simple generalization to any
primitive Pythagorean triple in \Cref{remark:pythagorean}. This use and generalization by
Pythagorean triples is an approach the author has not seen previously documented. The
material is based mainly on \cite{wagonbok} and \cite{kursbok}, while all other sources are cited where used.

A central concept in the work of Hausdorff, as well as Banach and Tarski, was the study of the
new field of measure theory, specifically in proving the existence of non-measurable sets, which
will be discussed in \Cref{sec:non-measurable}. A clear inspiration comes from a result by Giuseppe
Vitali, whose well-known \textit{Vitali set} had already shown such existence. The result of
Hausdorff is significantly stronger, but much of the construction shares significant similarities to
that of Vitali, who had successfully applied his reasoning to a circle obtaining the paradoxical
ability to duplicate it in much the same way as the other paradoxes, as shown in \Cref{ex:vitali}.

A common theme for all paradoxes that results in non-measurable sets is the necessity of the
\nameref{ax_choice}. Nevertheless, the are many other kinds of paradoxes in mathematics that do not require
it. A simple and well-known paradoxical result can be obtained from the concept of
\textit{irrational rotation}.

The following example might look out of place for an introduction. However, it is motivating to see that
we can achieve easy-to-understand paradoxical results without too much work, heavy
theory, or usage of the \nameref{ax_choice}. We will repeatedly return to this example
to illustrate, motivate and contrast various ideas and properties we will encounter. Furthermore,
it will even return in a generalized form as part of a vital lemma in \Cref{lm:irrational}.

	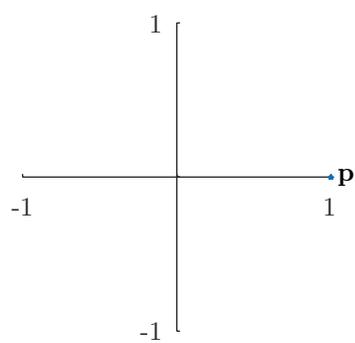
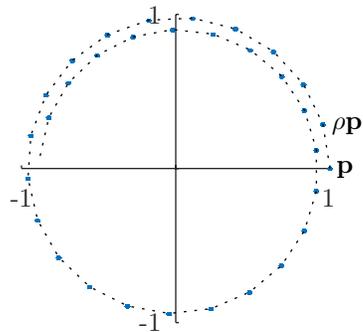
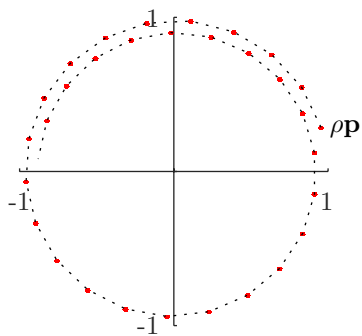
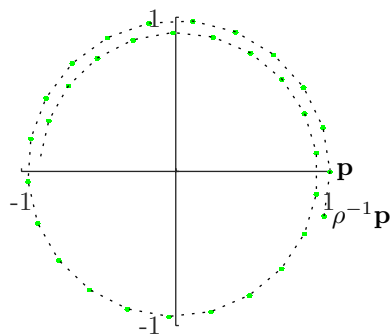
\begin{figure}
		\centering
		\begin{subfigure}{0.45\textwidth}
			\centering
			\input{figs/irrational1.tex}
			\caption{The point $\vect{p}$.}
			\label{fig:irrational1}
		\end{subfigure}
		\quad
		\begin{subfigure}{0.45\textwidth}
			\centering
			\input{figs/irrational2.tex}
			\caption{Construction of $X$ as the set of repeated rotations
			of $\vect{p}$ by $\rho$.}
			\label{fig:irrational2}
		\end{subfigure}

		\begin{subfigure}{0.45\textwidth}
			\centering
			\input{figs/irrational3.tex}
			\caption{$X$ rotated by $\rho$, $\rho X$, is identical to $X$
			except for the missing point $\vect{p}$.}
			\label{fig:irrational3}
		\end{subfigure}
		\quad
		\begin{subfigure}{0.45\textwidth}
			\centering
			\input{figs/irrational4.tex}
			\caption{$X$ rotated by $\rho^{-1}$, $\rho^{-1} X$, is identical to $X$
			with the additional point $\rho^{-1}\vect{p}$.}
			\label{fig:irrational4}
		\end{subfigure}
		\caption{Construction and rotation of $X$. Only the first 33 points are drawn and
		the points are spiraling towards the center only to improve visual clarity.}
	\end{figure}

\begin{example}[Irrational rotation]
	\label{ex:irrational}

	Let $\theta=2\pi a$ where $a$ is irrational. Let $\rho$ be the rotation by $\theta$
	around origo in the $xy$-plane in $\RE^2$. This type of transformation is known as an
	\textit{irrational rotation}. Informally we call $\theta$ an \textit{irrational angle}.

	Consider the point $\vect{p}=(1,0)$, as in \Cref{fig:irrational1}. Every repeated rotation
	of $\vect{p}$ by $\rho$ must result in a unique point. In other words, none of the points
	obtained will have periodic orbits in $\rho$: otherwise, there are $n,m\in\N$ such that
	$\rho^m\vect{p}=\rho^n\vect{p}$, which gives that $\rho^k\vect{p}=\vect{p}$ where
	$k=\abs{m-n}$, and thus $\theta k$ must be an integer multiple of $\pi$.
	This would contradict our choice of $\theta$.

	Now consider the set of all these points, $X=\bigcup_{n=0}^\infty\rho^n\vect{p}$,
	as illustrated in \Cref{fig:irrational2}. If we apply our rotation $\rho$ to the whole set
	$X$, this means each point gets moved to the position of the next, resulting in $\rho
	X=X\setminus\vect{p}$, as shown in \Cref{fig:irrational3}. Hence by doing nothing more than
	rotating $X$, we can remove points
	from it. Similarly, we can add new points by
	rotating $X$ in the opposite direction, by $\rho^{-1}$, as illustrated in \Cref{fig:irrational4}
\end{example}
This paradox marks our first mathematical deviation from intuition. Unlike
the Hausdorff and Banach-Tarski paradoxes, it does not depend on the \nameref{ax_choice}, but just
like them, it does depend on the notion of points and that sets of infinite size do exist, the
\textit{axiom of infinity}. These concepts turn out to be more noteworty for the paradoxical
results than the \nameref{ax_choice}.

%% file: figs/illustration.pdf_tex
\begingroup%
  \makeatletter%
  \providecommand\color[2][]{%
    \errmessage{(Inkscape) Color is used for the text in Inkscape, but the package 'color.sty' is not loaded}%
    \renewcommand\color[2][]{}%
  }%
  \providecommand\transparent[1]{%
    \errmessage{(Inkscape) Transparency is used (non-zero) for the text in Inkscape, but the package 'transparent.sty' is not loaded}%
    \renewcommand\transparent[1]{}%
  }%
  \providecommand\rotatebox[2]{#2}%
  \newcommand*\fsize{\dimexpr\f@size pt\relax}%
  \newcommand*\lineheight[1]{\fontsize{\fsize}{#1\fsize}\selectfont}%
  \ifx\svgwidth\undefined%
    \setlength{\unitlength}{574.56444769bp}%
    \ifx\svgscale\undefined%
      \relax%
    \else%
      \setlength{\unitlength}{\unitlength * \real{\svgscale}}%
    \fi%
  \else%
    \setlength{\unitlength}{\svgwidth}%
  \fi%
  \global\let\svgwidth\undefined%
  \global\let\svgscale\undefined%
  \makeatother%
  \begin{picture}(1,0.68173462)%
    \lineheight{1}%
    \setlength\tabcolsep{0pt}%
    \put(0,0){\includegraphics[width=\unitlength,page=1]{illustration.pdf}}%
  \end{picture}%
\endgroup%

%% file: figs/irrational1.tex
\setlength{\unitlength}{1pt}
\begin{picture}(0,0)
\includegraphics[scale=1]{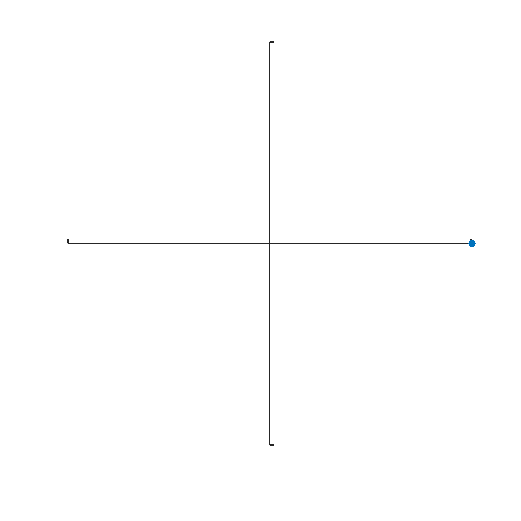}
\end{picture}%
\begin{picture}(150,150)(0,0)
\fontsize{10}{0}\selectfont\put(19.5,72.0612){\makebox(0,0)[t]{\textcolor[rgb]{0.15,0.15,0.15}{{-1}}}}
\fontsize{10}{0}\selectfont\put(135.75,72.0612){\makebox(0,0)[t]{\textcolor[rgb]{0.15,0.15,0.15}{{1}}}}
\fontsize{10}{0}\selectfont\put(72.4096,21.7593){\makebox(0,0)[r]{\textcolor[rgb]{0.15,0.15,0.15}{{-1}}}}
\fontsize{10}{0}\selectfont\put(72.4096,138.009){\makebox(0,0)[r]{\textcolor[rgb]{0.15,0.15,0.15}{{1}}}}
\fontsize{10}{0}\selectfont\put(138.656,79.8843){\makebox(0,0)[l]{\textcolor[rgb]{0,0,0}{{$\vect{p}$}}}}
\end{picture}

%% file: figs/irrational2.tex
\setlength{\unitlength}{1pt}
\begin{picture}(0,0)
\includegraphics[scale=1]{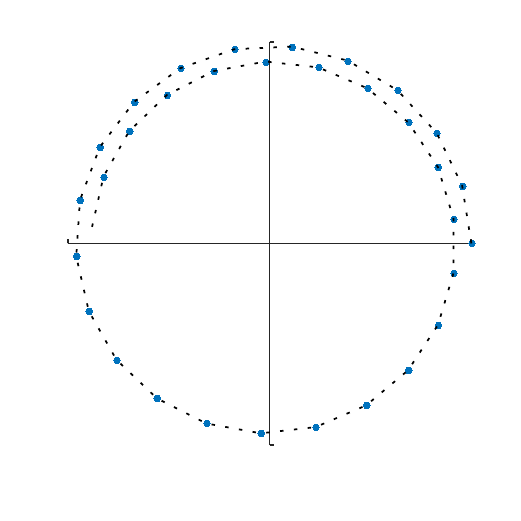}
\end{picture}%
\begin{picture}(150,150)(0,0)
\fontsize{10}{0}\selectfont\put(19.5,72.0612){\makebox(0,0)[t]{\textcolor[rgb]{0.15,0.15,0.15}{{-1}}}}
\fontsize{10}{0}\selectfont\put(135.75,72.0612){\makebox(0,0)[t]{\textcolor[rgb]{0.15,0.15,0.15}{{1}}}}
\fontsize{10}{0}\selectfont\put(72.4096,21.7593){\makebox(0,0)[r]{\textcolor[rgb]{0.15,0.15,0.15}{{-1}}}}
\fontsize{10}{0}\selectfont\put(72.4096,138.009){\makebox(0,0)[r]{\textcolor[rgb]{0.15,0.15,0.15}{{1}}}}
\fontsize{10}{0}\selectfont\put(138.656,79.8843){\makebox(0,0)[l]{\textcolor[rgb]{0,0,0}{{$\vect{p}$}}}}
\fontsize{10}{0}\selectfont\put(136.913,96.1593){\makebox(0,0)[l]{\textcolor[rgb]{0,0,0}{{$\rho \vect{p}$}}}}
\end{picture}

%% file: figs/irrational3.tex
\setlength{\unitlength}{1pt}
\begin{picture}(0,0)
\includegraphics[scale=1]{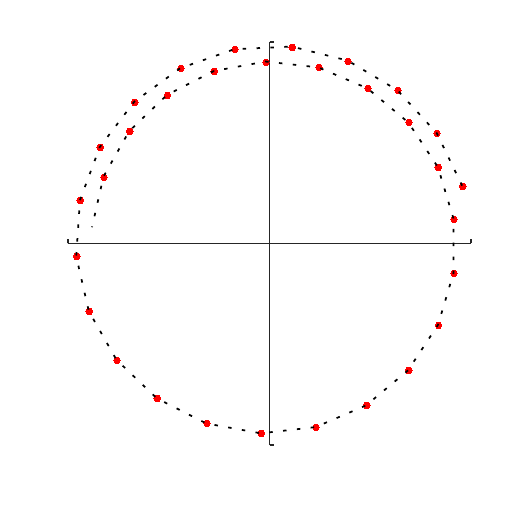}
\end{picture}%
\begin{picture}(150,150)(0,0)
\fontsize{10}{0}\selectfont\put(19.5,72.0612){\makebox(0,0)[t]{\textcolor[rgb]{0.15,0.15,0.15}{{-1}}}}
\fontsize{10}{0}\selectfont\put(135.75,72.0612){\makebox(0,0)[t]{\textcolor[rgb]{0.15,0.15,0.15}{{1}}}}
\fontsize{10}{0}\selectfont\put(72.4096,21.7593){\makebox(0,0)[r]{\textcolor[rgb]{0.15,0.15,0.15}{{-1}}}}
\fontsize{10}{0}\selectfont\put(72.4096,138.009){\makebox(0,0)[r]{\textcolor[rgb]{0.15,0.15,0.15}{{1}}}}
\fontsize{10}{0}\selectfont\put(136.913,96.1593){\makebox(0,0)[l]{\textcolor[rgb]{0,0,0}{{$\rho \vect{p}$}}}}
\end{picture}

%% file: figs/irrational4.tex
\setlength{\unitlength}{1pt}
\begin{picture}(0,0)
\includegraphics[scale=1]{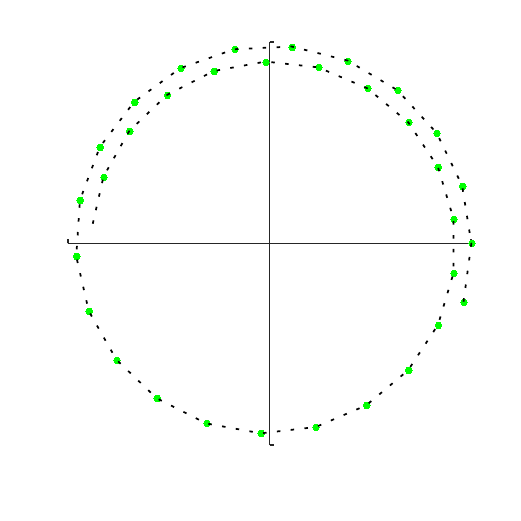}
\end{picture}%
\begin{picture}(150,150)(0,0)
\fontsize{10}{0}\selectfont\put(19.5,72.0612){\makebox(0,0)[t]{\textcolor[rgb]{0.15,0.15,0.15}{{-1}}}}
\fontsize{10}{0}\selectfont\put(135.75,72.0612){\makebox(0,0)[t]{\textcolor[rgb]{0.15,0.15,0.15}{{1}}}}
\fontsize{10}{0}\selectfont\put(72.4096,21.7593){\makebox(0,0)[r]{\textcolor[rgb]{0.15,0.15,0.15}{{-1}}}}
\fontsize{10}{0}\selectfont\put(72.4096,138.009){\makebox(0,0)[r]{\textcolor[rgb]{0.15,0.15,0.15}{{1}}}}
\fontsize{10}{0}\selectfont\put(138.656,79.8843){\makebox(0,0)[l]{\textcolor[rgb]{0,0,0}{{$\vect{p}$}}}}
\fontsize{10}{0}\selectfont\put(136.913,63.6093){\makebox(0,0)[l]{\textcolor[rgb]{0,0,0}{{$\rho^{-1} \vect{p}$}}}}
\end{picture}

%% file: preliminaries.tex
\section{Preliminaries}

The \ACC is needed in one proof, of \Cref{G_7}.
\begin{restatable}[axiom of choice]{axiom}{axchoice}
	\label{ax_choice}
	For any nonempty collection of disjoint nonempty sets, there is a set containing exactly
	one element from each set in the collection.
\end{restatable}
Of course, this is not the only axiom we use. We always use many axioms implicitly, such as
Zermelo–Fraenkel set theory. However, it deserves special attention as it often gets used without being
explicitly named and has a history of being questioned and often blamed for various controversial
results.

In the case of a finite collection, this property follows directly from Zermelo–Fraenkel set theory.
The collection of sets is allowed to be uncountably infinite. There are weaker versions, for
instance, restricted to only countably infinite collections, but they will be too weak for our use.

The \ACC is necessary to construct the Hausdorff paradox, and thus the Banach-Tarski paradox. It is
also needed when constructing a similar paradox by Vitali in \Cref{ex:vitali}. Even though it is
generally accepted, it is still common to see the
\ACC singled out as the reason for
these paradoxical results, often due to its non-constructive nature.

However, there is a far more
nuanced
discussion to be had on what assumptions make these paradoxes possible. We will see in
\Cref{disc_ac} that there are further, often overlooked assumptions that have received far less
scrutiny. And by changing them, we would be able to keep the \nameref{ax_choice} together with the critical
results it gives while still eliminating most paradoxes of this kind, including several paradoxes
that can be constructed without even using the \nameref{ax_choice}.

The rest of this section is group- and set-theoretic: we will describe objects as sets of points and
we need a proper mathematical description of the concept of moving and rotating them, known as a
\textit{rigid transformation}, which will be defined using \textit{group actions}. Furthermore, we
need to introduce the concept of \textit{equidecomposability} and \textit{paradoxicality} of sets,
and deduce some important properties of them. First, we define two objects that we will
use repeatedly:

\begin{definition}[The unit sphere and ball]
	The unit sphere $\SPH$ and ball $\BLL$, centered at origo in $\RE^3$, are defined by
	$\SPH\coloneqq\set{\vect{x}\in\RE^3 | \norm{\vect{x}}=1}$ and
	$\BLL\coloneqq\set{\vect{x}\in\RE^3 | \norm{\vect{x}}\le1}$.
\end{definition}

\subsection{Rigid Transformations: Rotations and Translations}

To begin, we need to define rotations around origo and will use a common definition:
\begin{definition}[The special orthogonal group in $\RE^3$]
	The $SO(3)$ is the group of all orthogonal $3\times3$ matrices in $\RE^3$ with determinant $1$.
\end{definition}
This group consists of precisely all rotations in $\RE^3$, together with the identity which describes
the trivial case of no rotation. The fact that $SO(3)$ is a group is easy to verify: the identity matrix is in
$SO(3)$, and the inverse of an orthogonal matrix is its inverse, which is still an orthogonal
matrix with determinant $1$ and also belongs to $SO(3)$. The remaining properties, such as
transitivity, follow directly from the usual matrix properties.

The name \textit{special orthogonal group} and notation $SO(3)$ comes from the fact that it is a subgroup of the \textit{orthogonal group}, $O(3)$, which
consists of all orthogonal $3x3$ matrices in $\RE^3$. Remember that an orthogonal matrix means
that its column vectors are orthogonal and normal, or equivalently defined on the row vectors.
Basic linear algebra tells us that the matrices in $O(3)$ must have either determinant $1$ or $-1$.
The matrices with determinant $-1$ result in a reflection, which we do not consider to
be a \textit{proper rotation}, but they belong to a broader concept of isometries and congruence in
geometry as part of the group $O(3)$.

\begin{remark}[SO(3) is the rotations in $\RE^3$]
A complete proof that $SO(3)$ contains all rotations $\RE^3$ and only these rotations can get long and
technical, but a couple of points are worth bringing up to help justify what might otherwise seem
like an arbitrary set of requirements: first, consider what a rotation does: it places a local
coordinate-system within the global.

Three base vectors give this local frame-of-reference as a coordinate system in the global frame-of-reference. Placing these base vectors as columns in a matrix
gives a transformation of coordinates from the local frame to the global frame, and it will be an orthogonal
matrix since the base vectors are orthonormal. And vice versa, all orthonormal matrices describe a
rotation with possible reflection. This is precisely all the $O(3)$ matrices without further requirements.

	For the determinant, consider $A$ to be a rotation and observe that since it's orthogonal
	its inverse is its transpose, $AA^T=I$, and since determinants remain the same after a
	matrix gets transposed, $\det(A^T)=\det(A)$ we see
	$\det(A)^2=\det(A)\det(A^T)=\det(AA^T)=\det(I)=1$, and thus $\det(A)=\pm 1$.

	Now for the sign: any proper rotation can also be performed halfway, and consider $B$ to be
	a rotation of half the angle of $A$, so that $A=BB$. Since $B$ is a rotation it also
	satisfies $\det(B)=\pm1$ and we see $\det(A)=\det(BB)=\det(B)\det(B)=(\pm1)^2=1$.

	We can also obtain that all improper rotations have determinant $-1$, by consider them as a
	one reflection and a combination of proper rotations, for instance the reflection along
	the $x$-axis:
	\[R=\begin{bmatrix}-1&0&0\\0&1&0\\0&0&1\end{bmatrix}\]
	from which we see that $\det(R)=-1$ and the determinant of the total combination of matrices
	is the product of $-1$ and some numbers of $1$.

	Thus $SO(3)$ contains all proper rotations and nothing more.
\end{remark}

From now on, the elements of $SO(3)$ will interchangeably be referred to by matrices or just
plainly described rotations. Another equivalent form of describing all of these rotations is provided
by the following famous result of Leonhard Euler dating back to 1775:
\begin{theorem}[Euler's rotation theorem]
	\label{thm:euler_rotation}
	All rotations from $SO(3)$ can be described by a rotation along an axis passing through
	origo.
\end{theorem}
This will be a useful and intuitive way of considering rotations. This classic result has
received many different proofs, see \cite{palaisart} for one of them, together with some history and
Euler's original geometric proof translated to English. Of course, we could
also use quaternions to describe rotations, but that will not be covered.

We can now easily build on $SO(3)$ by adding translations:
\begin{definition}[The group of rigid transformations in $\RE^3$]
	Let $G_3$ be the collection of all transformations $\vect{x}\mapsto A\vect{x}+\vect{b}$
	where $A\in SO(3)$ and $\vect{b}\in\RE^3$. It follows easily from the fact that $SO(3)$ is
	a group that $G_3$ must also be a group.
\end{definition}
Unlike $SO(3)$, the notation $G_3$ is not an established convention but seems to be a somewhat common
choice in this specific area, probably as a result of its usage in \cite{wagonbok}. It is often
called a \textit{rigid transformation}, or a \textit{rigid body transformation} since its
elements preserve distances between points and does not cause reflections.

Recall that applying group elements on sets is known as \textit{group actions}:
\begin{definition}[Group actions]
	\label{def:groupaction}
	We say that a group $G$ \textit{acts on a set} $X$ when there is an \textit{action}
	defining a mapping $G\times X\to X$, usually written as $gx$ for some $g\in G$ and $x\in X$, which satisfies
	the two axioms $ex=x$ and $f(gx)=(fg)x$ for all $x\in X$ and $f,g\in G$. We can also say
	that $X$ is a \textit{$G$-set}. For convenience we often define $gX\coloneqq\set{gx|x\in X}$ for some $g\in
	G$, $Gx\coloneqq\set{gx|g\in G}$ for some $x\in X$ and $GX\coloneqq\set{gx|g\in
	G,x\in X}$. The set $Gx$ is sometimes called the \textit{$G$-orbit} of $x$, or just \textit{orbit} if $G$ is understood.
\end{definition}
We will mainly let $G$ be $SO(3)$ or $G_3$ and $X=\RE^3$. In this situation, the group action will
always be taken as the obvious transformation of the vectors that these group elements represent.

\subsection{Equidecomposability}

Recall that a set can be \textit{partitioned} into subsets:
\begin{definition}[Partition, decomposition]
	A \textit{partition} of a set $A$ is a collection $P$ of pairwise disjoint and nonempty subsets of $A$
	such that $A=\bigcup_{A'\in P}A'$. When the sets are written as a union in this way, it is
	also called a \textit{decomposition} of $A$.
\end{definition}
The partitions and decompositions can consist of an infinite number of sets, but we will primarily
use a finite number. Furthermore, we allow the trivial case where the decomposition is a single set: the set being
its own decomposition, or equivalently $\set{A}$ being a partition of the set $A$.

We can now define a key concept:
\begin{definition}[Equidecomposability]
	Let $G$ be a group and $X$ a $G$-set, then $A,B\subseteq X$ are
	\textit{$G$\-/equidecomposable} if both sets have a decomposition in the same finite number
	of subsets, $A=A_1\cup A_2\cup\dots\cup A_n$ and $B=B_1\cup B_2\cup\dots\cup B_n$ and there
	are $g_i\in G$ so $B_i=g_iA_i$ for all $i=1\dots n$. We can also say that $A$ and $B$ are
	\textit{equidecomposable} if $G$ is understood. This is sometimes written as $A \sim_G B$, or
	just $A \sim B$ if there is no risk of misunderstanding.
\end{definition}

From the definition, it should be clear that the identity of the group $G$ means that any subset of a
$G$-set is equidecomposable with itself, and since each $g_i$ in the definition has an inverse such that
$A_i=g^{-1}_iB_i$ where $i=1\dots n$, from which we observe:
\begin{prop}[Reflexivity and symmetry of equidecomposability]
	\label{equidecomp_inverse}
	$G$\-/equidecomposability is reflexive: $A\sim A$, and symmetric: if $A\sim B$ then $B\sim
	A$, for all subsets of a $G$-set.
\end{prop}
We will see that equidecomposability is also transitive, in \Cref{equidecomp_transitive}, which
means it is an equivalence relation, justifying the use of the symbol $\sim$.

\begin{figure}[h]
	\centering
	\def\svgscale{0.5}
	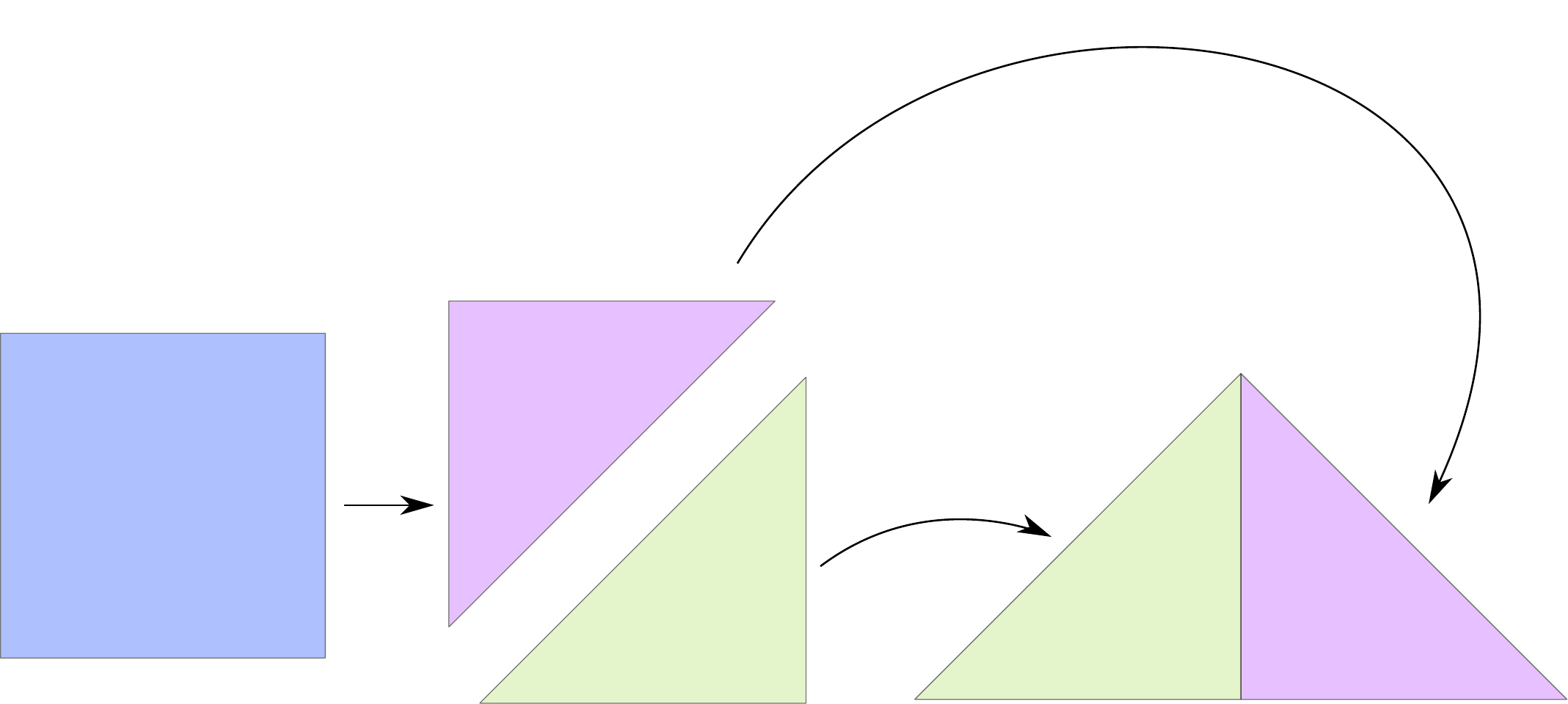
	\caption{The square $A$ can be decomposed into subsets $A_1$ and $A_2$, which are rotated
	and translated by $g_1$ and $g_2$ to form the triangle $B$. $A$ and $B$ are
	equidecomposable.}
	\label{fig:equid}
\end{figure}

\begin{example}
	The square $A$ can be decomposed into subsets $A_1$ and $A_2$, which are rotated
	and translated by $g_1$ and $g_2$ to form the triangle $B$. See \Cref{fig:equid}. Here
	$g_1,g_2\in G_2$ which is defined in the same spirit of $G_3$; the rotations and
	translations in $\RE^2$.
\end{example}

A less obvious but quite important example can be found in \Cref{ex:irrational}. This
construction gives us a simple example of an extreme case of equidecomposability:
\begin{example}
	The irrational rotation in \Cref{ex:irrational} also gives us an example of
	equidecomposability:
	let $\vect{p}$,
	$\rho$ and $X$ be as in \Cref{ex:irrational}, then the set $X$ and $X\setminus\vect{p}$ are
	$SO(2)$\-/equidecomposable, where $SO(2)$ is the group of rotations around origo in $\RE^2$,
	since $\rho X=X\setminus\vect{p}$ and $\rho\in SO(2)$. Similarly, from the opposite rotation, $X\sim_{SO(2)} X\cup{\rho^{-1}\vect{p}}$.
\end{example}
Observe that here we don not even need to decompose the original set $X$ into any subsets to get any
of the other sets. A single rotation of the entire set $X$ is enough.

We also see some immediate properties that help us apply these properties more generally.
	First note that a trivial decomposition of a set is the set itself. From this case, we see that
	an action from $G$ applied directly to a set gives a trivial form of equidecomposability,
	and due to the transitivity, two equidecomposable sets remain equidecomposable even when one
	of them is further transformed: 
\begin{prop}[Transformation invariance of equidecomposability]
	\label{equi_transf}
	If $A\sim_G B$ then $A\sim_G gB$ for any $g\in G$.
\end{prop}
For instance, two $G_3$\-/equidecomposable sets will remain $G_3$\-/equidecomposable even after one of
them gets rotated and translated. This is also true when both sets are transformed, as we will
prove transitivity for equidecomposability below in \Cref{equidecomp_transitive}). One more
property in a similar spirit can be deduced specifically for $\RE^3$ and $G_3$:

Remember that if $X\subseteq\RE^3$, the elements are vectors which allows \textit{scalar
multiplication} by elements like $k\in\RE$. Similarly to \Cref{def:groupaction}, but not to be
confused with it, we define $kX\coloneqq\set{kx|x\in X}$ and observe:
\begin{prop}[Scaling invariance of equidecomposability]
	\label{equi_scale}
	$G_3$\-/equidecomposability in $\RE^3$ is invariant under scaling, meaning that for any
	$A,B\subseteq \RE^3$ and $A\sim_{G_3}B$ it's also true that $kA\sim_{G_3}kB$ for any $k\in\RE, k>0$.
\end{prop}
\begin{proof}
	Let $A$, $B$ and $k$ be as stated. From the definition of $A\sim_{G_3} B$ we have, for some
	$n\ge1$, a decompositions $A=A_1\cup A_2 \cup \dots \cup A_n$, $B=B_1\cup B_2 \cup \dots
	\cup B_n$, and $g_i\in G_3$ for $i=1,2,\dots,n$, such that $A_i=g_iB_i$ for $i=1,2,\dots,n$.

	We see directly that we have decompositions $kA=kA_1\cup kA_2 \cup \dots kA_n$ and
	$kB=kB_1\cup kB_2 \cup \dots kB_n$, since all sets sets will still be pairwise disjoint after
	multiplication by $k$. Let $g'(x)\coloneqq kg(k^{-1}x)$. It follows that
	\[g_i'(kA_i)=kg_i(k^{-1}kA_i)=kg_i(A_i)=kB_i\text{,}\]
	so $g_i':kA_i\to kB_i$. We just need to verify that $g'\in G_3$. Remember that the
	operators $g_i\in G_{G_3}$ are defined to be of the form $g_i=\vect{t_i}+r_i(x)$, where
	$\vect{t_i}\in\RE^3$ gives the translation and $r_i\in SO(3)$ gives the rotation.

	Since $r_i$ is linear, we get:
	\[g'_i(x)=kg_i(k^{-1}x)=k(\vect{t_i}+r_i(k^{-1}x))=k\vect{t_i}+kk^{-1}r_i(x)=k\vect{t_i}+r_i(x)\]
	and since $k\vect{t_i}\in\RE^3$ it must be that all $g'_i\in G_3$ and we have our proof.
\end{proof}

\begin{remark}
	Observe that the last proposition only covers the group $G_3$ acting on $\RE^3$. It is not
	a general property for all groups and sets, but of course, it be deduced for some others
	like $SO(3)$ for $\RE^3$. And also $G_2$ and $SO(2)$ for $\RE^2$.
\end{remark}

Beyond just stating that two sets are equidecomposabile, it is often convenient to be able to
describe and work with the underlying decompositions and actions that give the
equidecomposability. An easy way of doing this is by certain piecewise functions:
\begin{definition}[Piecewise function in a group]
	\label{equidecomp_function}
	Let $A$ and $B$ be two $G$\-/equidecomposable sets. Then for some $n$ there exist
	decompositions $A=A_1\cup,A_2\cup\dots\cup A_n$ and $B=B_1\cup,B_2\cup\dots\cup B_n$ and
	$g_1,g_2,\dots,g_n\in G$ such that $A_i=g_iB_i$ for every $i=1,2,\dots,n$.
	An equivalent way of specifying the equidecomposability between two sets is through a
	piecewise-defined bijection $f:A\to B$ as:
	\[
		f(x)= \begin{cases}
			g_1x, & x\in A_1 \\
			g_2x, & x\in A_2 \\
			\vdots \\
			g_nx, & x\in A_n.
		\end{cases}
	\]
And vice-versa: if there exists a bijection that is piecewise-defined on the sets from a partition by operators from a group, then it describes the equidecomposability of the
partitioned set and its image.
	A function of this kind is sometimes said to be \textit{piecewise in G}. But note that this is
not a well-established name, even though the construction itself is.
\end{definition}
This form makes it easy to see and prove certain interesting properties of equidecomposability. Since
it is a bijection we might suspect, correctly, that the inverse function describes the original process of
equidecomposability one set to the other in reverse: both properties \Cref{equidecomp_inverse} are
easily represented in this form: reflexivity, $A\sim A$, is given by $f:A\to A$ as $f\coloneqq e$
and if $g:A\to B$ represents $A\sim B$, then $g^{-1}:B\to A$ gives the symmetry $B\sim A$.

\begin{figure}[h]
	\centering
	\def\svgscale{0.5}
	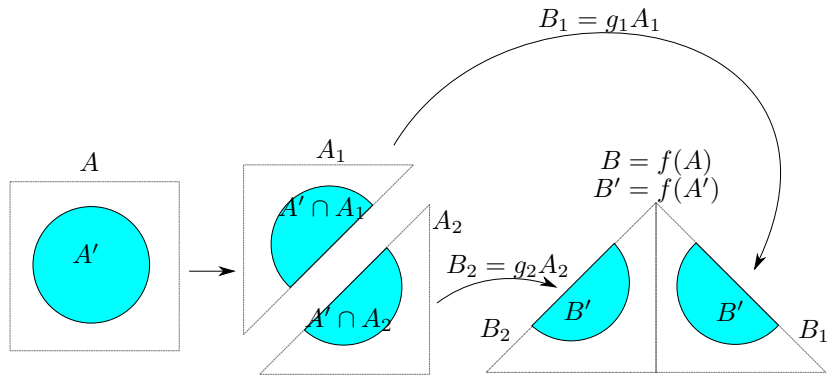
	\caption{The equidecomposability of $A$ and $B$ allows the same decomposition and
	transformations between the subset $A'\subseteq A$ and $B'\subseteq B$, so $A'\sim
	B'$.}
	\label{fig:equi_subset}
\end{figure}

Furthermore, we immediately see an important property: that the same decomposition and actions that
describes the equidecomposability between two sets also works for their subsets, as illustrated in
\Cref{fig:equi_subset}: A bijection $f:A\to B$ piecewise in a group describing $A\sim B$ can be
restricted to a subset $A'\subseteq A$ as $f':A'\to f(A)$ and it will remain a bijection, and $f'$
is defined on a partition of $A'$ instead of $A$. This partition is obtained by the intersection of
$A'$ and each set in the partition of $A$, as in \Cref{fig:equi_subset}.

\begin{prop}[Equidecomposability of subsets]
	\label{equidecomp_subsets}
	If $A\sim B$, then every a subset $A'\subset A$ is equidecomposable with a subset of $B$.
	If $f:A\to B$ describes the equidecomposability of $A$ and $B$, as in
	\Cref{equidecomp_function}, then we have $A'\sim f(A')$. Illustrated in \Cref{fig:equi_subset}.
\end{prop}

The function-based description of equidecomposability also allows us to prove that
equidecomposability is transitive.

\begin{prop}
	\label{equidecomp_transitive}
	Equidecomposability is transitive: if $A\sim_G B$ and $B\sim_G C$ for some group $G$, then $A\sim_G C$.
\end{prop}
\begin{proof}

\begin{figure}[h]
	\centering
	\begin{subfigure}{0.9\columnwidth}
		\centering
		\def\svgwidth{\columnwidth}
		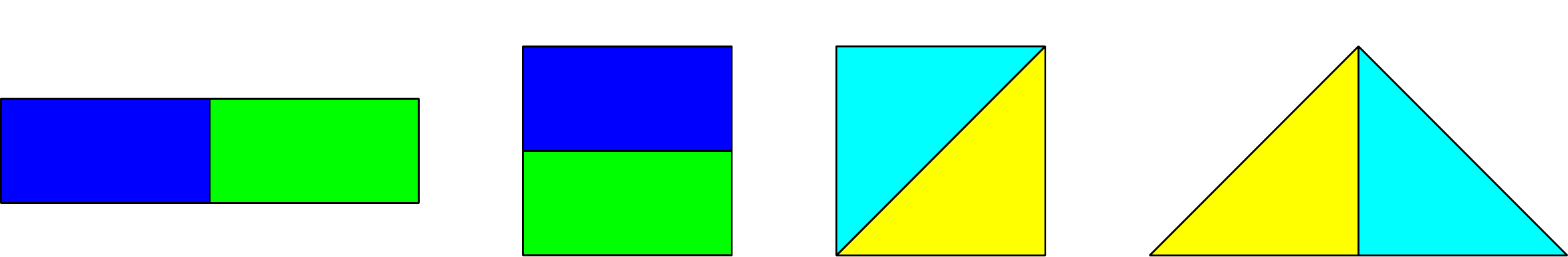
		\caption{The two steps of decomposition and transformation for $A\sim B$ and $B\sim
		C$ described by $f:A\to B$ and $g:B\to C$.}
		\label{fig:eq_rel1}
	\end{subfigure}

	\begin{subfigure}{0.9\columnwidth}
		\def\svgwidth{\columnwidth}
		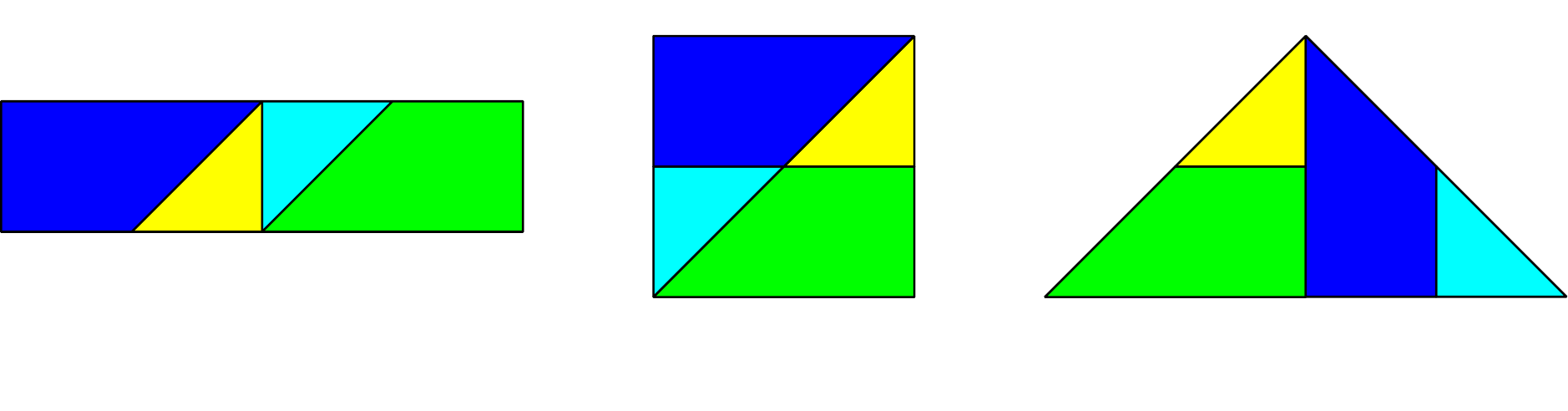
		\caption{$f$ and $g$ can be combined into $h:A\to C$ dividing $A$ into all
		necessary pieces immediately, showing that if $A\sim B$ and $B\sim C$ means $A\sim
		C$.}
		\label{fig:eq_rel2}
	\end{subfigure}
	\caption{Example of transitivity of equidecomposability.}
\end{figure}

	Let $A$, $B$ and $C$ satisfy $A\sim B$ and $B\sim C$, where the equidecomposabilities are
	described by $f:A\to B$ and $g:B\to C$, defined piecewise in the group $G$ as in \Cref{equidecomp_function}.
	Let $A=A_1\cup A_2\cup \dots\cup A_n$ be the decomposition of $A$ for $A\sim
	B$ 
	and $B=B_1\cup B_2\cup \dots\cup B_m$ the decomposition of $B$ for $B\sim C$. 
	See the illustration in \Cref{fig:eq_rel1}.

	Since $f$ and $g$ are bijections, we can composite them to obtain a bijection $h:A\to C$ by
	$h\coloneqq g\circ f$. We know that it satisfies $h(A)=g(f(A))=C$, and by transitivity of
	$G$, $h$ maps each element in $A$ to $C$ by an element from $G$. We just need to check that
	$h$ is defined on a finite partition of $A$.

	But this partition of $A$ consists of at most $n\cdot m$ pieces: for each
	$i=1,2,\dots,n$, the elements
	in $A_i$ are mapped by $f$ into $B$, where they are contained in one of $B_1,B_2,\dots,B_m$.

	For a different perspective: the new partition of $A$, as illustrated in
	\Cref{fig:eq_rel2}, can also be explicitly
	defined by intersecting each of the sets that partition $A$ with each of the preimages of
	the sets that partition of $B$. Defined these by $D_{i,j}\coloneqq A_i
	\cap f^{-1}(B'_j)$ for $i=1,2,\dots,n, j=1,2,\dots,m$.
	Some, but not all, of the $D_{i,j}$
	might be empty and can be discarded. There are at most $n\cdot m$ remaining nonempty sets,
	forming the new partition of $A$ for $h$. 
\end{proof}











The transitivity property will play an important role in several proofs. Together with
\Cref{equidecomp_inverse} we see that equidecomposability is an equivalence relation.

\begin{prop}
	For a group $G$ and a $G$-set $A$, equidecomposability is an equivalence relation on
	$\mathcal{P}(A)$, the collection of all subsets of $A$.
\end{prop}


\subsection{Paradoxicality}

We can now define the concept of a set being \textit{paradoxical}. This terminology does
not appear in the original papers by Hausdorff or Banach and Tarski. Its purpose is only to
ease repeated statements about a particular paradoxical property to make the reasoning more compact and clearer:
\begin{definition}[Paradoxical set]
	\label{def:paradoxical}
	Let $G$ be a group acting on a set $X$, and let $E\subseteq X$ be a non-empty subset. Then
	$E$ is \textit{G\-/paradoxical}, or \textit{paradoxical with respect to G}, if it has a
	\textit{paradoxical decomposition}: $E=A\cup B$ where $A\cap B=\emptyset$, $A\sim_G E$ and
	$B\sim_G E$. We can also say that $X$ is \textit{paradoxical} if the group is understood.
\end{definition}
Although the definition might appear convoluted, its form is suitable when being proved and used
inside proofs. The practical implication of it is that an object can be duplicated, like in
\Cref{quote2} on page \pageref{quote2} and the second illustration in \Cref{fig:illustration}:

\begin{example}[Duplication of paradoxical objects]
	\label{paradoxic_duplicate}
	By definition, a $G$\-/paradoxical set contains two subsets which are $G$\-/equidecomposable
	itself. If these two copies can be translated by elements of $G$ so that they are
	disjoint, which is the case for all $G_3$\-/paradoxical bounded sets in $\RE^3$, then their
	union is $G$\-/equidecomposable with the original set, as a result of transformation
	invariance of \Cref{equi_transf} and transitivity of \Cref{equidecomp_transitive}:
	\[A\sim g_1A\cup g_2A\]
	This result is not the case if the two copies share common points; taking their union would
	remove these duplicate points so parts of $A$ would be lost.
\end{example}

Observe that the description of duplicating a ball in \Cref{quote2} on page \pageref{quote2} is
inspired by the idea above, where the ball is $G_3$\-/paradoxical. That the ball in $\RE^3$ is
$G_3$\-/paradoxical is the definition of the Banach-Tarski paradox, or its weak form, and is stated
in \Cref{b-t2}. Furthermore, it is equivalent to the strong form of the Banach-Tarski paradox, in
\Cref{b-t1}. All of this will be stated and proven in the following sections.

\begin{remark}
There are other definitions of $G$\-/paradoxicality; for example, in the famous book on the
	subject, \cite{wagonbok}, the subsets need only be disjoint, with no requirements on their
	union. This turns out to be equivalent to our definition, which we will see in
	\Cref{bsb_cor}.
\end{remark}

Some further important properties follow from those of equidecomposability, the first two being
immediate results from the definition above, together with \Cref{equi_transf} and \Cref{equi_scale},
respectively, which gives us:
\begin{cor}[Transformation invariance of paradoxicality]
	\label{para_transf}
	If $A$ is $G$\-/paradoxical, then $gA$ is also $G$\-/paradoxical for all $g\in G$.
\end{cor}

\begin{cor}[Scaling invariance of paradoxicality]
	\label{para_scale}
	A $G_3$\-/paradoxical set $A\subseteq\RE^3$ is still paradoxical after being scaled: $cA$ is
	$G_3$\-/paradoxical for any $c\in\RE, c>0$.
\end{cor}
Just like in \Cref{equi_scale}, observe that this property is only stated for the group $G_3$ and
set $\RE^3$.

Much in the spirit of the reasoning behind the proof of \Cref{equidecomp_subsets}, as illustrated by \Cref{fig:equi_subset},
another powerful property can be proven: that paradoxicalness is a property inherited between
equidecomposable sets.
\begin{lemma}
	\label{paradox_inherit}
	Let $A$ and $B$ be two equidecomposable sets. Then $B$ is paradoxical if and only if $A$ is
	paradoxical.
\end{lemma}

\begin{figure}[h]
	\centering
	\begin{subfigure}{\columnwidth}
	\centering
		\def\svgwidth{\columnwidth}
		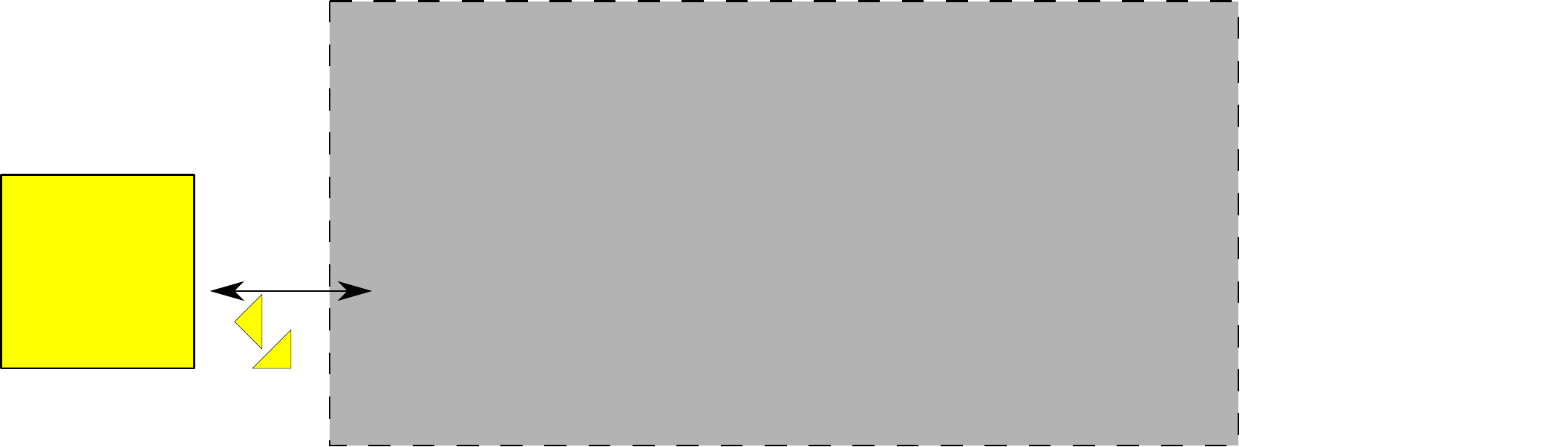
		\caption{$A$ and $B$ are equidecomposible and $B$ is paradoxical, being
		equidecomposable with its subsets $B_1$ and $B_2$. Transitivity shows that $A$ is
		equidecomposable with two copies of itself.}
		\label{fig:paradox_inherit1}
	\end{subfigure}

	\begin{subfigure}{0.47\columnwidth}
	\centering
		\def\svgwidth{\columnwidth}
		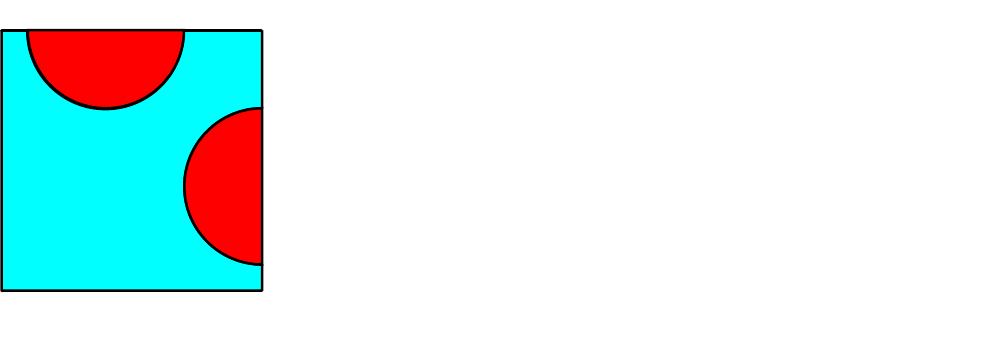
		\caption{From $A\sim B$ and $B=B_1\cup B_2$ a partition of $A$ can be made by
		creating subsets $A_1$ and $A_2$ such that $A_1\sim B_1$ and $A_2\sim B_2$.}
		\label{fig:paradox_inherit2}
	\end{subfigure}
	\quad
	\quad
	\begin{subfigure}{0.47\columnwidth}
	\centering
		\centering
		\def\svgwidth{\columnwidth}
		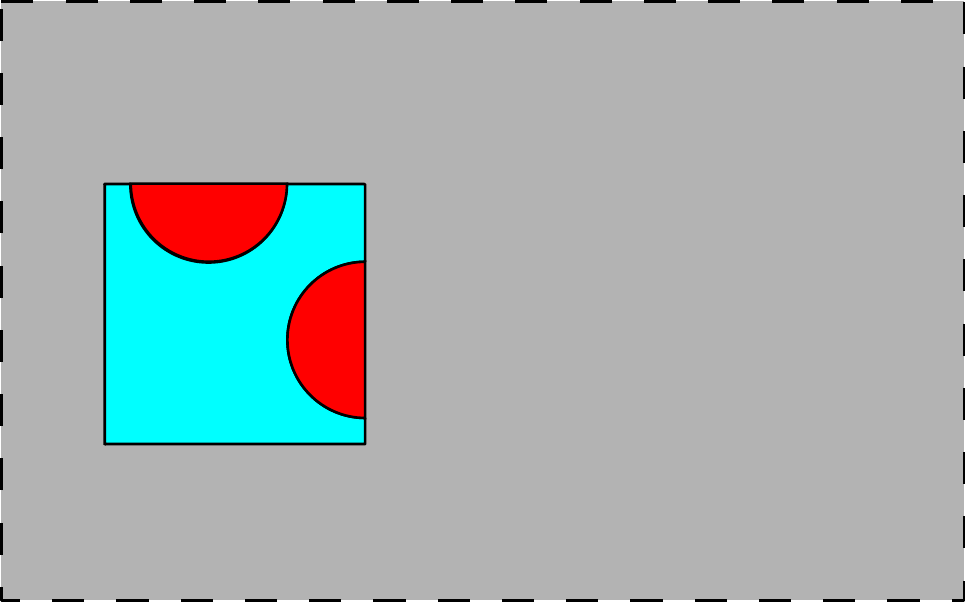
		\caption{$A$ is equidecomposable with its two subsets.}
		\label{fig:paradox_inherit3}
	\end{subfigure}
	\caption{Creation of a paradoxical decomposition of the set $A$, using the paradoxical set $B$
	it is equidecomposable with.}
	\label{fig:paradox_inherit}
\end{figure}

\begin{proof}
	Let $A$ and $B$ be as stated, with $B_1$ and $B_2$ being the paradoxical decomposition of
	$B$, that is: $A\sim B$, $B_1\sim B$ and $B_2\sim B$, where $B=B_1 \cup B_2$ and $B_1\cap
	B_2=\emptyset$. See \Cref{fig:paradox_inherit1} for an illustration. Since $A\sim
	B$, there is a bijection $f:B\to A$ describing the equidecomposability according to
	\Cref{equidecomp_function}.

	As noted in \Cref{equidecomp_subsets} we can apply this bijection to subsets of $B$ to
	obtain equidecomposable subsets in $A$: define
	$A_1\coloneqq f(B_1)$ and $A_2\coloneqq f(B_2)$ as shown in \Cref{fig:paradox_inherit2}.
	These subsets partition $A$:
	both $A=A_1 \cup A_2$ and $A_1 \cap A_2=\emptyset$ follows from the same properties of
	$B$, $B_1$ and $B_2$ and that $f$ is a bijection from $B$ to $A$.

	We now observe that the subsets satisfy $A_1 \sim B$ and $A_2\sim B$: from the construction
	$A_1 \sim B_1$ and we know that $B_1\sim B$ and $B\sim A$ from which it follows directly by
	transitivity of equidecomposability in \Cref{equidecomp_transitive} that $A_1\sim A$.
	Similarly, $A_2\sim A$, from which we see that $A$ is paradoxical.
\end{proof}

\begin{remark}
	Of course, the equidecomposability can also be shown directly by composing functions
	describing the two instances of equidecomposability: let the bijections $g_1:B\to B_1$ and $g_2:B\to B_2$
	describe the paradoxicality of $B$. Again, define $A_1\coloneqq f(B_1)$ and $A_2\coloneqq
	f(B_2)$. It follows that $f(g_1^{-1}(f^{-1}(A_1)))=f(g_1^{-1}(B_1))=f(B)=A$, thus $A \sim
	A_1$ and similarly $A \sim A_2$.
\end{remark}

We end this section by deducing a final property, which brings the machinery of the
\textit{Schröder-Bernstein theorem} from set theory to equidecomposability. Due to its close resemblance,
it has received the nickname the \textit{Banach-Schröder-Bernstein theorem}.
Banach published it in the same volume containing the article where he and Tarski
presented their famous paradox it enables, \cite{banachart}.

\begin{theorem}[The Banach-Schröder-Bernstein Theorem]
	\label{thm:b-s-b}
	Let $G$ be a group and $A,B\subseteq X$ for some $G$-set $X$. If $A$ is
	$G$\-/equidecomposable with a subset of $B$ and $B$ is $G$\-/equidecomposable with a subset of
	$A$, then $A$ and $B$ are $G$\-/equidecomposable.
\end{theorem}

This version of the proof is based on a well-known proof for the normal Schröder-Bernstein
theorem due to Birkhoff and MacLane. See \cite[pp.~29--30]{simmons1963} for more details and
\cite[pp.~419]{kursbok} for this proof in a briefer form.

\begin{proof}
	Let $A$ and $B$ be as stated, with subsets $A'\subseteq A$ and $B'\subseteq B$ such that
	$A\sim_G B'$ and $B\sim_G A'$. From \Cref{equidecomp_function}, let $f:A\to B'$
	and $g:B\to A'$ be the bijections describing the equidecomposabilities $A\sim B'$ and $B\sim A'$, respectively.


	The inverses of the bijections, $f^{-1}$ and $g^{-1}$, are defined on 
	$B'$ and $A'$, respectively, but undefined for any other elements. We can also consider $f$ and
	$g$ to be injections between $A$ and $B$, but not necessarily surjective, see
	\Cref{fig:bsb1} on the next page.

\begin{figure}[h]
	\centering
	\def\svgwidth{0.8\columnwidth}
	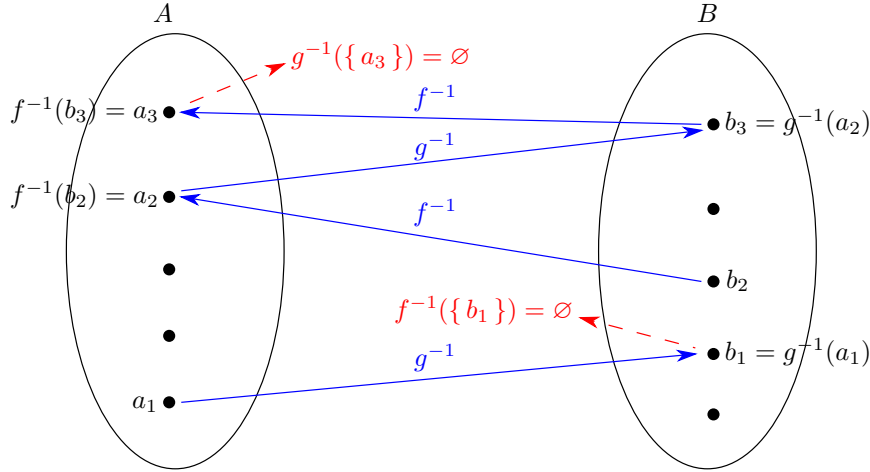
	\caption{$a_1$ and $b_3$ have one ancestor each: 
	$b_1$ in $B$ and $a_3$ in $A$, while $b_2$ has 3: its first $a_2$, second $b_3$ and
	third ancestor $a_3$ in $A$, $B$ and $A$, respectively. Since $b_1$ and $a_3$ are not in the range
	of $g$ and $f$, they have no ancestors, illustrated by the empty preimages.}
	\label{fig:b-s-b}
\end{figure}

	The first step is to partition $A$ and $B$ into three pairwise disjoint subsets each,
	between which we can find new bijections. The
	subsets will be defined using the
	concept of \textit{ancestors} of the elements, illustrated in \Cref{fig:b-s-b}: each element $a\in A$ has an ancestor
	$g^{-1}(a)$ if it is defined, that is, if $a$ is in the range of $g$. This ancestor
	is an element in $B$, which gets mapped to $a$ by $g$. The ancestors of $B$ are similarly
	defined. Since each ancestor, in turn, may have an ancestor, we can talk about the
	\textit{number of
	ancestors} for each element. The first ancestor of an element is referred to
	as its \textit{parent}. Since $f$ and $g$ are injective, each
	element has at most one parent.

	We can now partition $A$ and $B$ into three subsets, each containing elements based on how
	many ancestors they have: $A_o$, $A_e$ and $A_\infty$ containing those with an odd, even
	and infinite number of ancestors, respectively. And similarly, for $B$:
	$B_o$, $B_e$ and $B_\infty$. For the elements in $A_\infty$ and $B_\infty$, it is clear
	that we can consider $f$ as a bijection when restricted to $A_\infty\to B_\infty$, which
	describes a new equidecomposability $A_\infty\sim B_\infty$ as in \Cref{equidecomp_subsets}. For
	the remaining subsets, some more observations are needed.

	Every element in $B_o$ has a nonzero number of ancestors, so it has at least one parent; in
	$A$. For those with $1$ ancestor, their parents are precisely the elements in
	$A$ that lack parents, those having $0$ ancestors. Thus these parents are in $A_e$.
	Similarly, those of $B_o$ with $3$ ancestors have the elements in $A_e$ with $2$
	ancestors as parents. And so on, every element in $B_o$ has a parent in $A_e$. We see that $f$ maps
	some element of $A_e$ to every element in $B_o$, which means $f$ is a bijection when restricted to $A_e\to
	B_o$. Again, using \Cref{equidecomp_subsets}, we see that it describes $A_e\sim B_o$.

	Similarly, $g$ is bijective when restricted to $B_e\to A_o$, describing $A_o\sim B_e$. We can
	also describe it by $g^{-1}: A_o\to B_e$. \Cref{fig:bsb2} illustrate these subsets and
	their equidecomposabilities.

\begin{figure}[h]
	\centering
	\begin{subfigure}{0.4\textwidth}
	\centering
		\def\svgscale{1.3}
		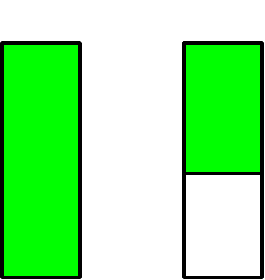
		\caption{Equidecomposability $A\sim B'$ by $f$. $f(A)=B'$ is not necessarily all of $B$.}
		\label{fig:bsb1}
	\end{subfigure}
	\quad
	\begin{subfigure}{0.4\textwidth}
	\centering
		\def\svgscale{1.3}
		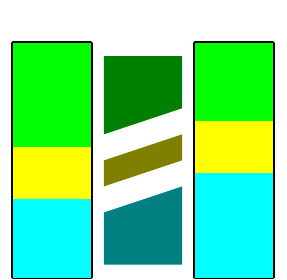
		\caption{Equidecomposabilities by $f$ and $g^{-1}$ on the subsets, together
		giving $A\sim B$.}
		\label{fig:bsb2}
	\end{subfigure}
	\caption{Construction of the new bijection, $h$.}
\end{figure}

	\Cref{fig:bsb2} also shows how these functions and subsets together define a new piecewise
	bijection $h:A\to B$, defined by:
	\[
		h(x)= \begin{cases}
			f(x), & x\in A_e \cup A_\infty \\
			g^{-1}(x), & x\in A_o\text{.}
		\end{cases}
	\]
	We can see that it is piecewise in $G$ if we consider the individual partitions of the
	sets $A_e$, $A_o$ and $A_\infty$ and the restricted bijections from $f$ and $g^{-1}$ that were
	created when we deduced $A_\infty\sim B_\infty$, $A_e\sim B_o$ and $A_o\sim B_e$: the
	bijections are still piecewise in $G$ and the partitions together gives us a new,
	finite, partition of $A$.

	By \Cref{equidecomp_function}, $h$ gives us that $A$ and $B$ are $G$\-/equidecomposable.
\end{proof}

The following corollary to \Cref{thm:b-s-b} helps make certain proofs easier:
\begin{cor}
	\label{bsb_cor}
	Let $G$ be a group acting on a set $X$. Then $A\subset X$ is $G$\-/paradoxical if it contains
	disjoint subsets $A_1$ and $A_2$, both equidecomposible with $A$.
\end{cor}
\begin{proof}
	Let $A$, $A_1$ and $A_2$ be as stated. Then $A\setminus A_1$ is equidecomposable with a
	subset of $A$, since it is a subset. And $A$ is equidecomposable with $A_2$, a
	subset of $A\setminus A_1$. Then theorem \Cref{thm:b-s-b} directly gives us that $A \sim_G
	A\setminus A_1$.

	And since $A_1$ and $A\setminus A_1$ are disjoint such that $A=A_1 \cup (A\setminus A_1)$,
	they form a paradoxical decomposition of $A$ and the
	result follows immediately from the definition of paradoxicality.
\end{proof}

This result allows us to relax the conditions in
\Cref{def:paradoxical}, by not requiring the subsets to form a partition. This can make the
process of proving that a set is paradoxical easier since we don not need to consider the union of
the subsets.

\begin{remark}
	\label{paradox_inherit_weak}
	An immediate result of \Cref{bsb_cor} is a weaker version of \Cref{paradox_inherit}: let
	$A\subseteq_G B$, where $A$ and $B$ are $G$\-/equidecomposable and $A$ is $G$\-/paradoxical, for some group $G$. Then $B$ is also
	paradoxical. Just observe that the two subsets in the paradoxical decomposition of $A$ are also subsets of $B$ and
	equidecomposable with $B$ by transitivity from \Cref{equidecomp_transitive}. The result
	now follows directly from \Cref{bsb_cor}.
\end{remark}

%% file: figs/ex_equidecomposable.pdf_tex
\begingroup%
  \makeatletter%
  \providecommand\color[2][]{%
    \errmessage{(Inkscape) Color is used for the text in Inkscape, but the package 'color.sty' is not loaded}%
    \renewcommand\color[2][]{}%
  }%
  \providecommand\transparent[1]{%
    \errmessage{(Inkscape) Transparency is used (non-zero) for the text in Inkscape, but the package 'transparent.sty' is not loaded}%
    \renewcommand\transparent[1]{}%
  }%
  \providecommand\rotatebox[2]{#2}%
  \newcommand*\fsize{\dimexpr\f@size pt\relax}%
  \newcommand*\lineheight[1]{\fontsize{\fsize}{#1\fsize}\selectfont}%
  \ifx\svgwidth\undefined%
    \setlength{\unitlength}{615.22311191bp}%
    \ifx\svgscale\undefined%
      \relax%
    \else%
      \setlength{\unitlength}{\unitlength * \real{\svgscale}}%
    \fi%
  \else%
    \setlength{\unitlength}{\svgwidth}%
  \fi%
  \global\let\svgwidth\undefined%
  \global\let\svgscale\undefined%
  \makeatother%
  \begin{picture}(1,0.44921)%
    \lineheight{1}%
    \setlength\tabcolsep{0pt}%
    \put(0,0){\includegraphics[width=\unitlength,page=1]{ex_equidecomposable.pdf}}%
    \put(0.09734775,0.12837537){\color[rgb]{0,0,0}\makebox(0,0)[t]{\smash{\begin{tabular}[t]{c}$A$\end{tabular}}}}%
    \put(0.7881893,0.22689109){\color[rgb]{0,0,0}\makebox(0,0)[t]{\smash{\begin{tabular}[t]{c}$B=B_1\cup B_2$\end{tabular}}}}%
    \put(0.34486851,0.19364203){\color[rgb]{0,0,0}\makebox(0,0)[t]{\smash{\begin{tabular}[t]{c}$A_1$\end{tabular}}}}%
    \put(0.44707862,0.04340554){\color[rgb]{0,0,0}\makebox(0,0)[t]{\smash{\begin{tabular}[t]{c}$A_2$\end{tabular}}}}%
    \put(0.72415403,0.05079423){\color[rgb]{0,0,0}\makebox(0,0)[t]{\smash{\begin{tabular}[t]{c}$B_2$\end{tabular}}}}%
    \put(0.86084466,0.05079423){\color[rgb]{0,0,0}\makebox(0,0)[t]{\smash{\begin{tabular}[t]{c}$B_1$\end{tabular}}}}%
    \put(0.72661688,0.43131122){\color[rgb]{0,0,0}\makebox(0,0)[t]{\smash{\begin{tabular}[t]{c}$B_1=g_1A_1$\end{tabular}}}}%
    \put(0.61209243,0.1296068){\color[rgb]{0,0,0}\makebox(0,0)[t]{\smash{\begin{tabular}[t]{c}$B_2=g_2A_2$\end{tabular}}}}%
  \end{picture}%
\endgroup%

%% file: figs/ex_equidecomposable_subset.pdf_tex
\begingroup%
  \makeatletter%
  \providecommand\color[2][]{%
    \errmessage{(Inkscape) Color is used for the text in Inkscape, but the package 'color.sty' is not loaded}%
    \renewcommand\color[2][]{}%
  }%
  \providecommand\transparent[1]{%
    \errmessage{(Inkscape) Transparency is used (non-zero) for the text in Inkscape, but the package 'transparent.sty' is not loaded}%
    \renewcommand\transparent[1]{}%
  }%
  \providecommand\rotatebox[2]{#2}%
  \newcommand*\fsize{\dimexpr\f@size pt\relax}%
  \newcommand*\lineheight[1]{\fontsize{\fsize}{#1\fsize}\selectfont}%
  \ifx\svgwidth\undefined%
    \setlength{\unitlength}{615.29130379bp}%
    \ifx\svgscale\undefined%
      \relax%
    \else%
      \setlength{\unitlength}{\unitlength * \real{\svgscale}}%
    \fi%
  \else%
    \setlength{\unitlength}{\svgwidth}%
  \fi%
  \global\let\svgwidth\undefined%
  \global\let\svgscale\undefined%
  \makeatother%
  \begin{picture}(1,0.43435582)%
    \lineheight{1}%
    \setlength\tabcolsep{0pt}%
    \put(0,0){\includegraphics[width=\unitlength,page=1]{ex_equidecomposable_subset.pdf}}%
    \put(0.09346551,0.1362858){\color[rgb]{0,0,0}\makebox(0,0)[t]{\smash{\begin{tabular}[t]{c}$A'$\end{tabular}}}}%
    \put(0.38289843,0.19495556){\color[rgb]{0,0,0}\makebox(0,0)[t]{\smash{\begin{tabular}[t]{c}$A' \cap A_1$\end{tabular}}}}%
    \put(0.41517191,0.05875573){\color[rgb]{0,0,0}\makebox(0,0)[t]{\smash{\begin{tabular}[t]{c}$A' \cap A_2$\end{tabular}}}}%
    \put(0.69583633,0.0672118){\color[rgb]{0,0,0}\makebox(0,0)[t]{\smash{\begin{tabular}[t]{c}$B'$\end{tabular}}}}%
    \put(0.882301,0.06576923){\color[rgb]{0,0,0}\makebox(0,0)[t]{\smash{\begin{tabular}[t]{c}$B'$\end{tabular}}}}%
    \put(0.09738571,0.24921759){\color[rgb]{0,0,0}\makebox(0,0)[t]{\smash{\begin{tabular}[t]{c}$A$\end{tabular}}}}%
    \put(0.39425027,0.26601529){\color[rgb]{0,0,0}\makebox(0,0)[t]{\smash{\begin{tabular}[t]{c}$A_1$\end{tabular}}}}%
    \put(0.53650412,0.1782091){\color[rgb]{0,0,0}\makebox(0,0)[t]{\smash{\begin{tabular}[t]{c}$A_2$\end{tabular}}}}%
    \put(0.79219594,0.25142299){\color[rgb]{0,0,0}\makebox(0,0)[t]{\smash{\begin{tabular}[t]{c}$B=f(A)$\\$B'=f(A')$\end{tabular}}}}%
    \put(0.72284184,0.4263375){\color[rgb]{0,0,0}\makebox(0,0)[t]{\smash{\begin{tabular}[t]{c}$B_1=g_1A_1$\end{tabular}}}}%
    \put(0.61017686,0.12405088){\color[rgb]{0,0,0}\makebox(0,0)[t]{\smash{\begin{tabular}[t]{c}$B_2=g_2A_2$\end{tabular}}}}%
    \put(0.59540123,0.04595315){\color[rgb]{0,0,0}\makebox(0,0)[t]{\smash{\begin{tabular}[t]{c}$B_2$\end{tabular}}}}%
    \put(0.98511088,0.04533747){\color[rgb]{0,0,0}\makebox(0,0)[t]{\smash{\begin{tabular}[t]{c}$B_1$\end{tabular}}}}%
  \end{picture}%
\endgroup%

%% file: figs/eq_relation.pdf_tex
\begingroup%
  \makeatletter%
  \providecommand\color[2][]{%
    \errmessage{(Inkscape) Color is used for the text in Inkscape, but the package 'color.sty' is not loaded}%
    \renewcommand\color[2][]{}%
  }%
  \providecommand\transparent[1]{%
    \errmessage{(Inkscape) Transparency is used (non-zero) for the text in Inkscape, but the package 'transparent.sty' is not loaded}%
    \renewcommand\transparent[1]{}%
  }%
  \providecommand\rotatebox[2]{#2}%
  \newcommand*\fsize{\dimexpr\f@size pt\relax}%
  \newcommand*\lineheight[1]{\fontsize{\fsize}{#1\fsize}\selectfont}%
  \ifx\svgwidth\undefined%
    \setlength{\unitlength}{675.64014754bp}%
    \ifx\svgscale\undefined%
      \relax%
    \else%
      \setlength{\unitlength}{\unitlength * \real{\svgscale}}%
    \fi%
  \else%
    \setlength{\unitlength}{\svgwidth}%
  \fi%
  \global\let\svgwidth\undefined%
  \global\let\svgscale\undefined%
  \makeatother%
  \begin{picture}(1,0.16326531)%
    \lineheight{1}%
    \setlength\tabcolsep{0pt}%
    \put(0,0){\includegraphics[width=\unitlength,page=1]{eq_relation.pdf}}%
    \put(0.13376203,0.12266145){\color[rgb]{0,0,0}\makebox(0,0)[t]{\smash{\begin{tabular}[t]{c}$A$\end{tabular}}}}%
    \put(0.40017603,0.1559632){\color[rgb]{0,0,0}\makebox(0,0)[t]{\smash{\begin{tabular}[t]{c}$B$\end{tabular}}}}%
    \put(0.59998651,0.1559632){\color[rgb]{0,0,0}\makebox(0,0)[t]{\smash{\begin{tabular}[t]{c}$B$\end{tabular}}}}%
    \put(0.86640058,0.1559632){\color[rgb]{0,0,0}\makebox(0,0)[t]{\smash{\begin{tabular}[t]{c}$C$\end{tabular}}}}%
    \put(0.06715853,0.05605795){\color[rgb]{0,0,0}\makebox(0,0)[t]{\smash{\begin{tabular}[t]{c}$A_1$\end{tabular}}}}%
    \put(0.20036553,0.05605795){\color[rgb]{0,0,0}\makebox(0,0)[t]{\smash{\begin{tabular}[t]{c}$A_2$\end{tabular}}}}%
    \put(0.56668476,0.08935969){\color[rgb]{0,0,0}\makebox(0,0)[t]{\smash{\begin{tabular}[t]{c}$B_1$\end{tabular}}}}%
    \put(0.62218776,0.03385678){\color[rgb]{0,0,0}\makebox(0,0)[t]{\smash{\begin{tabular}[t]{c}$B_2$\end{tabular}}}}%
    \put(0,0){\includegraphics[width=\unitlength,page=2]{eq_relation.pdf}}%
    \put(0.3002708,0.07825911){\color[rgb]{0,0,0}\makebox(0,0)[t]{\smash{\begin{tabular}[t]{c}$f$\end{tabular}}}}%
    \put(0.72209295,0.07825911){\color[rgb]{0,0,0}\makebox(0,0)[t]{\smash{\begin{tabular}[t]{c}$g$\end{tabular}}}}%
  \end{picture}%
\endgroup%

%% file: figs/eq_relation2.pdf_tex
\begingroup%
  \makeatletter%
  \providecommand\color[2][]{%
    \errmessage{(Inkscape) Color is used for the text in Inkscape, but the package 'color.sty' is not loaded}%
    \renewcommand\color[2][]{}%
  }%
  \providecommand\transparent[1]{%
    \errmessage{(Inkscape) Transparency is used (non-zero) for the text in Inkscape, but the package 'transparent.sty' is not loaded}%
    \renewcommand\transparent[1]{}%
  }%
  \providecommand\rotatebox[2]{#2}%
  \newcommand*\fsize{\dimexpr\f@size pt\relax}%
  \newcommand*\lineheight[1]{\fontsize{\fsize}{#1\fsize}\selectfont}%
  \ifx\svgwidth\undefined%
    \setlength{\unitlength}{540.64016989bp}%
    \ifx\svgscale\undefined%
      \relax%
    \else%
      \setlength{\unitlength}{\unitlength * \real{\svgscale}}%
    \fi%
  \else%
    \setlength{\unitlength}{\svgwidth}%
  \fi%
  \global\let\svgwidth\undefined%
  \global\let\svgscale\undefined%
  \makeatother%
  \begin{picture}(1,0.26183156)%
    \lineheight{1}%
    \setlength\tabcolsep{0pt}%
    \put(0,0){\includegraphics[width=\unitlength,page=1]{eq_relation2.pdf}}%
    \put(0.0561834,0.15559899){\color[rgb]{0,0,0}\makebox(0,0)[t]{\smash{\begin{tabular}[t]{c}$D_{1,1}$\end{tabular}}}}%
    \put(0.13248183,0.12091787){\color[rgb]{0,0,0}\makebox(0,0)[t]{\smash{\begin{tabular}[t]{c}$D_{1,2}$\end{tabular}}}}%
    \put(0.19490782,0.17363316){\color[rgb]{0,0,0}\makebox(0,0)[t]{\smash{\begin{tabular}[t]{c}$D_{2,1}$\end{tabular}}}}%
    \put(0.26427003,0.14172653){\color[rgb]{0,0,0}\makebox(0,0)[t]{\smash{\begin{tabular}[t]{c}$D_{2,2}$\end{tabular}}}}%
    \put(0,0){\includegraphics[width=\unitlength,page=2]{eq_relation2.pdf}}%
    \put(0.50010155,0.00300209){\color[rgb]{0,0,0}\makebox(0,0)[t]{\smash{\begin{tabular}[t]{c}$h=g\circ f$\end{tabular}}}}%
    \put(0.66657083,0.16947144){\color[rgb]{0,0,0}\makebox(0,0)[t]{\smash{\begin{tabular}[t]{c}$g$\end{tabular}}}}%
    \put(0.37524956,0.16947144){\color[rgb]{0,0,0}\makebox(0,0)[t]{\smash{\begin{tabular}[t]{c}$f$\end{tabular}}}}%
    \put(0.16716293,0.21108875){\color[rgb]{0,0,0}\makebox(0,0)[t]{\smash{\begin{tabular}[t]{c}$A$\end{tabular}}}}%
    \put(0.50010155,0.25270608){\color[rgb]{0,0,0}\makebox(0,0)[t]{\smash{\begin{tabular}[t]{c}$B$\end{tabular}}}}%
    \put(0.83304023,0.25270608){\color[rgb]{0,0,0}\makebox(0,0)[t]{\smash{\begin{tabular}[t]{c}$C$\end{tabular}}}}%
  \end{picture}%
\endgroup%

%% file: figs/paradox_inherit.pdf_tex
\begingroup%
  \makeatletter%
  \providecommand\color[2][]{%
    \errmessage{(Inkscape) Color is used for the text in Inkscape, but the package 'color.sty' is not loaded}%
    \renewcommand\color[2][]{}%
  }%
  \providecommand\transparent[1]{%
    \errmessage{(Inkscape) Transparency is used (non-zero) for the text in Inkscape, but the package 'transparent.sty' is not loaded}%
    \renewcommand\transparent[1]{}%
  }%
  \providecommand\rotatebox[2]{#2}%
  \newcommand*\fsize{\dimexpr\f@size pt\relax}%
  \newcommand*\lineheight[1]{\fontsize{\fsize}{#1\fsize}\selectfont}%
  \ifx\svgwidth\undefined%
    \setlength{\unitlength}{608.25113009bp}%
    \ifx\svgscale\undefined%
      \relax%
    \else%
      \setlength{\unitlength}{\unitlength * \real{\svgscale}}%
    \fi%
  \else%
    \setlength{\unitlength}{\svgwidth}%
  \fi%
  \global\let\svgwidth\undefined%
  \global\let\svgscale\undefined%
  \makeatother%
  \begin{picture}(1,0.28483495)%
    \lineheight{1}%
    \setlength\tabcolsep{0pt}%
    \put(0,0){\includegraphics[width=\unitlength,page=1]{paradox_inherit.pdf}}%
    \put(0.06226966,0.10542616){\color[rgb]{0,0,0}\makebox(0,0)[t]{\smash{\begin{tabular}[t]{c}$A$\end{tabular}}}}%
    \put(0.16707833,0.11159139){\color[rgb]{0,0,0}\makebox(0,0)[t]{\smash{\begin{tabular}[t]{c}$A\sim B$\end{tabular}}}}%
    \put(0,0){\includegraphics[width=\unitlength,page=2]{paradox_inherit.pdf}}%
    \put(0.3456913,0.09908018){\color[rgb]{0,0,0}\makebox(0,0)[t]{\smash{\begin{tabular}[t]{c}$B_1$\end{tabular}}}}%
    \put(0.27170879,0.06208888){\color[rgb]{0,0,0}\makebox(0,0)[t]{\smash{\begin{tabular}[t]{c}$B_2$\end{tabular}}}}%
    \put(0.34605496,0.02509749){\color[rgb]{0,0,0}\makebox(0,0)[t]{\smash{\begin{tabular}[t]{c}$B=B_1 \cup B_2$\end{tabular}}}}%
    \put(0,0){\includegraphics[width=\unitlength,page=3]{paradox_inherit.pdf}}%
    \put(0.49999991,0.1977236){\color[rgb]{0,0,0}\makebox(0,0)[t]{\smash{\begin{tabular}[t]{c}$B_1\sim B$\end{tabular}}}}%
    \put(0.49999994,0.09908012){\color[rgb]{0,0,0}\makebox(0,0)[t]{\smash{\begin{tabular}[t]{c}$B_2\sim B$\end{tabular}}}}%
    \put(0,0){\includegraphics[width=\unitlength,page=4]{paradox_inherit.pdf}}%
    \put(0.83292171,0.23489573){\color[rgb]{0,0,0}\makebox(0,0)[t]{\smash{\begin{tabular}[t]{c}$A\sim B$\end{tabular}}}}%
    \put(0,0){\includegraphics[width=\unitlength,page=5]{paradox_inherit.pdf}}%
    \put(0.83292171,0.07460009){\color[rgb]{0,0,0}\makebox(0,0)[t]{\smash{\begin{tabular}[t]{c}$A\sim B$\end{tabular}}}}%
    \put(0.33353909,0.25339118){\color[rgb]{0,0,0}\makebox(0,0)[t]{\smash{\begin{tabular}[t]{c}Paradoxicality of $B$\end{tabular}}}}%
    \put(0.65413046,0.18868087){\color[rgb]{0,0,0}\makebox(0,0)[t]{\smash{\begin{tabular}[t]{c}$B$\end{tabular}}}}%
    \put(0.65413046,0.04688089){\color[rgb]{0,0,0}\makebox(0,0)[t]{\smash{\begin{tabular}[t]{c}$B$\end{tabular}}}}%
    \put(0.93773037,0.21334174){\color[rgb]{0,0,0}\makebox(0,0)[t]{\smash{\begin{tabular}[t]{c}$A$\end{tabular}}}}%
    \put(0.93773037,0.05304612){\color[rgb]{0,0,0}\makebox(0,0)[t]{\smash{\begin{tabular}[t]{c}$A$\end{tabular}}}}%
  \end{picture}%
\endgroup%

%% file: figs/paradox_inherit2.pdf_tex
\begingroup%
  \makeatletter%
  \providecommand\color[2][]{%
    \errmessage{(Inkscape) Color is used for the text in Inkscape, but the package 'color.sty' is not loaded}%
    \renewcommand\color[2][]{}%
  }%
  \providecommand\transparent[1]{%
    \errmessage{(Inkscape) Transparency is used (non-zero) for the text in Inkscape, but the package 'transparent.sty' is not loaded}%
    \renewcommand\transparent[1]{}%
  }%
  \providecommand\rotatebox[2]{#2}%
  \newcommand*\fsize{\dimexpr\f@size pt\relax}%
  \newcommand*\lineheight[1]{\fontsize{\fsize}{#1\fsize}\selectfont}%
  \ifx\svgwidth\undefined%
    \setlength{\unitlength}{285.85942957bp}%
    \ifx\svgscale\undefined%
      \relax%
    \else%
      \setlength{\unitlength}{\unitlength * \real{\svgscale}}%
    \fi%
  \else%
    \setlength{\unitlength}{\svgwidth}%
  \fi%
  \global\let\svgwidth\undefined%
  \global\let\svgscale\undefined%
  \makeatother%
  \begin{picture}(1,0.35147925)%
    \lineheight{1}%
    \setlength\tabcolsep{0pt}%
    \put(0,0){\includegraphics[width=\unitlength,page=1]{paradox_inherit2.pdf}}%
    \put(0.10664492,0.26862878){\color[rgb]{0,0,0}\makebox(0,0)[t]{\smash{\begin{tabular}[t]{c}$A_1$\end{tabular}}}}%
    \put(0.23782825,0.15056402){\color[rgb]{0,0,0}\makebox(0,0)[t]{\smash{\begin{tabular}[t]{c}$A_1$\end{tabular}}}}%
    \put(0.10664492,0.13744567){\color[rgb]{0,0,0}\makebox(0,0)[t]{\smash{\begin{tabular}[t]{c}$A_2$\end{tabular}}}}%
    \put(0.13288156,0.00626196){\color[rgb]{0,0,0}\makebox(0,0)[t]{\smash{\begin{tabular}[t]{c}$A=A_1 \cup A_2$\end{tabular}}}}%
    \put(0,0){\includegraphics[width=\unitlength,page=2]{paradox_inherit2.pdf}}%
    \put(0.73709866,0.15056402){\color[rgb]{0,0,0}\makebox(0,0)[t]{\smash{\begin{tabular}[t]{c}$B_1$\end{tabular}}}}%
    \put(0.57890495,0.08497228){\color[rgb]{0,0,0}\makebox(0,0)[t]{\smash{\begin{tabular}[t]{c}$B_2$\end{tabular}}}}%
    \put(0,0){\includegraphics[width=\unitlength,page=3]{paradox_inherit2.pdf}}%
    \put(0.73709873,0.00626196){\color[rgb]{0,0,0}\makebox(0,0)[t]{\smash{\begin{tabular}[t]{c}$B=B_1 \cup B_2$\end{tabular}}}}%
    \put(0.42225879,0.33422042){\color[rgb]{0,0,0}\makebox(0,0)[t]{\smash{\begin{tabular}[t]{c}$A\sim B:$\\~$\begin{cases}A_1\sim B_1\\A_2\sim B_2\end{cases}$\end{tabular}}}}%
  \end{picture}%
\endgroup%

%% file: figs/paradox_inherit3.pdf_tex
\begingroup%
  \makeatletter%
  \providecommand\color[2][]{%
    \errmessage{(Inkscape) Color is used for the text in Inkscape, but the package 'color.sty' is not loaded}%
    \renewcommand\color[2][]{}%
  }%
  \providecommand\transparent[1]{%
    \errmessage{(Inkscape) Transparency is used (non-zero) for the text in Inkscape, but the package 'transparent.sty' is not loaded}%
    \renewcommand\transparent[1]{}%
  }%
  \providecommand\rotatebox[2]{#2}%
  \newcommand*\fsize{\dimexpr\f@size pt\relax}%
  \newcommand*\lineheight[1]{\fontsize{\fsize}{#1\fsize}\selectfont}%
  \ifx\svgwidth\undefined%
    \setlength{\unitlength}{278.0498337bp}%
    \ifx\svgscale\undefined%
      \relax%
    \else%
      \setlength{\unitlength}{\unitlength * \real{\svgscale}}%
    \fi%
  \else%
    \setlength{\unitlength}{\svgwidth}%
  \fi%
  \global\let\svgwidth\undefined%
  \global\let\svgscale\undefined%
  \makeatother%
  \begin{picture}(1,0.62316533)%
    \lineheight{1}%
    \setlength\tabcolsep{0pt}%
    \put(0,0){\includegraphics[width=\unitlength,page=1]{paradox_inherit3.pdf}}%
    \put(0.21637966,0.37901709){\color[rgb]{0,0,0}\makebox(0,0)[t]{\smash{\begin{tabular}[t]{c}$A_1$\end{tabular}}}}%
    \put(0.35124749,0.25763601){\color[rgb]{0,0,0}\makebox(0,0)[t]{\smash{\begin{tabular}[t]{c}$A_1$\end{tabular}}}}%
    \put(0.21637972,0.24414921){\color[rgb]{0,0,0}\makebox(0,0)[t]{\smash{\begin{tabular}[t]{c}$A_2$\end{tabular}}}}%
    \put(0.24335319,0.10928117){\color[rgb]{0,0,0}\makebox(0,0)[t]{\smash{\begin{tabular}[t]{c}$A=A_1 \cup A_2$\end{tabular}}}}%
    \put(0,0){\includegraphics[width=\unitlength,page=2]{paradox_inherit3.pdf}}%
    \put(0.540858,0.43256819){\color[rgb]{0,0,0}\makebox(0,0)[t]{\smash{\begin{tabular}[t]{c}$A\sim A_1$\end{tabular}}}}%
    \put(0.52737135,0.2572398){\color[rgb]{0,0,0}\makebox(0,0)[t]{\smash{\begin{tabular}[t]{c}$A\sim A_2$\end{tabular}}}}%
    \put(0,0){\includegraphics[width=\unitlength,page=3]{paradox_inherit3.pdf}}%
    \put(0.83677204,0.14974119){\color[rgb]{0,0,0}\makebox(0,0)[t]{\smash{\begin{tabular}[t]{c}$A$\end{tabular}}}}%
    \put(0.25684005,0.55434489){\color[rgb]{0,0,0}\makebox(0,0)[t]{\smash{\begin{tabular}[t]{c}Paradoxicality of $A$\end{tabular}}}}%
    \put(0.24335325,0.45993735){\color[rgb]{0,0,0}\makebox(0,0)[t]{\smash{\begin{tabular}[t]{c}$A$\end{tabular}}}}%
    \put(0.83677204,0.44645057){\color[rgb]{0,0,0}\makebox(0,0)[t]{\smash{\begin{tabular}[t]{c}$A$\end{tabular}}}}%
  \end{picture}%
\endgroup%

%% file: figs/banach-schroder-bernstein.pdf_tex
\begingroup%
  \makeatletter%
  \providecommand\color[2][]{%
    \errmessage{(Inkscape) Color is used for the text in Inkscape, but the package 'color.sty' is not loaded}%
    \renewcommand\color[2][]{}%
  }%
  \providecommand\transparent[1]{%
    \errmessage{(Inkscape) Transparency is used (non-zero) for the text in Inkscape, but the package 'transparent.sty' is not loaded}%
    \renewcommand\transparent[1]{}%
  }%
  \providecommand\rotatebox[2]{#2}%
  \newcommand*\fsize{\dimexpr\f@size pt\relax}%
  \newcommand*\lineheight[1]{\fontsize{\fsize}{#1\fsize}\selectfont}%
  \ifx\svgwidth\undefined%
    \setlength{\unitlength}{502.32282263bp}%
    \ifx\svgscale\undefined%
      \relax%
    \else%
      \setlength{\unitlength}{\unitlength * \real{\svgscale}}%
    \fi%
  \else%
    \setlength{\unitlength}{\svgwidth}%
  \fi%
  \global\let\svgwidth\undefined%
  \global\let\svgscale\undefined%
  \makeatother%
  \begin{picture}(1,0.57812789)%
    \lineheight{1}%
    \setlength\tabcolsep{0pt}%
    \put(0,0){\includegraphics[width=\unitlength,page=1]{banach-schroder-bernstein.pdf}}%
    \put(0.49169057,0.30682581){\color[rgb]{0,0,0.99607843}\makebox(0,0)[t]{\smash{\begin{tabular}[t]{c}$f^{-1}$\end{tabular}}}}%
    \put(0.49169048,0.12765805){\color[rgb]{0,0,0.99607843}\makebox(0,0)[t]{\smash{\begin{tabular}[t]{c}$g^{-1}$\end{tabular}}}}%
    \put(0,0){\includegraphics[width=\unitlength,page=2]{banach-schroder-bernstein.pdf}}%
    \put(0.15575116,0.5531813){\color[rgb]{0,0,0}\makebox(0,0)[t]{\smash{\begin{tabular}[t]{c}$A$\end{tabular}}}}%
    \put(0.82762987,0.5531813){\color[rgb]{0,0,0}\makebox(0,0)[t]{\smash{\begin{tabular}[t]{c}$B$\end{tabular}}}}%
    \put(0.93960957,0.13512345){\color[rgb]{0,0,0}\makebox(0,0)[t]{\smash{\begin{tabular}[t]{c}$b_1=g^{-1}(a_1)$\end{tabular}}}}%
    \put(0.13335518,0.07540088){\color[rgb]{0,0,0}\makebox(0,0)[t]{\smash{\begin{tabular}[t]{c}$a_1$\end{tabular}}}}%
    \put(0.42450257,0.50092404){\color[rgb]{1,0,0}\makebox(0,0)[t]{\smash{\begin{tabular}[t]{c}$g^{-1}(\set{a_3})=\varnothing$\end{tabular}}}}%
    \put(0,0){\includegraphics[width=\unitlength,page=3]{banach-schroder-bernstein.pdf}}%
    \put(0.49169053,0.44866684){\color[rgb]{0,0,0.99607843}\makebox(0,0)[t]{\smash{\begin{tabular}[t]{c}$f^{-1}$\end{tabular}}}}%
    \put(0,0){\includegraphics[width=\unitlength,page=4]{banach-schroder-bernstein.pdf}}%
    \put(0.49169053,0.38894429){\color[rgb]{0,0,1}\makebox(0,0)[t]{\smash{\begin{tabular}[t]{c}$g^{-1}$\end{tabular}}}}%
    \put(0.86495644,0.22470727){\color[rgb]{0,0,0}\makebox(0,0)[t]{\smash{\begin{tabular}[t]{c}$b_2$\end{tabular}}}}%
    \put(0.05870202,0.43373618){\color[rgb]{0,0,0}\makebox(0,0)[t]{\smash{\begin{tabular}[t]{c}$f^{-1}(b_3)=a_3$\end{tabular}}}}%
    \put(0.93960957,0.41880557){\color[rgb]{0,0,0}\makebox(0,0)[t]{\smash{\begin{tabular}[t]{c}$b_3=g^{-1}(a_2)$\end{tabular}}}}%
    \put(0.05870199,0.32922175){\color[rgb]{0,0,0}\makebox(0,0)[t]{\smash{\begin{tabular}[t]{c}$f^{-1}(b_2)=a_2$\end{tabular}}}}%
    \put(0.55141301,0.18738062){\color[rgb]{0.99607843,0,0}\makebox(0,0)[t]{\smash{\begin{tabular}[t]{c}$f^{-1}(\set{b_1})=\emptyset$\end{tabular}}}}%
    \put(0,0){\includegraphics[width=\unitlength,page=5]{banach-schroder-bernstein.pdf}}%
  \end{picture}%
\endgroup%

%% file: figs/banach-schroder-bernstein_function1.pdf_tex
\begingroup%
  \makeatletter%
  \providecommand\color[2][]{%
    \errmessage{(Inkscape) Color is used for the text in Inkscape, but the package 'color.sty' is not loaded}%
    \renewcommand\color[2][]{}%
  }%
  \providecommand\transparent[1]{%
    \errmessage{(Inkscape) Transparency is used (non-zero) for the text in Inkscape, but the package 'transparent.sty' is not loaded}%
    \renewcommand\transparent[1]{}%
  }%
  \providecommand\rotatebox[2]{#2}%
  \newcommand*\fsize{\dimexpr\f@size pt\relax}%
  \newcommand*\lineheight[1]{\fontsize{\fsize}{#1\fsize}\selectfont}%
  \ifx\svgwidth\undefined%
    \setlength{\unitlength}{76.03464758bp}%
    \ifx\svgscale\undefined%
      \relax%
    \else%
      \setlength{\unitlength}{\unitlength * \real{\svgscale}}%
    \fi%
  \else%
    \setlength{\unitlength}{\svgwidth}%
  \fi%
  \global\let\svgwidth\undefined%
  \global\let\svgscale\undefined%
  \makeatother%
  \begin{picture}(1,1.05808224)%
    \lineheight{1}%
    \setlength\tabcolsep{0pt}%
    \put(0,0){\includegraphics[width=\unitlength,page=1]{banach-schroder-bernstein_function1.pdf}}%
    \put(0.10544303,0.99319603){\color[rgb]{0,0,0}\makebox(0,0)[t]{\smash{\begin{tabular}[t]{c}$A$\end{tabular}}}}%
    \put(0.79591774,0.99319603){\color[rgb]{0,0,0}\makebox(0,0)[t]{\smash{\begin{tabular}[t]{c}$B$\end{tabular}}}}%
    \put(0,0){\includegraphics[width=\unitlength,page=2]{banach-schroder-bernstein_function1.pdf}}%
    \put(0.49999995,0.49999963){\color[rgb]{0,0,0}\makebox(0,0)[t]{\smash{\begin{tabular}[t]{c}$f$\end{tabular}}}}%
    \put(0.84523741,0.64795849){\color[rgb]{0,0,0}\makebox(0,0)[t]{\smash{\begin{tabular}[t]{c}$f(A)$\end{tabular}}}}%
    \put(0.50000002,0.89455697){\color[rgb]{0,0,0}\makebox(0,0)[t]{\smash{\begin{tabular}[t]{c}$\sim$\end{tabular}}}}%
  \end{picture}%
\endgroup%

%% file: figs/banach-schroder-bernstein_function2.pdf_tex
\begingroup%
  \makeatletter%
  \providecommand\color[2][]{%
    \errmessage{(Inkscape) Color is used for the text in Inkscape, but the package 'color.sty' is not loaded}%
    \renewcommand\color[2][]{}%
  }%
  \providecommand\transparent[1]{%
    \errmessage{(Inkscape) Transparency is used (non-zero) for the text in Inkscape, but the package 'transparent.sty' is not loaded}%
    \renewcommand\transparent[1]{}%
  }%
  \providecommand\rotatebox[2]{#2}%
  \newcommand*\fsize{\dimexpr\f@size pt\relax}%
  \newcommand*\lineheight[1]{\fontsize{\fsize}{#1\fsize}\selectfont}%
  \ifx\svgwidth\undefined%
    \setlength{\unitlength}{82.40043107bp}%
    \ifx\svgscale\undefined%
      \relax%
    \else%
      \setlength{\unitlength}{\unitlength * \real{\svgscale}}%
    \fi%
  \else%
    \setlength{\unitlength}{\svgwidth}%
  \fi%
  \global\let\svgwidth\undefined%
  \global\let\svgscale\undefined%
  \makeatother%
  \begin{picture}(1,0.97634089)%
    \lineheight{1}%
    \setlength\tabcolsep{0pt}%
    \put(0,0){\includegraphics[width=\unitlength,page=1]{banach-schroder-bernstein_function2.pdf}}%
    \put(0.49996451,0.64341059){\color[rgb]{0,0,0}\makebox(0,0)[t]{\smash{\begin{tabular}[t]{c}$f$\end{tabular}}}}%
    \put(0.49996451,0.38855731){\color[rgb]{0,0,0}\makebox(0,0)[t]{\smash{\begin{tabular}[t]{c}$f$\end{tabular}}}}%
    \put(0.49996451,0.09729694){\color[rgb]{0,0,0}\makebox(0,0)[t]{\smash{\begin{tabular}[t]{c}$g^{-1}$\end{tabular}}}}%
    \put(0.18139817,0.59790094){\color[rgb]{0,0,0}\makebox(0,0)[t]{\smash{\begin{tabular}[t]{c}$A_e$\end{tabular}}}}%
    \put(0.18139817,0.37035376){\color[rgb]{0,0,0}\makebox(0,0)[t]{\smash{\begin{tabular}[t]{c}$A_\infty$\end{tabular}}}}%
    \put(0.18139817,0.09729694){\color[rgb]{0,0,0}\makebox(0,0)[t]{\smash{\begin{tabular}[t]{c}$A_o$\end{tabular}}}}%
    \put(0.81853071,0.46137253){\color[rgb]{0,0,0}\makebox(0,0)[t]{\smash{\begin{tabular}[t]{c}$B_\infty$\end{tabular}}}}%
    \put(0.81853071,0.64341059){\color[rgb]{0,0,0}\makebox(0,0)[t]{\smash{\begin{tabular}[t]{c}$B_o$\end{tabular}}}}%
    \put(0.81853071,0.14280659){\color[rgb]{0,0,0}\makebox(0,0)[t]{\smash{\begin{tabular}[t]{c}$B_e$\end{tabular}}}}%
    \put(0.13588865,0.91646741){\color[rgb]{0,0,0}\makebox(0,0)[t]{\smash{\begin{tabular}[t]{c}$A$\end{tabular}}}}%
    \put(0.77302133,0.91646741){\color[rgb]{0,0,0}\makebox(0,0)[t]{\smash{\begin{tabular}[t]{c}$B$\end{tabular}}}}%
    \put(0.49996451,0.82544865){\color[rgb]{0,0,0}\makebox(0,0)[t]{\smash{\begin{tabular}[t]{c}$\sim$\end{tabular}}}}%
  \end{picture}%
\endgroup%

%% file: banach-tarski.tex
\section{The Banach-Tarski Paradox}
\label{sec:bt}
\subsection{The Strong and Weak forms of the Banach-Tarski Paradox}

We now have the foundation needed to define the Banach-Tarski paradox, which is often claimed to
have first been published in \cite{btart} by Stefan Banach and Alfred Tarski in 1924. However,
a less known fact is that an equivalent version of the Banach-Tarski paradox was already published
ten years earlier, in 1914 by Felix Hausdorff, in \cite[pp.~469--473]{hausdorffbok}. Banach and Tarski
still seem to be the creators of the form referred to as the \textit{strong Banach-Tarski paradox}:

\btprdx

Banach and Tarski note that this result is also valid for higher dimensions but not for lower.
By \textit{nonempty interior}, we mean that there exists an open set entirely contained within
the set. Specifically want there to exist an open ball
entirely contained within the set,
which can always be found if there is a nonempty interior: an open ball centered at any interior point is
contained if we choose a sufficiently small radius. In euclidean space, like $\RE^3$, the existence of
such a ball is equivalent to a nonempty interior. While this result is valid for all dimensions of three or higher, we will only
focus on the case of three dimensions.

Observe that \Cref{b-t1} results in both of the popular descriptions of the paradox stated in
\Cref{quote1,quote2} on page \pageref{quote1}, illustrated by \Cref{fig:illustration}. The first
states that a pea and the sun are $G_3$\-/equidecomposable, ignoring that the mathematical
description of the objects as sets of points fails to capture the aspects of their physical matter.
The second description, duplication of a
ball, follows similarly. However, note that the last result only needs $G_3$\-/paradoxicality of the ball,
as noted after \Cref{paradoxic_duplicate}. This seemingly weaker property is actually equivalent
and is usually the version referred to as \textit{the Banach-Tarski paradox}:

\begin{restatable}[The Banach-Tarski paradox]{theorem}{btprdxweak}
	\label{b-t2}
	$\BLL$, the closed unit ball at origo in $\RE^3$, is
	$G_3$\-/paradoxical.
\end{restatable}
The requirements on the radius and location of the ball are usually not expressed in the theorem, but it will help
reduce the size of the proofs. As we have already seen in \Cref{para_transf} and \Cref{para_scale}, this
result trivially generalizes to spheres of any radius and location.

As previously mentioned, an equivalent result was originally published in 1914 by Hausdorff in
\cite[pp.~469--473]{hausdorffbok}. He built it on top of what is now known
as the \textit{Hausdorff paradox}, using the exact same ideas and
reasoning as Banach and Tarski. This construction will be covered in \Cref{sec:proof}
and the Hausdorff paradox in \Cref{sec:hausdorff}. The result is instead typically credited to Banach
and Tarski, even though arguably the originality in their work was in deducing the strong form from it
and presenting further discussion on measurability. See \Cref{sec:hausdorff} for details on the
Hausdorff paradox, \Cref{sec:non-measurable} for more historical background and \Cref{sec:attribution} for
some discussion regarding the attribution of this result. It was
shown in 1947 that only five pieces are needed to duplicate the ball while any less is
impossible, in \cite{robinsonart}.

\begin{remark}

There is also a concept of sets being \textit{continuously equidecomposable}, where instead of
using elements of $G$ to
map the subsets directly to each other, we use paths in $G$. In $G_3$, this process of
continuous rigid transformation describes the entire process of moving and rotating the subsets, in
contrast to only describing the subsets before and after translation and rotation. Remarkably,
in 2005 it was shown that the underlying process of the equidecomposability giving the paradoxicality of the ball in \Cref{b-t2} can be performed by continuous
equidecomposability without any intersections. In other words: the ball can be decomposed into a finite number of pieces,
which are then moved and reassembled into two copies of the ball, without the pieces ever intersecting. This result was presented in \cite{wilsonart}.
\end{remark}

While the second version is less striking and might at first seem less powerful, the two statements
are equivalent. Showing the equivalence of these two forms of the paradox will be the immediate
next goal. Recall from \Cref{equi_transf}, that two $G_3$\-/equidecomposable sets remain $G_3$\-/equidecomposable
after any of them have been translated. And from the equivalent property for paradoxicality of sets in \Cref{para_transf} and
\Cref{para_scale}, we know that a $G_3$\-/paradoxical set remains paradoxical after being translated
and scaled. Furthermore, remember that multiple equidecomposability properties can be combined from the
transitivity shown in \Cref{equidecomp_transitive}.

\begin{restatable}[Equivalence of the Banach-Tarski paradox forms]{prop}{equivalencebt}
	\label{equivalence_bt}
	The strong and weak forms of the Banach-Tarski paradox, as
	presented in \Cref{b-t1,b-t2}, are equivalent.
\end{restatable}

\begin{proof}
	Suppose that \Cref{b-t1} holds, then clearly \Cref{b-t2} holds since the ball $\BLL$
	has a nonempty interior, and we can easily partition it into two subsets with
	nonempty interiors. From \Cref{b-t1} both are equidecomposable with $\BLL$ from which it follows that
	$\BLL$ is paradoxical. Now assume \Cref{b-t2} holds, and let $A$ and $B$ be as specified.
	Having a nonempty interior, both must contain an open ball inside which we can find a
	closed ball, $K$, by choosing a smaller radius than the open ball. The following steps in
	the proof are illustrated in \Cref{fig:bt_equi}.

	\Cref{para_scale} and \Cref{para_transf} allows us to apply \Cref{b-t2} to $K$ even
	when the ball is not of unit radius or centered at origo; $K$ is $G_3$\-/paradoxical. Just like in
	\Cref{paradoxic_duplicate} this means the ball can be duplicated: it has two subsets, each
	of which is equidecomposable with an identical copy the ball. We can repeatedly perform this
	duplication on the new balls to obtain any finite number of balls, $K_1, K_2, \dots, K_n$.

	Furthermore, since the balls are bounded, and translation is part of $G_3$, we
	can consider the balls to be translated to be pairwise disjoint using the transformation
	invariance of \Cref{equi_transf} and transitivity of \Cref{equidecomp_transitive}. Doing
	so, their union is equidecomposable with the original ball: $K\sim\bigcup_{i=1}^nK_i$. 

\begin{figure}[h]
	\centering
	\def\svgwidth{\columnwidth}
	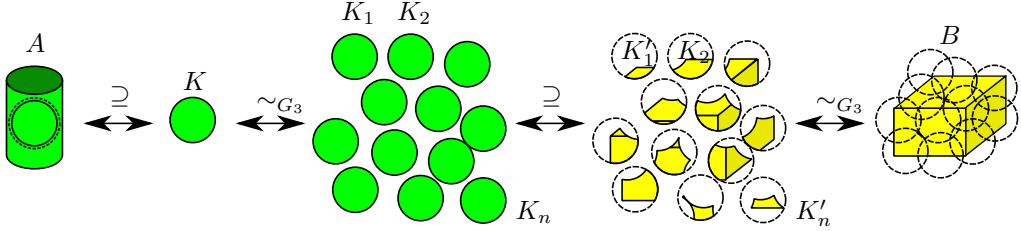
	\caption{The nonempty interior of $A$ contains a closed ball that can be duplicated
	repeatedly to cover $B$. A suitable subset of this covering gives us equidecomposability
	between $B$ and a subset of $A$.}
	\label{fig:bt_equi}
\end{figure}

	Since $B$ is bounded, there is a sufficiently large $n$ such that this finite collection of balls
	can be translated to fully cover $B$, and let these translations be $g_1,g_2,\dots,g_n\in
	G_3$, so that $B\subseteq \bigcup_{i=1}^ng_iK_i$. Note however that some of the translated
	balls might overlap, so this does not describe an equidecomposability between
	$\bigcup_{i=1}^nK_i$ and $\bigcup_{i=1}^ng_iK_i$.

	However, we can now define subsets $K'_i\subseteq K_i$ by removing all redundant points 
	mapped to the same points in $B$ to keep these subsets pairwise disjoint after the translations. 
	Furthermore, their union is made equidecomposable with $B$ by only including points mapped
	into $B$.
	Consider the points in $B$ that are contained in $g_1K_1$, and translate them back into a
	subset of $K_1$ by $K'_1\coloneqq g^{-1}_1(Y\cap g_1K_1)$. We now repeat this for $K_2$ and
	the remaining points $B\setminus g_1K_1$ by $K'_2\coloneqq g^{-1}_2((Y\setminus g_1K_1)\cap
	g_2K_2)$, for $K_3$ using the points $B\setminus (g_1K_1\cup g_2K_2)$ and so forth
	up until $K'_n\coloneqq g^{-1}_n((Y\setminus (g_1K_1\cup\dots\cup g_{n-1}K_{n-1}))\cap
	g_nK_n)$.

	See \Cref{fig:bt_equi} for an illustration of the process and result.
	The sets $K'_1,K'_2,\dots,K'_n$ are pairwise disjoint and this is also true for their
	translations $g_1K'_1,g_2K'_2,\dots,g_nK'_n$ satisfying $\bigcup_{i=1}^ng_iK'_i=B$. Thus we
	see that $B$ is equidecomposable with $\bigcup_{i=1}^nK'_i$,
	which is a subset of $\bigcup_{i=1}^nK_i$, which in turn is equidecomposable with
	$K\subseteq A$. By \Cref{equidecomp_subsets} and the transitivity in
	\Cref{equidecomp_transitive} means $B$ is equidecomposable with a subset of $A$.

	Repeating this process also shows that $A$ is equidecomposable within a subset of $B$. From
	this we can now apply \Cref{thm:b-s-b} to obtain equidecomposability of $A$ and $B$.
\end{proof}
Obtaining a proof of \Cref{b-t2} will be the goal of the following two sections.

\subsection{Non-measurable Sets and their History}
\label{sec:non-measurable}





Now that we have defined the Banach-Tarski paradox and understand its meaning, we can end the
section by showing its most important result: the existence of non-measurable sets, primarily
for the Lebesgue measure. This was the original goal of Hausdorff's construction. However, let us
first reflect on the historical significance of this result:

In the 19th century, the foundations of mathematics began to change significantly, leading to more
robust reasoning with a stricter, axiomatic approach. Towards the end of the century, the concept
of measures entered the mathematical world. By the shift of the century, the field of measure
theory was quickly reshaping analysis into its modern form with concepts such as the Lebesgue
integral and measure. Measures were still under active research.

In 1905 mathematician Giuseppe Vitali from Italy published in \cite{vitalibok} a remarkable construction allowing the interval
$[0,1)$ in $\RE$ to be split into a countably infinite number of disjoint subsets. Each one being identical except for being
translated along the real axis, now called the \textit{Vitali set}. Should these sets be measurable
by a countably additive measure that is invariant under translation, the measure of
$[0,1)$ would be equal to the infinite sum of the measures of these subsets. Since
each of these sets must have the same measure, the sum could only be 0 or infinite. This proved the existence of non-measurable
subsets of $\RE$ for measures like
the Lebesgue measure. But this construction also allowed for a curious paradox: not only could these
subsets be reassembled into $[0,1)$, but by further translating them, their union could also become
$[0,2)$. The construction also had another trick up its sleeve: this reasoning could also be
applied to a circle. This meant the circles in $\RE^2$ could also be duplicated, using an infinite
number of pieces.

Vitali's result would encourage more research into \textit{finite additivity} of measures. Surely, limiting
ourselves to only finite sums of measures would avoid this problem? Any hope that restricting measures to
be finitely additive would avoid these kinds of results would soon be shattered. In 1914
German mathematician Felix Hausdorff, taking inspiration from Vitali's construction
and using similar reasoning on a combination of different fields of mathematics, published a
groundbreaking result for $\RE^3$, in \cite[pp.~469--473]{hausdorffbok}. And this time, the paradoxical
outcome would receive full
attention: using only a finite number of subsets and purely rotating them, he was able to
duplicate a subset of a sphere, and from this result, he deduced what we call the Banach-Tarski
paradox: that the entire ball could also be duplicated! What
this meant for the research into measures was only outshined by the paradoxical result itself.

Hausdorff's work also laid the foundation for Polish mathematicians Stefan Banach and Alfred Tarski, who 1924
published \cite{btart} presenting a collection of various results regarding the decomposition of sets in different
dimensions. Of the results, the most well known has become known as the Banach-Tarski paradox, and
its strong form built on top of the results of Hausdorff, using a more rigorous definition of
equidecomposability and Banach's variant of the Schröder-Bernstein theorem applied to
equidecomposability.

To finish this
section, we can now observe one of the key results of the Banach-Tarski paradox: the existence of
non-measurable sets in $\RE^3$. This was the main result that Hausdorff wanted to show when
he created his version of the Banach-Tarski paradox, as in \Cref{b-t2}, and in the process, what is
now known as the Hausdorff paradox, the key
to prove the Banach-Tarski paradox and the goal of the next section.

Let us begin by showing how the Banach-Tarski paradox shows the existence of non-measurable subsets of
$\RE^3$, primarily for the Lebesgue measure, but also for many others:
\begin{example}[Existence of non-measurable subsets of $\RE^3$]
	\label{ex:unmeasurable}
	Let $\lambda$ be a measure on $\RE^3$, such as the Lebesgue measure, which is invariant under $G_3$ and at least
	\textit{finitely additive}. That is, for any measurable $X\subseteq
	RE^3$, $\lambda{X}=\lambda{gX}$ for any $g\in G_3$, and any disjoint and measurable
	$X_1,X_2,\dots,X_n\in \RE^3$, $\lambda(X_1\cup X_2\cup \dots\cup X_n)=
	\lambda(X_1)+\lambda(X_2)+\dots+\lambda(X_n)$.

	Now assume that all subsets of $\RE^3$ are measurable. Then two equidecomposable sets
	must have the same measure: assume $A$ and $B$ are equidecomposable, then from the
	definition there exist decompositions $A=A_1\cup A_2\cup\dots\cup A_n$ and $B=B_1\cup
	B_2\cup\dots\cup B_n$ together with $g_1,g_2,\dots,g_n\in G_3$ such that $A_i=g_iB_i$ for
	$i=1,2,\dots,n$. Then by the properties of $\lambda$ we get:
	\begin{align*}
		\begin{split}
			\lambda(B)&=\lambda(B_1\cup B_2\cup \dots\cup B_n)\\
			&=\lambda(B_1)+\lambda(B_2)+\dots+\lambda(B_n)\\
			&=\lambda(g_1B_1)+\lambda(g_2B_2)+\dots+\lambda(g_nB_n)\\
			&=\lambda(A_1)+\lambda(A_2)+\dots+\lambda(A_n)\\
			&=\lambda(A_1\cup A_2\cup \dots\cup A_n)=\lambda(A)\text{.}
		\end{split}
	\end{align*}

	However, the Banach-Tarski paradox, in \Cref{b-t2}, gives us the following paradoxical decomposition of the ball:
	$\BLL=A\cup B$, where $\BLL$ is equidecomposable with $A$ and $B$. This means
	\[\lambda(\BLL)=\lambda(A)+\lambda(B)=\lambda(\BLL)+\lambda(\BLL)\text{.}\]
	Assuming $\lambda(\BLL)$ is nonzero and not infinite, this is a contradiction, and it is
	impossible for all of the subsets to be measurable. This is the
	case for the Lebesgue measure: for a ball of radius $r$, its measure is its volume, $4/3\pi
	r^3$, which is nonzero and finite.
\end{example}
Some properties of the Lebesgue measure are covered in \cite{kursbok}. Of course, this result could
just as easily have been deduced from the strong form, \Cref{b-t1}, but this more closely reflects
Hausdorff's ambitions. Again, the weak form was presented by Hausdorff about ten years before
Banach and Tarski, with the explicit goal of showing this result. This result can be applied to
other measures than the Lebesgue, especially since the measure only needs to be finitely measurable.

A far more straightforward approach to showing the existence of non-measurable subsets,
for \textit{countably additive} measures such as the Lebesgue measure, is the famous Vitali set
previously mentioned. This is a well-known result, and it will not be covered here; see
\cite[Thm.~1.4.9]{kursbok} for details. While it is significantly weaker than the result of
Hausdorff, Banach and Tarski, by requiring a countably additive measure, it works in the one-dimensional case of $\RE$, instead of $\RE^3$ and higher.

When Giuseppe Vitali introduced his famous set in 1905, he also provided a lesser-known application of the
construction to a circle, which sometimes has received the informal nickname the \textit{Vitali paradox}.
Briefly, the classic Vitali set construction relies on the creation of a subset of the interval $[0,1)$
from which cosets are created by adding numbers, essentially moving the subset along the $x$-axis. If
the interval is instead turned into a circle, the movement becomes equivalent to a rotation. The
construction presented here is based on \cite[Thm.~1.5]{wagonbok} but follows Vitali's original idea.

This construction has similarities with \Cref{ex:irrational} and we could consider each rotation to
be either a \textit{rational rotation} or an \textit{irrational rotation}.
\begin{example}[The Vitali paradox]
	\label{ex:vitali}
	Let $SO(2)$ be the group of rotation around origo in $\RE^2$ and $\CIRC$ be the unit circle
	centered at origo in $\RE^2$. There is an equivalence relation on the points of
	$\CIRC$ such that $\vect{p_1}\sim\vect{p_2}$ if the central angle between them is a
	rational multiple of $2\pi$. See \Cref{fig:vitali}.

	\begin{figure}[h]
		\centering
		\scalebox{1.0}{\input{figs/vitali.tex}}
		\caption{The central angle $\theta$ between $\vect{p_1}$ and $\vect{p_2}$ on
		$\CIRC$. If $\theta$ is a rational multiple of $2\pi$ then $\vect{p_1}\sim\vect{p_2}$.}
		\label{fig:vitali}
	\end{figure}

	Each point on $\CIRC$ belongs to one of these orbits, and the collection of unique orbits
	forms a partition of $\CIRC$. We use the \nameref{ax_choice} to
	select one point from each orbit and call this set of points $C$. Each point on the circle
	can be obtained by taking the point of $C$ that belongs to the same orbit and rotating it
	by one specific rotation of a rational number in $[0,1)$ multiplied by $2\pi$. Observe that
	each point can be obtained by precisely one rotation like this.

	Now consider all rational numbers in $[0,1)$: since they are countable, we can enumerate
	them as $q_1,q_2,\dots$. For each $q_n$ there is a rotation $r_n$ in $SO(2)$ by $2\pi q_n$
	radians.
	By defining every rotation of $C$ by $q_n$ as $C_n\coloneqq q_nC$ we have:
	\[\CIRC=C_1\cup C_2\cup C_3\cup\dots\text{.}\]

	We can utilize the countable infinity we have obtained: we can recreate $\CIRC$ using only
	the sets of even or odd index by applying a suitable rotation to them:
	\[\CIRC=C_1\cup C_3+(q_2-q_3)\cup C_5+(q_3-q_5)\cup\dots\text{,}\]
	and
	\[\CIRC=C_2+(-q_2)\cup C_3+(q_2-q_3)\cup C_5+(q_3-q_5)\cup\dots\text{.}\]

	Thus $\CIRC$ has a decomposition into two subsets, $C_1\cup C_3\cup\dots$ and $C_2\cup
	C_4\cup\dots$,
	such that each of them is \textit{countably $SO(2)$\-/equidecomposable} with
	$\CIRC$, from which $\CIRC$ is \textit{countably $SO(2)$\-/paradoxical}.
\end{example}
Beyond the remarkable result, we will see a surprising similarity between the steps above and in
the construction of the \textit{Hausdorff paradox} in the next section. There, a subset of
$\mathbb{S}^2$, the unit sphere, will be rotated by
elements of $SO(3)$. This includes a similar use of the \nameref{ax_choice} to select points from
orbits. However, Hausdorff goes one step further by carefully selecting rotations that cause a curious pattern, from which certain unions can be made on a countable number of subsets, filling the same role as
the subsets $C_n$ above. This reduces them to
only a finite number of subsets resulting in finite paradoxicality, which is much stronger.

The work by Vitali is a clear source of inspiration for Hausdorff. And it will be worthwhile to
keep this example in mind when reading the next section.

%% file: figs/bt_equi.pdf_tex
\begingroup%
  \makeatletter%
  \providecommand\color[2][]{%
    \errmessage{(Inkscape) Color is used for the text in Inkscape, but the package 'color.sty' is not loaded}%
    \renewcommand\color[2][]{}%
  }%
  \providecommand\transparent[1]{%
    \errmessage{(Inkscape) Transparency is used (non-zero) for the text in Inkscape, but the package 'transparent.sty' is not loaded}%
    \renewcommand\transparent[1]{}%
  }%
  \providecommand\rotatebox[2]{#2}%
  \newcommand*\fsize{\dimexpr\f@size pt\relax}%
  \newcommand*\lineheight[1]{\fontsize{\fsize}{#1\fsize}\selectfont}%
  \ifx\svgwidth\undefined%
    \setlength{\unitlength}{474.98909895bp}%
    \ifx\svgscale\undefined%
      \relax%
    \else%
      \setlength{\unitlength}{\unitlength * \real{\svgscale}}%
    \fi%
  \else%
    \setlength{\unitlength}{\svgwidth}%
  \fi%
  \global\let\svgwidth\undefined%
  \global\let\svgscale\undefined%
  \makeatother%
  \begin{picture}(1,0.22829682)%
    \lineheight{1}%
    \setlength\tabcolsep{0pt}%
    \put(0.11303736,0.12255719){\color[rgb]{0,0,0}\makebox(0,0)[t]{\smash{\begin{tabular}[t]{c}$\supseteq$\end{tabular}}}}%
    \put(0,0){\includegraphics[width=\unitlength,page=1]{bt_equi.pdf}}%
    \put(0.53936291,0.12255719){\color[rgb]{0,0,0}\makebox(0,0)[t]{\smash{\begin{tabular}[t]{c}$\supseteq$\end{tabular}}}}%
    \put(0,0){\includegraphics[width=\unitlength,page=2]{bt_equi.pdf}}%
    \put(0.27158947,0.11893656){\color[rgb]{0,0,0}\makebox(0,0)[t]{\smash{\begin{tabular}[t]{c}$\sim_{G_3}$\end{tabular}}}}%
    \put(0,0){\includegraphics[width=\unitlength,page=3]{bt_equi.pdf}}%
    \put(0.82423378,0.11893656){\color[rgb]{0,0,0}\makebox(0,0)[t]{\smash{\begin{tabular}[t]{c}$\sim_{G_3}$\end{tabular}}}}%
    \put(0.18710954,0.13563839){\color[rgb]{0,0,0}\makebox(0,0)[t]{\smash{\begin{tabular}[t]{c}$K$\end{tabular}}}}%
    \put(0.34816582,0.20511367){\color[rgb]{0,0,0}\makebox(0,0)[t]{\smash{\begin{tabular}[t]{c}$K_1$\end{tabular}}}}%
    \put(0.40343023,0.20511367){\color[rgb]{0,0,0}\makebox(0,0)[t]{\smash{\begin{tabular}[t]{c}$K_2$\end{tabular}}}}%
    \put(0.62448797,0.16563909){\color[rgb]{0,0,0}\makebox(0,0)[t]{\smash{\begin{tabular}[t]{c}$K'_1$\end{tabular}}}}%
    \put(0.67975247,0.16563909){\color[rgb]{0,0,0}\makebox(0,0)[t]{\smash{\begin{tabular}[t]{c}$K_2$\end{tabular}}}}%
    \put(0.02921119,0.17353401){\color[rgb]{0,0,0}\makebox(0,0)[t]{\smash{\begin{tabular}[t]{c}$A$\end{tabular}}}}%
    \put(0.93238981,0.18142893){\color[rgb]{0,0,0}\makebox(0,0)[t]{\smash{\begin{tabular}[t]{c}$B$\end{tabular}}}}%
    \put(0.521854,0.00774074){\color[rgb]{0,0,0}\makebox(0,0)[t]{\smash{\begin{tabular}[t]{c}$K_n$\end{tabular}}}}%
    \put(0.79817615,0.00774074){\color[rgb]{0,0,0}\makebox(0,0)[t]{\smash{\begin{tabular}[t]{c}$K'_n$\end{tabular}}}}%
  \end{picture}%
\endgroup%

%% file: figs/vitali.tex
\begin{tikzpicture}[thick]
	\path (0,0) node[circle,draw,inner sep=0pt,minimum width=100](circle){};
	\coordinate (A) at (0,0);
	\coordinate (B) at (circle.east);
	\coordinate (C) at (circle.north east);

	\draw (B) -- (A) -- (C);
	\draw [circle, inner sep=0pt, minimum size=5.0] (B) node [fill, label=right:$\vect{p}_1$] {} -- (A)
	-- (C) node [fill, label=right:$\vect{p}_2$] {};

	\node at (-1.5,1.5) {$\CIRC$};

	\tkzMarkAngle[mark=none,color=cyan,size=0.6](B,A,C)
	\tkzLabelAngle[pos=0.8,color=cyan](B,A,C){$\theta$}
\end{tikzpicture}

%% file: hausdorff.tex
\section{The Hausdorff Paradox}
\label{sec:hausdorff}


An essential part of the proof of the Banach-Tarski paradox is the \textit{Hausdorff paradox}.
Even though the result is not as immediately striking, it is far more groundbreaking and uses
much more advanced reasoning. It was published by Felix Hausdorff in 1914 with the goal of showing
the existence of a stronger form of non-measurable sets compared to the previous results by Vitali,
by the idea of \Cref{ex:unmeasurable}. In doing so, he also provides all the reasoning we will see
in \Cref{sec:proof} where the weak Banach-Tarski paradox of \Cref{b-t2} is deduced from the
Hausdorff paradox of this chapter.

The result and the proof of this paradox were published the same year in both
\cite[pp.~469--473]{hausdorffbok} and \cite{hausdorffart}. The latter does not include the further
reasoning of \Cref{sec:proof}. In short, the original formulation translates to:
\begin{theorem}[The Hausdorff paradox]
	\label{hausdorff_original}
	The sphere, $K$, can be divided into the disjoint sets $A$, $B$, $C$ and $Q$ such that
	$A$,$B$ and $C$ are congruent with each other and with $B\cup C$, and $Q$ is countable.
\end{theorem}
In Hausdorff's own words, this effectively means that \enquote{a half} and \enquote{a third of the
sphere} are
congruent. The term \textit{congruent} here means one set can be transformed into the
other by some element of $G_3$ together with a possible reflection. As the proof does not use
reflections and only rotations as in $SO(3)$, this would seem unnecessarily broad and probably reflects
common terminology at the time. The concept of equidecomposability and paradoxicality had not
been defined at this time.

From this result, it is easy to see that the sphere with the points of $Q$ removed is
$SO(3)$\-/paradoxical: let $\sim$ represent congruence, then $K\setminus Q=A\cup B\cup C$ and
$A\cup B\sim A\cup B\cup C$ together with $C\sim B\cup
C\sim A\cup C\sim A\cup B\cup C$ shows us that $A\cup B$ and $C$ gives a paradoxical decomposition
of $K\setminus Q$. While the original result was used by Hausdorff as well as Banach and Tarski, the
paradoxical property is the only result we need. By simplifying the theorem we also simplify its
proof.
\begin{restatable}[The simplified Hausdorff paradox]{theorem}{hausdorffprdx}
	\label{hausdorff}
	Let $\SPH$ be the unit sphere centered at origo in $\RE^3$.
	There is a
	countable subset $D$ of $\SPH$ such that $\SPHMD$ is $SO(3)$\-/paradoxical.
\end{restatable}
To reduce the size of the proof, the statement is restricted to a unit sphere centered at origo. As we have already seen in \Cref{para_transf,para_scale}, this result trivially generalizes to spheres of any radius and location. While there are differences, both results and proofs share the same core ideas, including the need to remove a countable subset from the sphere. In contrast to the original, our proof will be divided into several individual steps and utilize modern concepts and terminology.

The work by Vitali is a clear source of inspiration to Hausdorff. It is worthwhile to keep the construction of Vitali from \Cref{ex:vitali} in mind as many of the core ideas will return throughout this section. Working in $\RE^3$ or higher dimensions is vital to this result, and we will see why it is impossible in lower dimensions. Just like in Vitali's construction, we will require the \nameref{ax_choice}, \Cref{ax_choice}, to pick points of certain orbits in the proof of \Cref{G_7}.

Where Vitali's construction in \Cref{ex:vitali} gave us a countable collection of sets with
paradoxical properties, Hausdoff's construction takes it one step further by
putting these sets back together into just two: a paradoxical decomposition. What enables this is a
group of rotations with the curious property of being paradoxical on itself.

\begin{figure}[h]
	\centering
	\def\svgwidth{\columnwidth}
	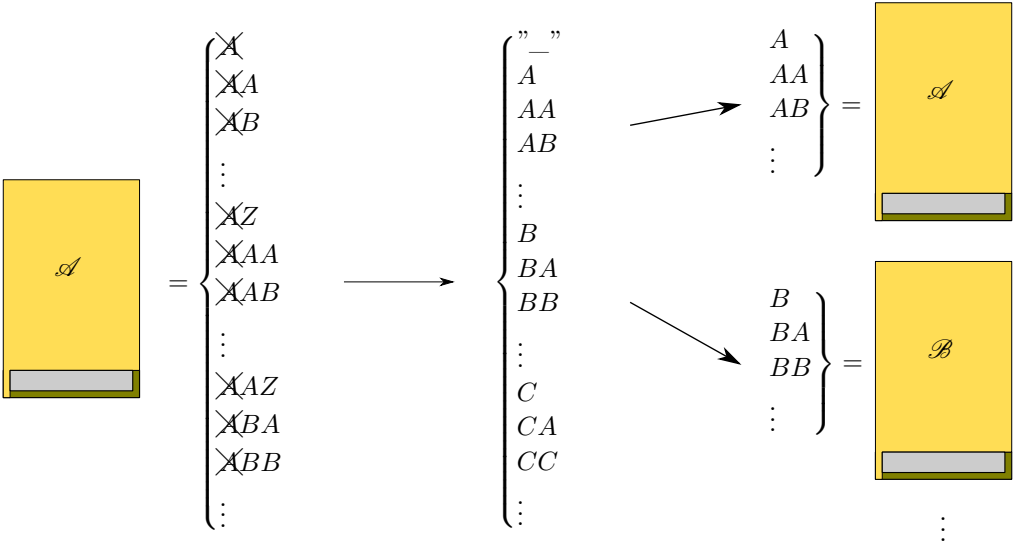
	\caption{The words in volume \textit{A} can be modified to become all words
	in the entire Hyperwebster. Here $\_$ represents a kind of empty word that remains of the word \WRD{A}.}
	\label{fig:hyperwebster}
\end{figure}

What enables this property on the group is that it is a \textit{free group}, covered in
\Cref{def:free_group}. The existence of this group is given by \Cref{G_6}, and its proof is by far
the biggest and arguably most important in this thesis. Free groups, in general, have received much
research and offer many exciting properties beyond our usage case. Free groups will require some
work to define, but the following example, inspired by \cite[pp.~175--176]{mysbok}, can help introduce the property and demonstrate its usefulness:

\begin{example}[The Hyperwebster]
	\label{ex:hyperwebster}
	The \textit{Hyperwebster} is an overambitious work that provides a complete list of all words in
	the form of every possible combination of letters of any finite length. It is a countably
	infinite collection of all letter combinations and not just the words we know. Sadly
	it does not attempt to assign a definition to the words. It is divided into 26 volumes,
	\WRD{A} to \WRD{Z}, each one contains all words that start with the volume's unique letter:
	\textit{A} contains \WRD{ALFRED}, \WRD{AALFRED}, \WRD{ABANACH} and similar.

	We now vandalize the first volume, \textit{A}, by crossing over the first letter of each word.
	However, the resulting list of words still contains all of its original
	words: the word \WRD{ALFRED} returns when the first letter of the word \WRD{AALFRED} gets crossed
	over. However, we also see that the words from all the other volumes also appear: the word
	\WRD{BANACH} from volume \WRD{B} appears when the first letter of the word \WRD{ABANACH} gets crossed
	over. See \Cref{fig:hyperwebster}.
	Remarkably, by removing the first letter of each word in one volume, we have managed to
	recreate all the other volumes and thus the entire Hyperwebster.
\end{example}

If we ignore for a moment that the volumes are of infinite length, the result is truly paradoxical to our intuition. However, just like the irrational rotation in \Cref{ex:irrational} we are once again dealing with infinite collections, this time of the words. The Hyperwebster and each of its volumes are countably infinite. It is still a powerful result and will be the essential inspiration for this section. Note that the minimum number of letters required for this paradox is two; with only one letter, we would only have one volume. Any attempt to add and remove letters to these words will not accomplish anything useful.

We have seen a way of modifying a subset of words into the entire set. This has a striking resemblance to our mathematical definition of paradoxicality, \Cref{def:paradoxical}. Under certain conditions, this procedure can be applied almost identically to certain groups: the alphabet and dictionary are instead a set of group elements and the group it spawns. The words represent compositions of the elements from the set, from which all elements from the group are obtained. Like in the example, we create subsets of the group based on the leftmost element in the composition. By manipulating the sets using group actions, we can significantly increase the elements they contain. Just like in the example, we will need at least two letters.

Based on the formulation of the Hausdorff paradox in \Cref{hausdorff}, the reader might already have guessed that our end goal is to use elements from $SO(3)$ as our dictionary. However, we first need to give a precise mathematical definition of \textit{letters} and \textit{words}, which is the goal of the following subsection. There will be some requirements on our choice of letters; we want our words to represent unique elements in the group. Otherwise, the process of adding or removing letters to words will not work as in the example.

These definitions will require some work which will be the goal in the following subsection. However, this work is not without some rewards. As a practical upshot, we obtain a proper definition of an \textit{empty word}, which solves the one leftover element we obtained in our example when we removed the letter \WRD{A} from the word \WRD{A}, which we simply discarded. Furthermore, having inverses as part of our letters gives us a mathematical way of removing the leftmost letter in words: we just apply the inverse of the letter. This would otherwise not have been an obvious approach.

\subsection{Words and Free Groups}
We now define a mathematical version of words, letters and similar for groups and their elements.

\begin{definition}[Generated group]
	Let $G$ be a group with a subset $S \subset G$. The smallest
	subgroup $F\subseteq G$ that includes $S$ is called the group, or subgroup,
	\textit{generated} by $S$. $S$ is a \textit{generating set} of $G$ and, the elements in $S$
	are the \textit{generators} of $G$.
\end{definition}

\begin{example}
	Let $G$ be a group. Then $G$ itself is a generating set of $G$.
\end{example}

\begin{definition}[Word]
	\label{def:word}
	Let $S=\set{a,b,c\dots}$ be a finite or countable set of \textit{symbols}, then a
	\textit{word} in $S$ is a finite sequence of \textit{letters}, $(s_1,s_2,\dots,s_n)$, where
	each $s_i$ is one of $a,b,c,\dots$ or $a^{-1},b^{-1},c^{-1},\dots$.
	The number of letters in the word, $n$, is its \textit{length}. The set $S$ is called an
	\textit{alphabet} and its elements are \textit{symbols}. Each symbol gives us two letters.

	Words are often written without parentheses and commas between letters. To avoid
	confusion, words and letters will always be written stylized as
	$\texttt{s}_\texttt{1}\texttt{s}_\texttt{2}\dots\texttt{s}_\texttt{n}$.
\end{definition}

\begin{definition}[Representation]
	Each symbol in $S$ is defined to \textit{represent} some element in a group $G$ by a
	mapping $S \to G$.
	If \WRD{a} represents $g\in G$, then \WRD{a\INV} represents the inverse: $g^{-1}\in G$.
	The word $\texttt{s}_\texttt{1}\texttt{s}_\texttt{2}\dots\texttt{s}_\texttt{n}$ represents the result of the
	composition of the group elements represented by the letters,
	$g_1g_2\dots g_n\in G$. The special case when $n=0$ is called the
	\textit{empty word} and is always defined to represent $e\in G$.
\end{definition}

We will always consider $S$ as a subset of a group, $S\subseteq G$, where each
letter is implicitly defined to represent the same element in the group. Note that $S$ is a set and
not necessarily a subgroup.
Since the words allow us to represent all combinations of the group elements and their inverses,
including the identity element, we see the following:
\begin{remark}[Group element representations by words]
	Each element in the group generated by $S$ can be represented as a word in $S$.
	The element $e$ can be represented by the empty word, in case $e\not\in S$.
\end{remark}
$S\subseteq G$ might contain both an element and its inverse from the group $G$. It might also contain
the identity element. In these cases, it is vital to understand the difference between a letter
like \WRD{a\INV} based on the \textit{symbol} $a$ from our alphabet $S$, and when the
\textit{symbol} $a^{-1}$ itself has been taken directly from $S$, assuming both are in $S$. These two letters are different and result in entirely different words, even though they represent the same element in $G$. Similarly, both the word \WRD{e} and the empty word are different words representing the same element. Apart from the awkward ambiguity of writing \WRD{a\INV} in these cases, which is merely a syntactical issue, we will find that identifying and avoiding these kinds of redundant words will be important. To clarify, consider:
\begin{example}[Generators of integers]
	\label{ex:integers}
	Let $S=\set{-1,0,1}\subset \Z$, where $\Z$ is the group of integers under the operation of
	addition. Then it is easy to see that $S$ generates $G$; for example, $3\in\Z$ can be
	represented by the word \WRD{111}. However, words like \WRD{-1} and \WRD{1\INV} are different
	words in $S$ even though they represent the same element, $-1$, in $Z$. Similarly, words
	like \WRD{1-1}, \WRD{-11}, \WRD{11\INV}, \WRD{1\INV1}, \WRD{0} and the empty word are all different
	words that represent the same identity element $0\in Z$. And even they are not the only ones.

	In fact, it would have been enough to have just $S=\set{1}$ or $S=\set{-1}$ to generate
	$\Z$. We will revisit this in \Cref{ex:integer_generators}.
\end{example}

Consider a word representing an element of a group $G$ generated by $S$. Any sequence of letters of the form \WRD{ss\INV} and \WRD{s\INV s} in the word is redundant: they do not affect the element of $G$ represented since $ss\INV=s\INV s=e$ in the group. Consider \WRD{1-1} which represents $0\in G$ in the example above. They can be removed from the word without changing the element it represents, even if the result is the empty word. This observation gives us the concept of \textit{reduction} of words:
\begin{definition}[Reduced word]
	\label{def:reduced_word}
	A word in $S$ is \textit{reduced} if a symbol and its inverse never occur next to each
	other. That is, the word contains no consecutive letters of the
	form \WRD{ss\INV} or \WRD{s\INV s} for any $s\in S$.
	Any word in $S$ can be turned into reduced form by removing all such pairs. This
	process is called \textit{reduction}, and the reduced word will still represent the same
	element in $G$.
\end{definition}

Note that even if we only use reduced words, they may still not represent unique elements in the group. As noted above, if $S$ generates a group and happens to contain both an element and its inverse in $G$, let us say call them $a$ and $b$, then \WRD{ab} is a reduced word but still represents the same element in $G$ as the empty word, as illustrated by $a=1$ and $b=-1$ in \Cref{ex:integers}. Similarly, if $c\in S$ represents the identity element of the group, then \WRD{c} and the empty word are both
reduced words that represent the identity element, as illustrated by $c=0$ in \Cref{ex:integers}. This also occurs for words like \WRD{aa} where $a$ is its own inverse in $G$.

However, even if $S$ does not contain $e$ or both a group element and its inverse, there can still be group elements that are
represented by several different reduced words. For example, say we have $ab=c$ in the group generated by $S$. If $a,b,c\in S$ then both the words \WRD{ab} and \WRD{c} are reduced words in $S$ that represent the same element in $G$. Equivalently, \WRD{abc\INV} represents the identity, just like the empty word. Similarly $S$ might contain a cyclic element, $a$, such that $aaa=e$ in $G$, meaning the word \WRD{aaa} and the empty word represent the same element. 

If we want to use our mathematical words with similar results to what we did with the Hyperwebster in \Cref{ex:hyperwebster}, each reduced word needs to have a unique meaning. Otherwise, we get the situation as in \Cref{ex:integers}. This important concept is given a mathematical definition in \textit{free groups} and \textit{free generators}. The definition might not at first appear strong enough to prevent all cases mentioned above, but \Cref{free_remark} will show that it really is.

\begin{definition}[Free group]
	\label{def:free_group}
	Let $G$ be a group generated by the set $S$. Then $G$ is \textit{free} if no nonempty
	reduced word in $S$ represents the identity element in $G$. We might also say that $G$ is
	\textit{freely generated} by $S$ or that $S$ is a \textit{free generating set} of $G$. The
	elements of $S$ are \textit{free generators} of $G$, and the number of free generators, the
	cardinality of $S$, is the \textit{rank} of the group. $G$ is said to be \textit{free on
	$n$ generators}, where $n$ is its rank.
\end{definition}

\begin{definition}[Equivalent definition of free group]
	\label{free_remark}
	An equivalent definition of $G$ being freely generated by $S$ is that every element
	in $G$ can be represented by only one reduced word; that is, no two different reduced words
	represent the same element in $G$.
\end{definition}
\begin{proof}[Proof of the equivalence of \Cref{def:free_group,free_remark}]
	Assume that we have two reduced words in $S$ describing the same element in $G$, that is:
	\[s_1s_2\dots s_n x =s'_1s'_2\dots s'_m x \quad \forall x\in X,\]
	then by picking one side and repeatedly applying the inverse of the leftmost element
	to both sides, we get:
	\[x=s^{-1}_ns^{-1}_{n-1}\dots s^{-1}_1s'_1s'_2\dots s'_m x \quad \forall x\in X,\]
	and the nonempty word on the other side can further be reduced to a nonempty word;
	otherwise, the
	original two words would be equal. Thus we have a reduced nonempty word which represents the identity element. And vice versa: if
	there is a reduced nonempty word such that:
	\[x=s_1s_2\dots s_n x \quad \forall x\in X,\]
	then we immediately have two words, the empty word, and $\texttt{s}_\texttt{1}\texttt{s}_\texttt{2}\dots\texttt{s}_\texttt{n}$, representing
	the identity element. Thus \Cref{def:free_group,free_remark} are equivalent.
\end{proof}

While the precise definition and terminology for words and free groups vary, the underlying idea is
always the same. We want to find operators from a group that generates the group, or a subgroup, so
each element can only be reached uniquely from reduced compositions of these generating operators.
For an attempt at a simpler definition, see \cite[pp.~421]{kursbok}. However, be aware that this
approach requires some additional explicit requirements on the generating set, a cost to be paid in
exchange for simplicity. This approach is mainly based on \cite{combibok}.

\begin{remark}
	\label{remark:freeresults}
	The area of free groups has received much research, not at least due to the works
	by Hausdorff, Banach and Tarski. More on the subject can be found in \cite{wagonbok} and
	\cite{combibok}. Among their many interesting properties is that any two free groups are
	isomorphic if, and only if, they have the same rank. We will not have a use for this.

	Another property that is somewhat related to our usage is that all free groups of
	rank $2$ or higher must be non-abelian. Otherwise, changing the order of two letters in
	a reduced word would result in a different word representing the same element.
\end{remark}

\begin{example}[Free generators of integers]
	\label{ex:integer_generators}
	Recall \Cref{ex:integers}. As we observed, the subset $\set{1}\subset \Z$ generates $\Z$.
	Similarly, the subset $\set{2}\subset \Z$ generates $2\Z$, the group of even integers under
	the operation of addition. These words would be of the form \WRD{222} and
	\WRD{2\INV 2\INV} representing the elements $6$ and $-4$ in $2Z$). Clearly, $\set{1}$ and
	$\set{2}$ are free generators of $\Z$ and $2\Z$, respectively. Meaning both $\Z$ and $2\Z$
	are free of rank $1$. Another, and
	only other, free generating subset of $\Z$ is $\set{-1}$. As we observed in
	\Cref{ex:integers}, $\Z$ is not \textit{freely} generated by $S=\set{-1,0,1}$.

	Since these two groups are free of rank $1$, they must also be isomorphic, as noted in
	\Cref{remark:freeresults}. This
	can be verified by the classic isomorphism $x \mapsto 2x$, which maps $Z\to2\Z$. However, properties
	of isomorphism will not be of further interest to us.
\end{example}

Another example, which will be important, can be based on the irrational rotation in \Cref{ex:irrational}:
\begin{example}
	Recall the operator $\rho$ of irrational rotation from \Cref{ex:irrational}. The collection of all repeated
	rotations by $\rho$ and its inverse, the counter-rotation by $\rho^{-1}$, together with the
	identity, no rotation by $\rho^0=e$, can be described by $\set{\rho^n|n\in\Z}$. This is a subgroup of all rotations
	around origo in $\RE^2$. 

	Furthermore, this group is freely generated by the irrational rotation, $\rho$, using the
	carefully selected angle $\theta$. This observation will be useful later in this section.
	\label{ex:irrationalgen}
\end{example}

\subsection{Proof of the Hausdorff Paradox}
Now we are ready to prove the Hausdorff paradox. We begin by proving two crucial propositions.

We begin with a result in the spirit of the Hyperwebster of \Cref{ex:hyperwebster}, using reduced words and free groups. Note that in the following theorem, similarly to the example, our alphabet needs at least two symbols:

\begin{prop}
	\label{G_5}
	If the group $G$ is free of rank 2, then $G$ is $G$\-/paradoxical. Here $G$ is seen as acting
	on itself.
\end{prop}
\begin{proof}
	Let $G$ be freely generated by $S=\{\rho$,$\tau\}$ and define $G_g$ for each of
	$g=\rho,\tau,\rho^{-1},\tau^{-1}$ as the set of all elements from $G$ represented by
	reduced words in $S$ having the element $g$ as the first, or leftmost, letter. By
	\Cref{free_remark},
	we know the sets $\{e\}$, $G_\rho$, $G_{\tau}$, $G_{\rho^{-1}}$ and $G_{\tau^{-1}}$ are
	pairwise disjoint and form a partition of $G$. Here $e$ is the identity element in $G$.

	We will show that $G \sim_G G_\rho \cup G_{\rho^{-1}}$. Observe that, with the exception of
	$\rho$, all the elements in $G_{\rho}$ can be further separated by their second letter in
	their word, the one after the initial \GWRD{r}. And since the words are reduced, the second letter is
	never \GWRD{r\INV}, but must be one of \GWRD{r}, \GWRD{t} or \GWRD{t\INV}. These three sets of
	words can also be obtained from the words describing $G_\rho$, $G_\tau$ and $G_\tau^{-1}$,
	respectively, by appending the letter \GWRD{r} to the beginning of all words. That is:
	\[G_\rho=\{\rho\} \cup \rho G_{\rho} \cup \rho G_{\tau} \cup \rho G_{\tau^{-1}}\text{.}\]

	Inspired by what we did with the Hyperwebster in \Cref{ex:hyperwebster} we now remove
	the first letter, \GWRD{r}, from the words that represent these elements. We do this by applying
	$\rho^{-1}$ to the elements in the set $G_\rho$ to obtain:
	\[\rho^{-1}G_\rho=\{e\} \cup G_{\rho} \cup G_\tau \cup G_{\tau^{-1}}\text{.}\]
	And since $G=\{e\} \cup G_{\rho} \cup G_\tau \cup G_{\tau^{-1}} \cup G_{\rho^{-1}}$, we see
	that $G=\rho^{-1}G_\rho\cup G_{\rho^{-1}}$, meaning $G$ is $G$\-/equidecomposable with
	$G_\rho\cup G_{\rho^{-1}}$. Remember that $\rho^{-1}\in G$.

	Similarly, it can be shown that $G \sim_G G_\tau \cup G_{\tau^{-1}}$. Since $G$ is
	$G$\-/equidecomposable with its two disjoint subsets $G_\rho \cup G_{\rho^{-1}}$ and
	$G_\tau \cup G_{\tau^{-1}}$, \Cref{bsb_cor} now gives us that $G$ is $G$\-/paradoxical.
\end{proof}
Observe the striking resemblance to the example of the Hyperwebster in \Cref{ex:hyperwebster}. We
can consider $G$ as the entire dictionary, and for all four choices of $g$, the subsets $G_g$ can
be thought of as being the volume \WRD{g}. The only significant difference is that we take a
specific subset, $G_{\rho^{-1}}$, and leave it unchanged during the reasoning. And just like in the
example, we need at least two symbols in our alphabet: a rank of $2$ or more. Let us instead say
$G$ is freely generated by just $\set{\rho}$, then all reduced words except the empty could only
consist of one of the two letters repeated, such as \GWRD{rrr} and \GWRD{r\INV r\INV}.
Any word containing \GWRD{rr\INV} or \GWRD{r\INV r} is not reduced. A direct equivalence
in our example of the Hyperwebster would be to only use words in the single letter \WRD{A}, which, as
already stated, would not allow us to achieve the paradoxical result.

\begin{remark}
	In a more general context, what we truly need is a non-abelian free group. And as noted in
	\Cref{remark:freeresults}, all free groups of rank $2$ or higher are non-abelian.
\end{remark}

If the paradoxical group acts on a set, we can use this paradoxicality to show that the set also is paradoxical, under one condition. It is here, and only here, that we need the \nameref{ax_choice}, which is stated in \Cref{ax_choice}. We say that an action, $g$, from a group $G$ acting on a $G$-set $X$ has a \textit{fixed point} in a $G$-set $X$ if some element in $X$ satisfies $gx=x$. Recall from \Cref{def:groupaction} that an \textit{orbit} of an element $x\in X$ is the set of all elements it gets mapped to by all actions of $G$, written $Gx$. Observe that any two orbits either share all points or are disjoint, so the collection of all unique orbits is a partition of $X$.

\begin{prop}
	\label{G_7}
	Let $G$ be a $G$\-/paradoxical group, acting on a set $X$. If only the identity in $G$ has
	any fixed points in $X$, then $X$ is also $G$\-/paradoxical.
\end{prop}
Informally, we might consider this to \textit{lift} the paradoxical decomposition of $G$ to the $G$-set. The proof essentially creates an isomorphism from the decomposition and actions between the group and the $G$-set.
\begin{proof}
	Let $G$ and $X$ be as stated and consider the orbits of the elements in $X$. From the collection
	of all orbits of the elements in $X$, we use the \nameref{ax_choice}, \Cref{ax_choice}, to obtain a set,
	$C$, containing one element from each unique orbit.

	We see that applying all group actions of $G$ on an element from $C$ results in all elements
	in that element's orbit. Thus applying all group actions on all elements of $C$ gives us the original $X$:
	\[GC=\set{gc | g\in G, c\in C}=\set{Gc | c\in C}=X.\]

	Since $G$ is $G$\-/paradoxical, there exists a paradoxical decomposition $G=A\cup B$ such that
	$G \sim_G A$ and $G \sim_G B$. Just as $X=GC$, we can apply the group actions of $A$ and $B$ on $C$ to
	obtain a partition of $X$ by $Y\coloneqq AC$ and $Z\coloneqq BC$. We see that $X=Y\cup Z$ directly from
	$X=GC=(A\cup B)C=AC\cup BC=Y\cup Z$, and disjointness follows from requirements on fixed
	points in $X$ and the construction of $C$: assume there is a common element in $Y$ and $Z$, meaning there are $a\in
	A$, $b\in B$ and $c_1,c_2\in C$ such that $ac_1=bc_2$. But then $c_1=a^{-1}bc_2$, so both
	$c_1$ and $c_2$ belong to the same orbit. Since $A$ and $B$ are disjoint $a\ne b$ and
	$a^{-1}b$ is not the identity. Combined with our assumption on $X$, this means $c_1\ne c_2$ because
	otherwise it would be a fixed point for $a^{-1}b$. Thus $C$ contains more than one element from
	an orbit, which is a contradiction.

	We wish to show that $Y$ and $Z$ form a paradoxical decomposition of $X$. It only remains to check
	that $X \sim_G Y$ and $X \sim_G Z$: since $G\sim_G A$ the definition of
	paradoxicality means there are decompositions $G=\bigcup_{i=1}^n G_i$ and
	$A=\bigcup_{i=1}^n A_i$ with elements $g_1,g_2,\dots,g_n\in G$ such that $A_i=g_iG_i$ for
	$i=1,2,\dots,n$. Using this, we define decompositions $X=\bigcup_{i=1}^nX_i$ and
	$Y=\bigcup_{i=1}^nY_i$ by $X_i\coloneqq G_iC$ and $Y_i\coloneqq A_iC$.
	That these are decompositions follows from the same reasoning as for $X=Y\cup Z$ above and
	keeping in mind that $X=GC$ and $Y=AC$.

	It's now easy to show that $Y_i=g_iX_i$ for all $i=1,2,\dots,n$. Let $i$ be any of these
	integers. Since $A_i=g_iG_i$, we get
	\begin{align*}
		Y_i=A_iC =(g_iG_i)C	& =\set{hc|h\in g_iG_i,c\in C} \\
					& =\set{g_ig'c|g'\in G_i,c\in C} \\
					& =\set{g_ix | x\in G_iC} =g_i(G_iC)=g_iX_i,
	\end{align*}
	and we have shown the $G$\-/equidecomposability of $Y$ and $A$ using the same actions
	from $G$ that gives $G\sim_G A$. $X \sim_G Z$ follows similarly and we see that $X$ is
	$G$\-/paradoxical.
\end{proof}

Observe that the paradoxicality of $X$ is built from the decompositions and elements $g_i$ that describe the paradoxicality of $G$. In order to do this we the \nameref{ax_choice}, \Cref{ax_choice}. In many cases, the collection of orbits is uncountable, which means we need the full strength of the axiom since no weaker version is sufficient.

Once again, we see the close resemblance to Vitali's construction in \Cref{ex:vitali}, both in the usage of the \nameref{ax_choice} to pick one element from each unique orbit and how group actions are applied to this set to restore the original orbits. This is also where the key difference appears: we use the properties of the group to put these subsets back together into two paradoxical subsets and obtain finite paradoxicality instead of the infinite paradoxicality of Vitali's approach.


We are now ready to prove the Hausdorff paradox. Recall that it says:
\hausdorffprdx*

Based on the results of \Cref{G_5,G_7} our goal is to find a subgroup of $SO(3)$ free on two generators, from which \Cref{G_5} makes it paradoxical on itself. And then find a subset $D$ of $\SPH$ that will allow \Cref{G_7} to transfer the paradoxical property to $\SPHMD$ itself.

Even though \Cref{G_7} requires the \nameref{ax_choice}, the result from \Cref{G_5} already gives us a paradoxical result without the use of the \nameref{ax_choice}. Still, the latter is often blamed for the paradoxical outcomes we are about to see.

We begin with finding the free subgroup of $SO(3)$. This forms the core of the Hausdorff paradox, which itself is the core of the Banach-Tarski paradox.
\begin{theorem}[Existence of a free rotation group]
	\label{G_6}
	The group $SO(3)$ has a subgroup that is free on two generators.
\end{theorem}

The following is by far the single biggest proof in this thesis and arguably the most important. It ties into the more extensive study of free groups. This proof will utilize modular arithmetics, a later approach not seen in the original proof. It is unclear where it first originated, but the idea has been expressed in various forms in \cite{kursbok}, \cite{taoart} and the later edition of \cite{wagonbok}. The version presented here also utilizes a Pythagorean triple, which lends itself to a simple generalization in \Cref{remark:pythagorean}. This approach does not seem to have been previously documented.



Remember the example of irrational rotation in \Cref{ex:irrational}. As observed in
\Cref{ex:irrationalgen}, we know that a single rotation operator by any \textit{irrational angle},
an irrational multiple of $2\pi$, will give us a free generator. And a rational angle would not. Unfortunately, this method and result do not generalize cleanly to multiple rotation operators in $\RE^3$, but it motivates us to limit our focus to irrational angles.

However, it is worth noting that the result we seek \textit{can} be obtained using any irrational angle, and surprisingly any angle, depending on the axes of rotation. In fact, Hausdorff's original proof used rotations of $\ang{180}$ and $\ang{120}$. However, all of these results require more complicated reasoning. See \Cref{remark:anyangle} for more details.

\begin{proof}[Proof of \Cref{G_6}]
Let us pick two rotations around irrational angles. To make the reasoning as simple as possible, we
	pick two perpendicular axes to rotate around: one along the $x$-axis and the other along
	the $y$-axis. For further simplicity, we use the same rotation angle for both.
	Define $\sigma$ and $\tau$ as rotations around the $x$-axis and $y$-axis by
	$\theta=\arcsin(3/5)=\arccos(4/5)$. This angle is taken from the right-angled
	triangle of sides $3$, $4$ and $5$, in \Cref{fig:triangle}, from the well-known Pythagorean triple. We can observe that this must be an irrational angle, for instance, by Niven's theorem, which shows this for any angle like this from a Pythagorean triple except for the right angle. This choice gives us the following rotation matrices from the well-known formulas for rotation matrices:

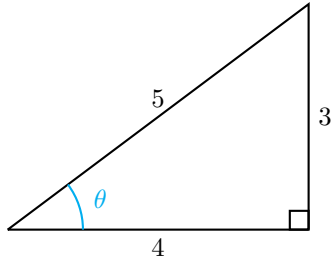
\begin{figure}[h]
    \centering
	\scalebox{1}{\input{figs/triangle.tex}}
	\caption{Angle $\theta$ from the primitive Pythagorean triangle.}
    \label{fig:triangle}
\end{figure}

\[
	\sigma=\begin{bmatrix}
		1 & 0 & 0 \\
		0 & \cos(\theta) & -\sin(\theta) \\
		0 & \sin(\theta) & \cos(\theta)
	\end{bmatrix}
	=\begin{bmatrix}
		5 & 0 & 0 \\
		0 & \frac{4}{5} & \frac{-3}{5} \\
		0 & \frac{3}{5} & \frac{4}{5}
	\end{bmatrix}
	=\frac{1}{5}\begin{bmatrix}
		5 & 0 & 0 \\
		0 & 4 & -3 \\
		0 & 3 & 4
	\end{bmatrix}
	\text{,}
\]
and similarly:
\[
	\tau=\begin{bmatrix}
		\cos(\theta) & 0 & \sin(\theta) \\
		0 & 1 & 0 \\
		-\sin(\theta) & 0 & \cos(\theta)
	\end{bmatrix}
	=\frac{1}{5}
	\begin{bmatrix}
		4 & 0 & 3 \\
		0 & 5 & 0 \\
		-3 & 0 & 4
	\end{bmatrix}
	\text{.}
\]
Their inverses, their counter-rotations, being:
\[
	\sigma^{-1}=\frac{1}{5}
	\begin{bmatrix}
		5 & 0 & 0 \\
		0 & 4 & 3 \\
		0 & -3 & 4
	\end{bmatrix}
	\quad\text{and}\quad
	\tau^{-1}=\frac{1}{5}
	\begin{bmatrix}
		4 & 0 & -3 \\
		0 & 5 & 0 \\
		3 & 0 & 4
	\end{bmatrix}
	\text{.}
\]

Note that the choice of $\theta$ and axes gives easy-to-read matrices. All
matrices have rational elements, so we can also consider them as mappings $\Q^3\to\Q^3$, and the only
denominator that occurs in their elements is $5$. Instead of matrices, we can describe them as
the following linear mappings:
\begin{align}
\begin{split}
	\sigma(x,y,z)&=1/5\cdot(5x,4y-3z,3y+4z) \text{,}\\
	\sigma^{-1}(x,y,z)&=1/5\cdot(5x,4y+3z,-3y+4z) \text{,}\\
	\tau(x,y,z)&=1/5\cdot(4x+3z,5y,-3x+4z) \text{,}\\
	\tau^{-1}(x,y,z)&=1/5\cdot(4x-3z,5y,3x+4z) \text{.}
	\label{eq:0}
\end{split}
\end{align}

\begin{figure}[h]
	\centering
	\def\svgscale{0.8}
	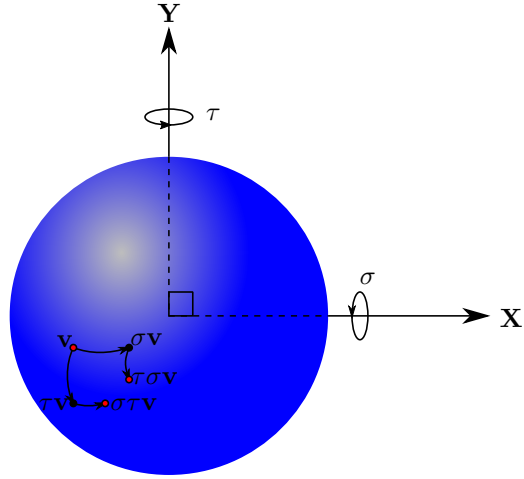
	\caption{$\sigma\tau$ and $\tau\sigma$ move $\vect{v}$ to two different
	points, and more generally, no two different reduced words represent the same rotation. As
	in \Cref{free_remark}, this is equivalent to no nonempty reduced word representing
	the identity, which keeps all points stationary.}
	\label{fig:hausdorff_ball}
\end{figure}

To show that these rotations freely generate a subset of $SO(3)$, we want to show that none of
their compositions, represented by a nonempty reduced word in $\set{\sigma,\tau}$,
represents the identity mapping, as in \Cref{def:free_group}. \Cref{fig:hausdorff_ball} illustrates
an equivalent perspective.
Note that we will consider our mappings to be $\Q^3\to\Q^3$, and work
with vectors in $\Q^3$ instead of $\RE^3$. Remember that the
identity maps all vectors to themselves, so for any other mapping there exists some vector that does not get mapped to
itself.

We will show that for any mapping represented by a nonempty reduced word in
$\set{\sigma,\tau}$, there is a rational vector that does not get mapped to itself. For
consistency with later parts of the proof, we will avoid the terminology of words:

\begin{claim}
	\label{claim1}


	For any mapping $\rho:\Q^3\to\Q^3$ by a finite composition of mappings
	from \Cref{eq:0}, not containing
	$\sigma\sigma^{-1}$, $\sigma^{-1}\sigma$, $\tau\tau^{-1}$ or $\tau^{-1}\tau$,
	there is a $\vect{v}\in\Q^3$ such that $\rho\vect{v}\ne\vect{v}$.
\end{claim}

By proving \Cref{claim1} we obtain a proof of \Cref{G_6}.

Since all mappings in \Cref{claim1} are rotations around origo, we observe that they do not change
the norm of vectors, only their directions. Any two nonzero vectors with the same direction will
still share the same direction after transformation: just observe that a mapping $\rho$ as in
\Cref{claim1} is linear and that $\rho(k\vect{v})=k\rho\vect{v}$ for all $\vect{v}\in\Q^3$ and $k\in
Q$. This motivates us to make some simplifications to work with integers
instead of rationals:

First, observe that for every vector in $\Q^3$, there is a vector in $\Z^3$ with the same direction:
for any $\vect{v}\in\Q^3$, let $k\in\Z$ be the product of the denominators in the components of
$\vect{v}$, or their least common multiple, for example. Then $k\vect{v}$ has integer components.
For any $\rho$ as in \Cref{claim1}, $\rho(k\vect{v})=k\rho\vect{v}=k\vect{v}$ if and only if
$\rho\vect{v}=\vect{v}$, since $k\ne0$. Thus we can limit ourselves to only vectors in $Z^3$ when
finding a $\vect{v}$ for \Cref{claim1}.

Furthermore, not all mappings $\rho$ will transform vectors from $\Z^3$ into vectors in $\Z^3$
since the mappings in \Cref{eq:0} contain rational coefficients, some of which have denominators
of $5$. However, similarly to how we scaled the vectors, the mappings can also,
in a way, be scaled: we define these mappings from them by, informally, multiplying them by $5$
and observing that the new mappings only have integer coefficients. Thus they map vectors from
$\Z^3$ to $\Z^3$.
We define these new $\Z^3\to\Z^3$ mappings as:
\begin{align}
\begin{split}
	s_+(a,b,c)&\coloneqq(5a,4b-3c,3b+4c) \text{,}\\
	s_-(a,b,c)&\coloneqq(5a,4b+3c,-3b+4c) \text{,}\\
	t_+(a,b,c)&\coloneqq(4a+3c,5b,-3a+4c) \text{,}\\
	t_-(a,b,c)&\coloneqq(4a-3c,5b,3a+4c) \text{.}
	\label{eq:1}
\end{split}
\end{align}


From any mapping $\rho:\Q^3\to\Q^3$ as in \Cref{claim1}, we can obtain a mapping
$\rho':\Z^3\to\Z^3$ by 
replacing all occurrences of mappings $\sigma$, $\sigma^{-1}$, $\tau$ and $\tau^{-1}$ in the
composition with $s_+$, $s_-$, $t_+$ and $t_-$, respectively.

The new mapping satisfies $\rho'\vect{v}=5^n\rho\vect{v}$, for any $\vect{v}\in\Z^3$, where $n\ge1$
is the number of mappings in the composition; remember that each mapping in \Cref{eq:1}
increases the vector's norm by $5$.

We combine these two simplifications: if there is a $\rho$ as in \Cref{claim1} and
a $\vect{v}\in\Q^3$ such that
$\rho\vect{v}=\vect{v}$, then there also exists a $\vect{v'}\in\Z^3$ and $\rho':\Z^3\to\Z^3$, as
above, such that
$\rho'\vect{v}=5^n\vect{v}$ for some $n\ge1$. Now, observe its contrapositive: given any $\rho$
with $\rho'$ obtained from it as above, we see that if there is a
$\vect{v'}\in\Z^3$ such that $\rho'\vect{v'}\ne5^k\vect{v'}$ for all $k\ge1$, then it must also be
that $\rho\vect{v'}\ne\vect{v'}$.


Since the composition for any $\rho$ of \Cref{claim1} contains none of $\sigma\sigma^{-1}$,
$\sigma^{-1}\sigma$, $\tau\tau^{-1}$ and $\tau^{-1}\tau$, we get that the composition for $\rho'$
above never contains $s_+s_-$, $s_-s_+$, $t_+t_-$ or $t_-t_+$.
Observe that mappings like $s_+$ and $s_-$ are no longer each other's inverses
because they also scale vectors:
$s_+s_-\vect{v}=5\sigma5\sigma^{-1}\vect{v}=5^2\vect{v}$, for any $\vect{v}\in\Z^3$. Nevertheless,
they can still be considered in the same spirit, and avoiding these pairs in the composition will
be vital.

We can thus simplify our reasoning in \Cref{claim1} from $\RE^3$, or rather $\Q^3$, to $\Z^3$,
using only vectors with integer components and considering mappings from \Cref{eq:1} instead of
\Cref{eq:0}. This turns \Cref{claim1} into:

\begin{claim}
	\label{claim2}


	For any mapping $\rho:\Z^3\to\Z^3$ by a finite composition of the mappings
	from \Cref{eq:1}, not containing
	$s_+s_-$, $s_-s_+$, $t_+t_-$ or $t_-t_+$,
	there is a $\vect{v}\in\Z^3$ such that $\rho\vect{v}\ne5^k\vect{v}$ for all $k\ge1$.
\end{claim}
If we can prove \Cref{claim2}, then \Cref{claim1} follows. Note that $k\ge1$ and never $0$. This is because $\rho$ consists of at least one transformation.

So far, we have simplified our space without changing what we can or cannot prove, but we can make
one more significant simplification by considering all integers in modular form; this includes both
scalars and each component in the vectors. Since the condition in \Cref{claim2} contains a
multiple of $5$, it would be practical to work modulo $5$: observe that $5^k\vect{v}\equiv_5\vect{0}
\mod5$ for all $k\ge1$ and $\vect{v}\in\Z^3$, where $\vect{0}=(0,0,0)$.

Some notation: we will use vectors from $\Z_5^3$ and indicate congruence modulo $5$ by $a\equiv_5
b$ as a shorter form for $a \equiv b \mod 5$, where $a$ and $b$ can be either scalars or vectors. Two
vectors are considered congruent component-wise: if each component of one vector is congruent to
the matching component of the other. Our existing linear maps from \Cref{eq:1}
can be used within this context, considered as mappings $\Z_5^3\to\Z_5^3$, instead of
$\Z^3\to\Z^3$:
\begin{align}
\begin{split}
	s'_+(a,b,c)&\equiv_5(0,4b+2c,3b+4c) \text{,}\\
	s'_-(a,b,c)&\equiv_5(0,4b+3c,2b+4c) \text{,}\\
	t'_+(a,b,c)&\equiv_5(4a+3c,0,2a+4c) \text{,}\\
	t'_-(a,b,c)&\equiv_5(4a+2c,0,3a+4c) \text{.}
	\label{eq:2}
\end{split}
\end{align}
Here the negative coefficients from \Cref{eq:1} have been recalculated modulo $5$, and terms with
coefficients of $5$ are always congruent with $0$.

It turns out that this simplification from $\Z^3$ to $\Z_5^3$ will still allow us to prove what
we want but also make it much easier. Of course, any result made modulo $5$ will be made in a
smaller number of objects: $5^5=125$ different vectors in $\Z_5^3$ compared to the infinite number
of vectors in $\Z^3$. Each vector in $\Z_5^3$ simultaneously represents a countably
infinite number of vectors in $\Z^3$. In contrast, going from $\Q^3$ to $\Z^3$ above does not cause
these kinds of potential problems.

However, as long as we can show that our rotations do not give the identity in $\Z_5^3$, then it cannot
possibly occur in $\Z^3$ either, and we have our proof: consider a mapping $\rho$ as in
\Cref{claim2}, if $\rho\vect{v}=5^k\vect{v}$ for some $k\ge1$ and $\vect{v}\in\Z^3$, then we must
also have $\rho\vect{v}\equiv_5\vect{0}$. So if there exists a $\vect{v}\in\Z_5^3$
such that $\rho\vect{v}\not\equiv_5\vect{0}$, then $\rho\vect{v}\ne5^k\vect{v}$ for all $k\ge1$ and \Cref{claim2}
follows.

\begin{claim}
	\label{claim3}



	For any mapping $\rho:\Z_5^3\to\Z_5^3$ by a finite composition of the mappings
	from \Cref{eq:2}, not containing
	$s'_+s'_-$, $s'_-s'_+$, $t'_+t'_-$ or $t'_-t'_+$,
	there is a $\vect{v}\in\Z_5^3$ such that $\rho\vect{v}\not\equiv_5\vect{0}$.
\end{claim}

If we can prove \Cref{claim3}, then \Cref{claim2} follows, by which \Cref{claim1} and
\Cref{G_6} is proven. Now we begin observing useful properties of our mappings and prove
\Cref{claim3}.

\begin{proof}[Proof of \Cref{claim3}]
	\renewcommand\qedsymbol{$\blacksquare$}
	Let $s'_+$, $s'_-$, $t'_+$ and $t'_-$ be as in \Cref{eq:2}.
	We know $s'_+$ maps any vector $(a,b,c)\in Z_5^3$ to $(0,4b+2c,3b+4c)$. Basic algebra
	shows that $4b+2c \equiv_5 3(3b+4c)$, which means all these vectors can be described
	by $(0,3m,n)$ for all $m,n\in Z_5$. This gives us the \textit{range} of $s'_+$:
	$range(s'_+)=\set{(0,3m,n)\in Z_5^3}$.

	Of the vectors $(a,b,c)$, those that get mapped to $\vect{0}$ by $s'_+$ are precisely those
	that
	satisfy $(0,4b+2c,3b+4c)\equiv_5(0,0,0)$, which means $a$ can be of any value and
	$4b+2c\equiv_50$, or equivalently $3b+4c\equiv_50$ because $4b+2c \equiv_5
	3(3b+4c)$. From any of the congruences for $b$ and $c$, we obtain $c\equiv_53b$.
	Similar to the range, this gives us the \textit{kernel} of $s'_+$:
	$ker(s'_+)=\set{(m,n,3n)\in Z_5^3}$.
	In the same way, we calculate the range and kernel for all mappings:
\begin{align*}
\begin{split}
	range(s'_+)&=\set{(0,3m,m)\in Z_5^3}\text{ and }ker(s'_+)=\set{(m,n,3n)\in Z_5^3}\text{,} \\
	range(s'_-)&=\set{(0,m,3m)\in Z_5^3}\text{ and }ker(s'_-)=\set{(m,3n,n)\in Z_5^3}\text{,} \\
	range(t'_+)&=\set{(m,0,3m)\in Z_5^3}\text{ and }ker(t'_+)=\set{(3m,n,m)\in Z_5^3}\text{,} \\
	range(t'_-)&=\set{(3m,0,m)\in Z_5^3}\text{ and }ker(t'_-)=\set{(m,n,3m)\in Z_5^3}\text{.}
\end{split}
\end{align*}

	We see that the range of $s'_+$ is contained in the kernel of $s'_-$,
	$range(s'_+)\subset ker(s'_-)$, and vice versa, $range(s'_-)\subset ker(s'_+)$.
	This is not surprising since they represent $\sigma$ and
	$\sigma^{-1}$, with additional multiplications of $5$:
	$s'_+s'_-\vect{v}=5\sigma5\sigma^{-1}\vect{v}=5^2\vect{v}\equiv_5 \vect{0}$ for any $\vect{v}\in\Z^3$.
	And similarly for $t'_+$ and $t'_-$.
	But more importantly, no further overlappings of kernels and
	ranges occur for other pairs of mappings, except for $\vect{0}$, which
	belongs to all of them. For example: $range(s_+)\cap ker(t_-)=\set{\vect{0}}$.



	Since $\rho$ in \Cref{claim3} is given by a composition of these mappings where none of
	$s'_+s'_-$, $s'_-s'_+$, $t'_+t'_-$ or $t'_-t'_+$ occur, the only common element
	between the kernel and range of each mapping in the composition is the zero vector. If a
	vector is not in
	the kernel of the first mapping, it gets mapped to a vector that is not in the kernel of the
	next. This continues for all mappings in the finite composition until the last mapping
	gives a nonzero vector, and we would have our proof.

	All that remains is to find a nonzero vector $\vect{v}\in\Z_5^3$ that is not in the kernel
	of the first transformation, for each of the four mappings. We can easily do this, but we
	can even go one step further and find a single vector outside the kernel of all
	mappings, meaning it always gets mapped to a nonzero vector by all transformations. One
	such is $(0,0,1)$ which represents, among others, the point $(0,0,1)$ on $\SPH$ in $\RE^3$.
	To see that this vector does not fit in any of the kernels, observe that the two components
	that are congruent to zero would mean the last component also would need to be
	congruent to zero for it to match the form of any kernel. However, it is $1$.
\end{proof}
Having proven \Cref{claim3}, the earlier discussion also gives \Cref{claim2} and finally
\Cref{claim1}, from which we see that $\set{\sigma,\tau}$
freely generates a subgroup of $SO(3)$. We have proven \Cref{G_6}.
\end{proof}

\begin{remark}
A shorter but less insightful approach would have been to immediately choose the
vector $\vect{v}=(0,0,1)$ as mentioned above and, by induction on the length of the sequence, show
that any mapping of $\vect{v}$ by the sequence in \Cref{claim3} would result in a different vector.
This is closer to the approach taken in \cite[pp.422]{kursbok} and does not need any reasoning
regarding the kernels of the mappings, only their ranges. However, it offers less insight as we
	merely verify a suitable guess.
\end{remark}

\begin{remark}
	\label{remark:anyangle}
	Regarding the choice of rotation angles, Hausdorff's original proof uses angles of
	$\ang{180}$ and $\ang{120}$ degrees, around two non-orthogonal axes. It might at first
	seem to be an impossible choice since these rotations can be seen on their own to generate
	cyclic subgroups of rotations, isomorphic to $\Z_2$ and $\Z_3$, and cyclical groups
	themselves can not be free. But a clever composition of the two can create new elements
	that instead freely generate a subgroup of the rotations. This subgroup can be described by
	a \textit{free product}, $\Z_2*\Z_3$, which is isomorphic to the free group of rank $2$,
	not to be confused with the usual product of groups.
	Further investigation of this idea and relaxing of the requirements of the rotations can be
	found in \cite[pp.~15--16]{wagonbok} and, of course, the original works by Hausdorff in
	\cite{hausdorffbok} and \cite{hausdorffart}, which do not explicitly use this terminology.
\end{remark}

\begin{remark}[Generalization of the choice of Pythagorean triple]
	\label{remark:pythagorean}
	The reasoning above seems to generalize to any angle obtained from a \textit{primitive
	Pythagorean triple}, let one such be $(a,b,c)$. \textit{Primitive} means all values are
	relatively prime, and here $c$ is taken to be the length of the hypothenuse.

	The same reasoning follows, with the slight generalization towards the end when we
	determine images and kernels; for example, we need to be sure that we can find exactly one
	$k\in\Z_c$ such that $k(ax+by)\equiv_c -bx+ax$ holds for all $x,y\in Z_c$, and similarly
	for the other equations. That is $ka\equiv_c-b$ and $kb\equiv_ca$. Should none or more than
	one $k$ like this exist, then the image and kernel spaces will not have the necessary
	properties.

	However, $a^2+b^2=c^2$ and $a^2\equiv_c-b^2$. And since both of $a$ and $b$ are relatively
	prime to $c$, they have inverses. From this, we see that $ab^{-1}\equiv_c-ba^{-1}$ and let
	this be $k$. This value of $k$ satisfies both $ka\equiv_c-b$ and $kb\equiv_ca$.

	%

	After this deduction for all kernels and images, the resulting properties of disjoint and
	suitably subsets of images and kernels follows similarly to
	the original reasoning.

	It is worth noting that if $c$ is not prime, then $Z_c$ and $Z_c^3$ are not fields.
	However, this
	should not affect our calculations, as long as all values $a$, $b$, and $c$ are relatively
	prime, which is guaranteed by the use of primitive triples.

	Finally, it is curious that all Pythagorean triples can be obtained from two
	elements, in the form of matrices, that on their own also generate a free group of rank $2$,
	as noted in \cite{alperinart}.
\end{remark}

We now have all results we need to prove the Hausdorff paradox:
\hausdorffprdx*

\begin{proof}
	By \Cref{G_6}, we can obtain a subgroup $G$ of $SO(3)$ that is free on two generators. All
	elements of $G$ are rotations around origo, so when they are applied to points of $\SPH$,
	the resulting points will also be in $\SPH$. That is, $\SPH$ is \textit{closed}
	under $G$. That means we can consider $G$ to act on $\SPH$.

	We now want to apply \Cref{G_7} to $\SPH$, but $\SPH$ contains fixed points
	for elements of $G$ that are not the identity. Our goal is to construct $D$ so that $\SPHMD$
	does not contain these fixed points and remains closed under $G$, so that $\SPHMD$ can be
	seen as a $G$-set.

	From \nameref{thm:euler_rotation}, \Cref{thm:euler_rotation}, every non-identity element in $G$
	can be seen as a rotation around a line through origo. Each such line
	will intersect $\SPH$ in precisely two points: the two points
	that will stay fixed when this rotation is applied.

	Let $D$ be the collection of all these fixed points. Since all
	elements of $G$ can be represented by words in $\set{\rho,\tau}$, and these words are
	countable, $G$ is countable. Therefore $D$ must also be countable.

	It remains to show that $\SPHMD$ is closed under $G$. But if that was not the case, there would
	be a $\vect{p} \in \SPHMD$ and $\rho \in G$ such that $\rho \vect{p} \in D$. But that means $\rho
	\vect{p}$ is a fixed point for some other nontrivial rotation $\tau \in G$, so $\tau
	\rho \vect{p} = \rho \vect{p}$. Applying $\rho^{-1}$ to the left of both sides
	gives $\rho^{-1} \tau \rho \vect{p} = \vect{p}$. And since $\vect{p}$ cannot be a fixed
	point for any other element in $G$ but the trivial identity element, this means $\rho^{-1} \tau \rho = e$,
	and thus $\tau=\rho \rho^{-1} =e$, which contradicts that $\tau$ is nontrivial.

	From \Cref{G_5}, we know that the actions of $G$ on $G$ is paradoxical. And we have constructed
	$D$ such that $\SPHMD$ is closed under $G$ and hence is a $G$-set without any fixed points,
	except for the trivial identity element of $G$. We can
	now apply \Cref{G_7} to obtain that $\SPHMD$ is $G$\-/paradoxical. Since G is a subset of
	$SO(3)$, $\SPHMD$ is also $SO(3)$\-/paradoxical.
\end{proof}
Note that $D$ is not all of $\SPH$ since the latter is uncountable, so $\SPHMD$ is non-empty. The
elements in $D$ are sometimes called the \textit{poles} of the rotations since each point is kept
stationary under some rotation. This result remains valid in higher dimensions but not in lower
dimensions as we can at most find free groups of rank $1$; such as the subgroups of the rotation around origo in $\RE^2$.

In this proof, the full strength of the \nameref{ax_choice}, \Cref{ax_choice}, is needed: the number of orbits is uncountably infinite. However, this is not always the case:
\begin{remark}[Hausdorff paradox without the axiom of choice]
	\label{hausdorff_noac}
	We have seen that the group of rotations used above maps any rational point in $\SPH$ to a
	rational point in $\SPH$, which means the group also acts on the rational subset of the
	sphere, the \textit{rational sphere} $\SPH\cap\Q^3$. Thus, \Cref{G_7} can be applied to the
	rational sphere. We obtain a result that could be described as the Hausdorff paradox
	for the rational sphere: $(\SPH\cap\Q^3)\setminus D$ is $SO(3)$\-/paradoxical for some countable set $D$.

	However, there are similar groups of rotations closed in $\SPH\cap\Q^3$ that
	also have no fixed points in the set, except for the trivial identity. This means the
	entire rational sphere is $SO(3)$ paradoxical. A proof of the existence of rotations like
	this can be found in \cite[pp.~18--20]{wagonbok}, where the result is named the \textit{Satô
	paradox} in $\SPH\cap\Q^3$.

	Of the two groups of rotations mentioned, the second gives a significantly more powerful and
	exciting result: not only is $\SPH\cap\Q^3$ dense in the sphere, as shown
	in \cite[Prop.~2.7]{wagonbok}, but according to Wagon, the Satô paradox has the remarkable
	property of showing that the rational sphere is $SO(3)$\-/paradoxical without the need to
	use the \nameref{ax_choice}.
\end{remark}

We have now constructed a paradoxical subset of the sphere. The next section will build on this result by recovering the points we excluded and turn the sphere into a ball to obtain the famous result.

%% file: figs/hyperwebster.pdf_tex
\begingroup%
  \makeatletter%
  \providecommand\color[2][]{%
    \errmessage{(Inkscape) Color is used for the text in Inkscape, but the package 'color.sty' is not loaded}%
    \renewcommand\color[2][]{}%
  }%
  \providecommand\transparent[1]{%
    \errmessage{(Inkscape) Transparency is used (non-zero) for the text in Inkscape, but the package 'transparent.sty' is not loaded}%
    \renewcommand\transparent[1]{}%
  }%
  \providecommand\rotatebox[2]{#2}%
  \newcommand*\fsize{\dimexpr\f@size pt\relax}%
  \newcommand*\lineheight[1]{\fontsize{\fsize}{#1\fsize}\selectfont}%
  \ifx\svgwidth\undefined%
    \setlength{\unitlength}{557.31456806bp}%
    \ifx\svgscale\undefined%
      \relax%
    \else%
      \setlength{\unitlength}{\unitlength * \real{\svgscale}}%
    \fi%
  \else%
    \setlength{\unitlength}{\svgwidth}%
  \fi%
  \global\let\svgwidth\undefined%
  \global\let\svgscale\undefined%
  \makeatother%
  \begin{picture}(1,0.53935682)%
    \lineheight{1}%
    \setlength\tabcolsep{0pt}%
    \put(0,0){\includegraphics[width=\unitlength,page=1]{hyperwebster.pdf}}%
    \put(0.06472464,0.26881526){\color[rgb]{0,0,0}\makebox(0,0)[t]{\smash{\begin{tabular}[t]{c}$\A$\end{tabular}}}}%
    \put(0,0){\includegraphics[width=\unitlength,page=2]{hyperwebster.pdf}}%
    \put(0.92868941,0.00691584){\color[rgb]{0,0,0}\makebox(0,0)[t]{\smash{\begin{tabular}[t]{c}$\vdots$\end{tabular}}}}%
    \put(0,0){\includegraphics[width=\unitlength,page=3]{hyperwebster.pdf}}%
    \put(0.4844827,0.25632621){\color[rgb]{0,0,0}\makebox(0,0)[lt]{\smash{\begin{tabular}[t]{l}$\begin{dcases}"\_"\\A\\AA\\AB\\\vdots\\B\\BA\\BB\\\vdots\\C\\CA\\CC\\\vdots\end{dcases}$\end{tabular}}}}%
    \put(0,0){\includegraphics[width=\unitlength,page=4]{hyperwebster.pdf}}%
    \put(0.1632746,0.25469354){\color[rgb]{0,0,0}\makebox(0,0)[lt]{\smash{\begin{tabular}[t]{l}$=\begin{dcases}\xcancel{A}\\\xcancel{A}A\\\xcancel{A}B\\\vdots\\\xcancel{A}Z\\\xcancel{A}AA\\\xcancel{A}AB\\\vdots\\\xcancel{A}AZ\\\xcancel{A}BA\\\xcancel{A}BB\\\vdots\end{dcases}$\end{tabular}}}}%
    \put(0,0){\includegraphics[width=\unitlength,page=5]{hyperwebster.pdf}}%
    \put(0.75363035,0.43127233){\color[rgb]{0,0,0}\makebox(0,0)[lt]{\smash{\begin{tabular}[t]{l}$\begin{rcases}A\\AA\\AB\\\vdots\end{rcases}=$\end{tabular}}}}%
    \put(0.92599773,0.4437614){\color[rgb]{0,0,0}\makebox(0,0)[t]{\smash{\begin{tabular}[t]{c}$\A$\end{tabular}}}}%
    \put(0,0){\includegraphics[width=\unitlength,page=6]{hyperwebster.pdf}}%
    \put(0.92599773,0.188071){\color[rgb]{0,0,0}\makebox(0,0)[t]{\smash{\begin{tabular}[t]{c}$\B$\end{tabular}}}}%
    \put(0,0){\includegraphics[width=\unitlength,page=7]{hyperwebster.pdf}}%
    \put(0.75363043,0.17558186){\color[rgb]{0,0,0}\makebox(0,0)[lt]{\smash{\begin{tabular}[t]{l}$\begin{rcases}B\\BA\\BB\\\vdots\end{rcases}=$\end{tabular}}}}%
  \end{picture}%
\endgroup%

%% file: figs/triangle.tex
\begin{tikzpicture}[thick]
	\coordinate (A) at (0,0);
	\coordinate (B) at (4,0);
	\coordinate (C) at (4,3);

	\draw (A) -- node[below] {$4$} (B) -- node[right] {$3$} (C) -- node[above] {$5$} (A);
	\tkzMarkRightAngle(A,B,C)
	\tkzMarkAngle[mark=none,color=cyan](B,A,C)
	\tkzLabelAngle[pos=1.3,color=cyan](B,A,C){$\theta$}
\end{tikzpicture}

%% file: figs/hausdorff_ball.pdf_tex
\begingroup%
  \makeatletter%
  \providecommand\color[2][]{%
    \errmessage{(Inkscape) Color is used for the text in Inkscape, but the package 'color.sty' is not loaded}%
    \renewcommand\color[2][]{}%
  }%
  \providecommand\transparent[1]{%
    \errmessage{(Inkscape) Transparency is used (non-zero) for the text in Inkscape, but the package 'transparent.sty' is not loaded}%
    \renewcommand\transparent[1]{}%
  }%
  \providecommand\rotatebox[2]{#2}%
  \newcommand*\fsize{\dimexpr\f@size pt\relax}%
  \newcommand*\lineheight[1]{\fontsize{\fsize}{#1\fsize}\selectfont}%
  \ifx\svgwidth\undefined%
    \setlength{\unitlength}{253.3110946bp}%
    \ifx\svgscale\undefined%
      \relax%
    \else%
      \setlength{\unitlength}{\unitlength * \real{\svgscale}}%
    \fi%
  \else%
    \setlength{\unitlength}{\svgwidth}%
  \fi%
  \global\let\svgwidth\undefined%
  \global\let\svgscale\undefined%
  \makeatother%
  \begin{picture}(1,0.86330052)%
    \lineheight{1}%
    \setlength\tabcolsep{0pt}%
    \put(0,0){\includegraphics[width=\unitlength,page=1]{hausdorff_ball.pdf}}%
    \put(0.10916834,0.23897776){\color[rgb]{0,0,0}\makebox(0,0)[t]{\smash{\begin{tabular}[t]{c}$\vect{v}$\end{tabular}}}}%
    \put(0,0){\includegraphics[width=\unitlength,page=2]{hausdorff_ball.pdf}}%
    \put(0.27602986,0.16876485){\color[rgb]{0,0,0}\makebox(0,0)[t]{\smash{\begin{tabular}[t]{c}$\tau\sigma\vect{v}$\end{tabular}}}}%
    \put(0.23753962,0.12139227){\color[rgb]{0,0,0}\makebox(0,0)[t]{\smash{\begin{tabular}[t]{c}$\sigma\tau\vect{v}$\end{tabular}}}}%
    \put(0,0){\includegraphics[width=\unitlength,page=3]{hausdorff_ball.pdf}}%
    \put(0.38557884,0.66025532){\color[rgb]{0,0,0}\makebox(0,0)[t]{\smash{\begin{tabular}[t]{c}$\tau$\end{tabular}}}}%
    \put(0.67277521,0.35529433){\color[rgb]{0,0,0}\makebox(0,0)[t]{\smash{\begin{tabular}[t]{c}$\sigma$\end{tabular}}}}%
    \put(0.93924597,0.27535308){\color[rgb]{0,0,0}\makebox(0,0)[t]{\smash{\begin{tabular}[t]{c}\textbf{X}\end{tabular}}}}%
    \put(0.30267689,0.84382407){\color[rgb]{0,0,0}\makebox(0,0)[t]{\smash{\begin{tabular}[t]{c}$\textbf{Y}$\end{tabular}}}}%
    \put(0,0){\includegraphics[width=\unitlength,page=4]{hausdorff_ball.pdf}}%
    \put(0.26122593,0.24278453){\color[rgb]{0,0,0}\makebox(0,0)[t]{\smash{\begin{tabular}[t]{c}$\sigma\vect{v}$\end{tabular}}}}%
    \put(0.0895003,0.12139227){\color[rgb]{0,0,0}\makebox(0,0)[t]{\smash{\begin{tabular}[t]{c}$\tau\vect{v}$\end{tabular}}}}%
  \end{picture}%
\endgroup%

%% file: proof.tex
\section{Proof of the Banach-Tarski Paradox}
\label{sec:proof}

This section will present the final steps deducing the weak Banach-Tarski paradox, \Cref{b-t2},
from the Hausdorff paradox in \Cref{hausdorff}. As we saw by \Cref{equivalence_bt}, this result then
immediately gives us the strong Banach-Tarski paradox, \Cref{b-t1}.

Although Banach and Tarski are typically credited
for this result, an equivalent statement by the exact same ideas used here had already been
presented by Hausdorff ten years prior to Banach and Tarski, in \cite[pp.~469--473]{hausdorffbok} in
1914. To Hausdorff, it appears to have been a simple and an immediate result from what has now instead
become known as the Hausdorff paradox, covered in \Cref{sec:hausdorff}. See \Cref{sec:attribution}.

Our goal is to turn the paradoxicality from the subset of the sphere into paradoxicality of the
entire ball.
To give a rough outline: we will begin by obtaining the paradoxical result on the unit sphere
$\SPH$ by adding the missing points from \Cref{hausdorff} in \Cref{bt_proof_1}. Then expand the
paradoxical property on the sphere to a ball, without its center point, by extruding the subsets on
the sphere into cones in \Cref{bt_proof_2}. Finally, we add the missing point in the center of the
ball and prove \Cref{b-t2}. Fundamental to the first and last of these steps is the familiar
reasoning by irrational rotation in \Cref{lm:irrational} as a generalization of \Cref{ex:irrational}.

In the first and last of these steps, equidecomposability is shown between the different stages of
the ball's construction. We then apply \Cref{paradox_inherit} to transfer paradoxicality between
them. Alternatively, instead of \Cref{paradox_inherit}, the weaker result of
\Cref{paradox_inherit_weak} can also be used.

\subsection{Prerequisites}

Remember our irrational rotation in \Cref{ex:irrational}. Not only has it served to illustrate
various properties and ideas, but we will now create a lemma as a direct generalization which essentially allows
us to recover missing points. The generalization consists of moving from $\RE^2$ to
$\RE^3$, allowing an arbitrary axis of rotation and replacing the rotation of a single point with a,
possibly countably infinite, set of points. The only requirement is that no point stays stationary
by the rotation, which is equivalent to the no point being located along the axis of rotation.

A quick clarification on notation: here powers of $\rho$ means the number of times it, or its
inverse, are repeated. That is $\rho^0=e$, $\rho^n=\rho\rho\dots\rho$ and
$\rho^{-n}=\rho^{-1}\rho^{-1}\dots\rho^{-1}$.

\begin{lemma}[Generalized irrational rotation]
	\label{lm:irrational}
	Let the set $D\subset\RE^3$ be finite or countably infinite. If there is a line $L$ not passing
	through any point in $D$, then there exists a rotation $\rho\in G_3$ with $L$ as its axis,
	such that for $A\coloneqq\bigcup_{n=0}^\infty \rho^n D$ we have:
	\[\rho A=A\setminus D\text{.}\]
	If $L$ passes through origo, then $\rho$ also belongs to $SO(3)$.
\end{lemma}
Note that if the line does not pass through origo, we need the larger group $G_3$, which also
includes translations and thus allows for rotations
around any point. If $D$ consists of a single point, we could use any irrational angle for $\rho$,
and the result is almost identical to \Cref{ex:irrational}.

\begin{proof}
	Let $\rho_\theta$ be a rotation around $L$ by the angle $\theta$, where $\rho_\theta\in SO(3)$
	if $L$ passes through origo, otherwise $\rho_\theta\in G_3$. Consider the values
	of $\theta\in(0,2\pi]$ when $\rho_\theta$ is applied to $D$:

	First, we observe that only a countable number of angles can give $\rho_\theta D \cap D \ne
	\emptyset$. That is, any of the points in $D$ remain a point in $D$ after the
	rotation. This follows from $D$ being at most countable: for each point in $D$, there can
	only be a countable
	number of angles that make $\rho_\theta$ map the point to any point in
	$D$. Repeating this for each point in $D$, we get a set of angles like these
	for each point. We have at most a countable collection of
	sets, each containing a countable number of angles. Since the union of a countable
	collection of countable sets is countable, these angles are together at most
	countable.

	This result can be extended to repeated rotations: only a countable number of
	angles give $\rho_\theta^nD\cap D\ne\emptyset$, for some $n=1,2,\dots$. By repeating
	the reasoning above with $\rho_\theta^n$ instead of $\rho$, for every $n=1,2,\dots$, we get
	a countable collection of sets containing at most a countable number of angles. The union
	of all these sets is countable.

	Since we have an uncountable number of angles to choose from, we can choose $\theta$ to be
	different from all angles deduced above and define $\rho\coloneqq\rho_\theta$ for that
	value of $\theta$. This means $\rho^n D\cap D\ne\emptyset$ for all $n=1,2,\dots$. We also note that $\rho^nD \cap \rho^{n+m}D = \emptyset$ for all
	$n=0,1,\dots$
	and $m=1,2,\dots$; otherwise, we can apply $\rho^{-n}$ to both sides to obtain $D \cap
	\rho^mD \ne \emptyset$, contradicting our choice of $\theta$. Thus the sets in
	$\set{\rho^n D}_{n=1}^\infty$ are disjoint.

	Let $A\coloneqq\bigcup_{n=0}^\infty \rho^n D$, and by the disjointness noted above, we
	obtain:
	\[
		\rho A=\rho(\bigcup_{n=0}^\infty \rho^{n})=\bigcup_{n=0}^\infty \rho^{n+1} D
		=\bigcup_{n=1}^\infty \rho^n D=(\bigcup_{n=0}^\infty \rho^n D)\setminus D
		=A\setminus D\text{.}
	\]
\end{proof}

This result provides an equidecomposability property of a larger set containing $D$ and its
rotations, which will be vital for some of the remaining proofs:
\begin{cor}[Equidecomposability by irrational rotation]
	\label{cor:irrational}
	Let $X$ be a subset of $R^3$ with a finite or countably infinite subset $D$. If there is a
	line $L$ not passing through any point in $D$, and $D$ remains in $X$ after any
	rotation around $L$, 
	then $X\setminus D$ is $G_3$\-/equidecomposable
	with $X$. Should $L$ pass through origo, then $X\setminus D$ is also
	$SO(3)$\-/equidecomposable with $X$.
\end{cor}
See \Cref{fig:ball1} for a visualization of this kind of rotation, where $X$ is the unit sphere
$\SPH$. Used with \Cref{paradox_inherit}, this effectively allows us to add certain
additional points to a paradoxical set to obtain a larger set that remains paradoxical.

\begin{proof}
	Let $X$, $D$ and $L$ be as stated. By \Cref{lm:irrational}, there is a rotation $\rho$
	around $L$ such that for $A\coloneqq\bigcup_{n=0}^\infty \rho^n D$, we have
	$\rho A=A\setminus D$. Since $\rho^n D\subset X$ for all $n=1,2,\dots$ it follows
	that $A\subset X$.

	Now let $B\coloneqq X\setminus A$ so that $A$ and $B$ partition $X$. Since $\rho
	A=A\setminus D$ and $B\cap D=\emptyset$ the decomposition $\rho A \cup B = A\setminus D$ follows, and we
	have shown the equidecomposability. If $L$ passes through origo $\rho\in SO(3)$ they are
	$SO(3)$\-/equidecomposable; otherwise, they must be $G_3$\-/equidecomposable.
\end{proof}

$SO(3)$\-/equidecomposability also means $G_3$\-/equidecomposability, since $SO(3)\subset G_3$.

\subsection{Formal Proof}

We now have everything we need to prove the Banach-Tarski paradox.

\begin{figure}[h]
	\centering
	\def\svgscale{0.7}
	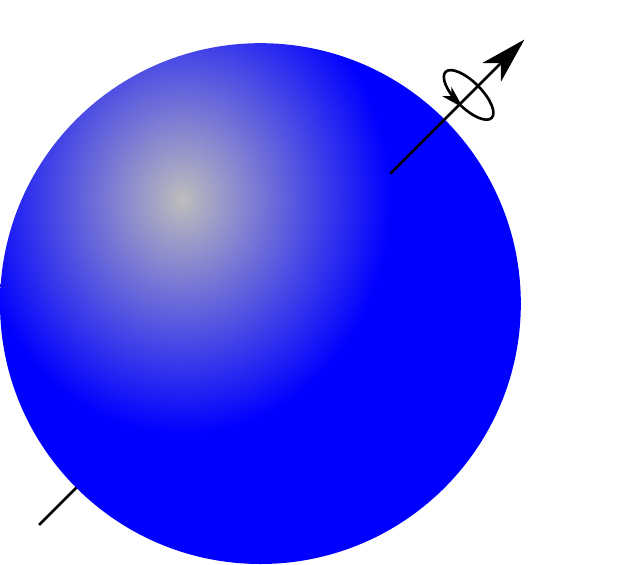
	\caption{The set $D$ of countably infinite points in the sphere $\SPH$ gets rotated around
	$L$ and remains in $\SPH$.}
	\label{fig:ball1}
\end{figure}

\begin{prop}
	\label{bt_proof_1}
	$\SPH$, the unit sphere centered at origo in $\RE^3$,
	is $SO(3)$\-/paradoxical.

\end{prop}
\begin{proof}
	Let $\SPH$ and $D$ be as in \Cref{hausdorff}: $\SPHMD$ is $SO(3)$\-/paradoxical. Select a
	line
	$L$ passing through origo but not intersecting any point in $D$. Note that there is an
	uncountable number of choices for lines passing
	through the center of $S$. But only a countable number of these will pass through any point
	in
	$D$ since $D$ is countable. Thus such a line $L$ can always be found.

	Since rotations of $D$ around $L$ does not change the distance between the points and
	origo, they will remain inside $\SPH$. See \Cref{fig:ball1}. By applying
	\Cref{cor:irrational} on $S$, $D$ and $L$ we have that $\SPHMD$ and $\SPH$ are
	$SO(3)$\-/equidecomposable. Since $\SPHMD$ is $SO(3)$\-/paradoxical, by \Cref{hausdorff}, we
	can now apply \Cref{paradox_inherit} to see that $\SPH$ also is $SO(3)$\-/paradoxical.
\end{proof}

\begin{prop}
	\label{bt_proof_2}
	Let $\BLL$ be the closed unit ball centered at origo in $\RE^3$.
	The ball with its center removed, $\BLLMC$, is $SO(3)$\-/paradoxical.
\end{prop}

The proof is based on the idea of turning the decompositions and actions from subsets on
the sphere to subsets of the ball in an isomorphic way. This informally \textit{lifts} the
paradoxicality from the sphere to the ball.

\begin{proof}
	From \Cref{bt_proof_1}, we know that $\SPH$ is $SO(3)$\-/paradoxical. For each subset $E
	\subset \SPH$, define the following conical shape:
	\[
		c(E)\coloneqq\set{t\vect{x}| \vect{x}\in E,t\in(0,1]}\text{,}
	\]
	The resulting mapping is illustrated in \Cref{fig:ball2}. Observe that
	\[ \BLLMC=\set{\vect{x}\in\RE^3 | 0<\norm{\vect{x}}\le1} \]
	So $c(\SPH)=\BLLMC$ as illustrated in \Cref{fig:ball_mapped} and $c(E) \subseteq \BLLMC$ for all $E
	\subseteq \SPH$.

	For any two disjoint subsets $A,B \subset \SPH$, their images by $c$ are also disjoint:
	$c(A) \cap c(B) = \emptyset$, otherwise they would share some point on $\SPH$ and not be
	disjoint, $A \cap B \ne \emptyset$. See \Cref{fig:ball3}. This generalizes to any finite
	number of pairwise disjoint sets.

\begin{figure}
	\centering
	\begin{subfigure}{0.3\textwidth}
		\def\svgscale{0.7}
		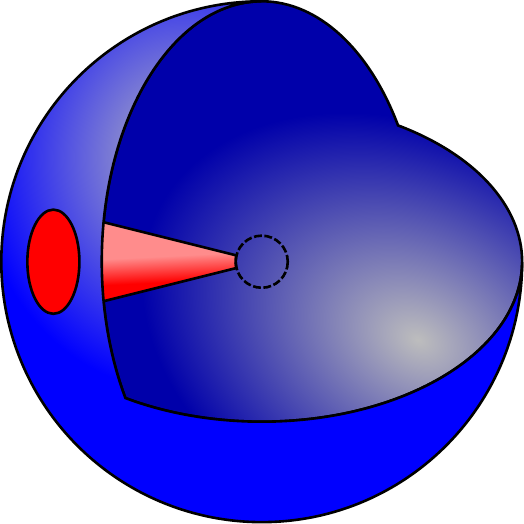
		\caption{The set $E$ in $\SPH$ gets mapped to the conical shape $c(E)$ without
		the point in origo.}
		\label{fig:ball2}
	\end{subfigure}
	\quad
	\begin{subfigure}{0.3\textwidth}
		\def\svgscale{0.7}
		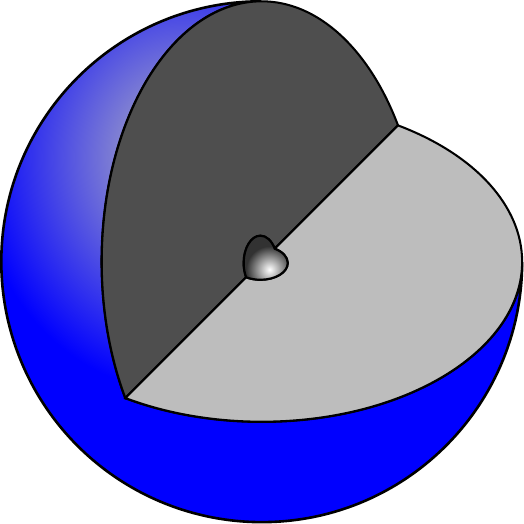
		\caption{$\SPH$ gets mapped to the unit ball without the point in origo,
		$c(\SPH)=\BLLMC$.}
		\label{fig:ball_mapped}
	\end{subfigure}
	\quad
	\begin{subfigure}{0.3\textwidth}
		\def\svgscale{0.7}
		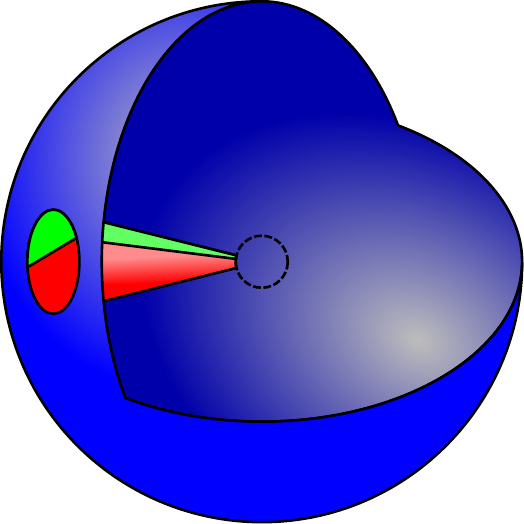
		\caption{The sets $A$ and $B$ get mapped to disjoint conical shapes $c(A)$ and
		$c(B)$.}
		\label{fig:ball3}
	\end{subfigure}
	\caption{Construction of conical shapes inside $\BLLMC$ from subsets of $\SPH$, the surface
	of $\BLL$.}
\end{figure}

	Since $\SPH$ is paradoxical, let $\SPH=C \cup D$ be a paradoxical decomposition into two
	disjoint subsets that are $SO(3)$\-/equidecomposable with $\SPH$. Note that this makes
	$c(\SPH)=c(C\cup D)=c(C)\cup c(D)$ a decomposition of $\BLLMC$. We will show that these
	subsets are $SO(3)$\-/equidecomposable with $c(\SPH)$, making $\BLLMC$ $SO(3)$\-/paradoxical:

	From the equidecomposability of $\SPH$ and $C$ there is a bijection $f:\SPH \to C$ defined
	piecewise in $SO(3)$, as in \Cref{equidecomp_function}. And using it we can define a
	mapping $g:c(\SPH) \to c(C)$ by extending the actions of $f$ to work on the conical shapes:
	\[
		g(\vect{x})=\|\vect{x}\|\cdot f(\vect{x}/\|\vect{x}\|),
	\]
	which remains defined piecewise in the same group actions used in $f$.
	$g$ is bijective: it is surjective since
	\[g(c(\SPH))=\set{tf(\vect{x})|\vect{x}\in\SPH, 0< t \le 1}=\set{t\vect{y}|\vect{y}\in C,
	0< t \le 1}=c(C)\]
	and injectivity follows by observing that if $g(\vect{x})=g(\vect{y})$ for any nonzero
	$\vect{x}$ and $\vect{y}$ then $\|\vect{x}\|=\|\vect{y}\|$ and
	$f(\vect{x}/\|\vect{x}\|)=f(\vect{y}/\|\vect{x}\|)$ from which
	$\vect{x}/\|\vect{x}\|=\vect{y}/\|\vect{x}\|$ and $\vect{x}=\vect{y}$ follows.

	Thus by
	\Cref{equidecomp_function} $g$ shows that $c(\SPH)=\BLLMC$ is equidecomposable with
	$c(C)$. The same reasoning also applies to $\BLLMC$ and $c(D)$.

	Since $\BLLMC$ is equidecomposable with both $c(C)$ and $c(D)$ whose union is
	$c(C)\cup c(D)=c(\SPH)=\BLLMC$, we see that $\BLLMC$ is paradoxical by definition.
\end{proof}

\begin{figure}[h]
	\centering
	\begin{subfigure}{0.3\textwidth}
		\def\svgscale{0.7}
		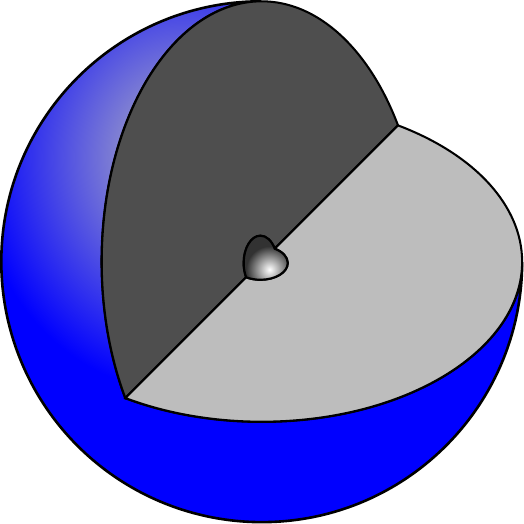
		\caption{The unit ball without the point in its center, $\BLLMC$.}
		\label{fig:ball4}
	\end{subfigure}
	\quad
	\begin{subfigure}{0.3\textwidth}
		\def\svgscale{0.7}
		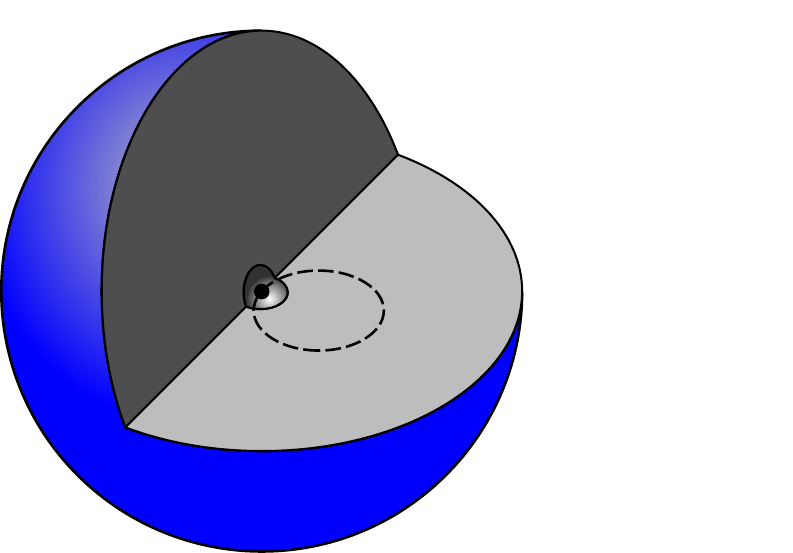
		\caption{$\vect{0}$ remains inside $\BLL$ after any rotation around $L$.}
		\label{fig:ball5}
	\end{subfigure}
	\quad
	\begin{subfigure}{0.3\textwidth}
		\def\svgscale{0.7}
		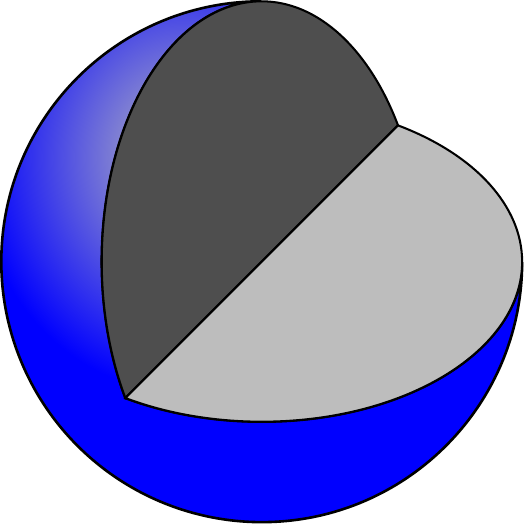
		\caption{The solid ball $\BLL$ after the center has been added.}
		\label{fig:ball6}
	\end{subfigure}
	\caption{The process of adding the missing point to $\BLLMC$.}
\end{figure}

Finally, using the same steps as \Cref{bt_proof_1}, the point in the center of the ball can also be
absorbed through equidecomposability, this time a simpler application of \Cref{lm:irrational}, which is almost identical to the example of irrational rotation in \Cref{ex:irrational}, but
applied in $R^3$. This final detail was omitted by Hausdorff, most likely since his interest was in measurability, or lack thereof, and a single point would have zero measure.
\btprdxweak*
\begin{proof}
	From \Cref{bt_proof_2}, we know that $\BLLMC$ in \Cref{fig:ball4}.is $SO(3)$\-/paradoxical.
	Let $L$ be a line that passes within the distance $1/2$ of origo, the
	center of $\BLL$, but without intersecting it. For example, we can let $L$ be the line along the
	$z$-axis intersecting the point $(1/3,0,0)$. Let $D$ be a set containing only
	the point in origo; $D=\{\vect{0}\}$.

	Note that all rotations of $D$ around $L$ move the point along a circle of
	radius less than $1/2$ centered within the distance $1/2$ from origo. Thus $D$ remains a
	subset of $\BLL$ for all rotations. See \Cref{fig:ball5}. These are rotations in $G_3$ and
	not $SO(3)$ since $L$ does not pass through origo. By applying \Cref{cor:irrational} on
	$\BLL$,
	$D$ and $L$, we obtain that $\BLL$ and $\BLLMC$ are $G_3$\-/equidecomposable.

	Since $\BLLMC$ is $SO(3)$\-/paradoxical, it is also $G_3$\-/paradoxical as $SO(3)$ is a
	subgroup of $G_3$. From \Cref{paradox_inherit}, we see that $\BLL$ is also
	$G_3$\-/paradoxical. See \Cref{fig:ball6}.
\end{proof}

We have now proven the Banach-Tarski paradox. From \Cref{equivalence_bt}, we also get its strong
form:
\btprdx*

%% file: figs/ball1.pdf_tex
\begingroup%
  \makeatletter%
  \providecommand\color[2][]{%
    \errmessage{(Inkscape) Color is used for the text in Inkscape, but the package 'color.sty' is not loaded}%
    \renewcommand\color[2][]{}%
  }%
  \providecommand\transparent[1]{%
    \errmessage{(Inkscape) Transparency is used (non-zero) for the text in Inkscape, but the package 'transparent.sty' is not loaded}%
    \renewcommand\transparent[1]{}%
  }%
  \providecommand\rotatebox[2]{#2}%
  \newcommand*\fsize{\dimexpr\f@size pt\relax}%
  \newcommand*\lineheight[1]{\fontsize{\fsize}{#1\fsize}\selectfont}%
  \ifx\svgwidth\undefined%
    \setlength{\unitlength}{182.38631407bp}%
    \ifx\svgscale\undefined%
      \relax%
    \else%
      \setlength{\unitlength}{\unitlength * \real{\svgscale}}%
    \fi%
  \else%
    \setlength{\unitlength}{\svgwidth}%
  \fi%
  \global\let\svgwidth\undefined%
  \global\let\svgscale\undefined%
  \makeatother%
  \begin{picture}(1,0.89060221)%
    \lineheight{1}%
    \setlength\tabcolsep{0pt}%
    \put(0,0){\includegraphics[width=\unitlength,page=1]{ball1.pdf}}%
    \put(0.88822474,0.47289754){\color[rgb]{1,0,0}\makebox(0,0)[t]{\smash{\begin{tabular}[t]{c}$D$\end{tabular}}}}%
    \put(0.92934641,0.34953302){\color[rgb]{0,1,0}\makebox(0,0)[t]{\smash{\begin{tabular}[t]{c}$\rho D$\end{tabular}}}}%
    \put(0,0){\includegraphics[width=\unitlength,page=2]{ball1.pdf}}%
    \put(0.84299039,0.69906599){\color[rgb]{0,0,0}\makebox(0,0)[t]{\smash{\begin{tabular}[t]{c}$\rho$\end{tabular}}}}%
    \put(0.39065398,0.86355193){\color[rgb]{0,0,0}\makebox(0,0)[t]{\smash{\begin{tabular}[t]{c}$\SPH$\end{tabular}}}}%
    \put(0.82242979,0.8429911){\color[rgb]{0,0,0}\makebox(0,0)[t]{\smash{\begin{tabular}[t]{c}$L$\end{tabular}}}}%
  \end{picture}%
\endgroup%

%% file: figs/ball2.pdf_tex
\begingroup%
  \makeatletter%
  \providecommand\color[2][]{%
    \errmessage{(Inkscape) Color is used for the text in Inkscape, but the package 'color.sty' is not loaded}%
    \renewcommand\color[2][]{}%
  }%
  \providecommand\transparent[1]{%
    \errmessage{(Inkscape) Transparency is used (non-zero) for the text in Inkscape, but the package 'transparent.sty' is not loaded}%
    \renewcommand\transparent[1]{}%
  }%
  \providecommand\rotatebox[2]{#2}%
  \newcommand*\fsize{\dimexpr\f@size pt\relax}%
  \newcommand*\lineheight[1]{\fontsize{\fsize}{#1\fsize}\selectfont}%
  \ifx\svgwidth\undefined%
    \setlength{\unitlength}{150.75001504bp}%
    \ifx\svgscale\undefined%
      \relax%
    \else%
      \setlength{\unitlength}{\unitlength * \real{\svgscale}}%
    \fi%
  \else%
    \setlength{\unitlength}{\svgwidth}%
  \fi%
  \global\let\svgwidth\undefined%
  \global\let\svgscale\undefined%
  \makeatother%
  \begin{picture}(1,1.0000003)%
    \lineheight{1}%
    \setlength\tabcolsep{0pt}%
    \put(0,0){\includegraphics[width=\unitlength,page=1]{ball2.pdf}}%
    \put(0.10199003,0.6111732){\color[rgb]{0,0,0}\makebox(0,0)[t]{\smash{\begin{tabular}[t]{c}$E$\end{tabular}}}}%
    \put(0.40049746,0.5862972){\color[rgb]{0,0,0}\makebox(0,0)[t]{\smash{\begin{tabular}[t]{c}$c(E)$\end{tabular}}}}%
  \end{picture}%
\endgroup%

%% file: figs/ball_mapped.pdf_tex
\begingroup%
  \makeatletter%
  \providecommand\color[2][]{%
    \errmessage{(Inkscape) Color is used for the text in Inkscape, but the package 'color.sty' is not loaded}%
    \renewcommand\color[2][]{}%
  }%
  \providecommand\transparent[1]{%
    \errmessage{(Inkscape) Transparency is used (non-zero) for the text in Inkscape, but the package 'transparent.sty' is not loaded}%
    \renewcommand\transparent[1]{}%
  }%
  \providecommand\rotatebox[2]{#2}%
  \newcommand*\fsize{\dimexpr\f@size pt\relax}%
  \newcommand*\lineheight[1]{\fontsize{\fsize}{#1\fsize}\selectfont}%
  \ifx\svgwidth\undefined%
    \setlength{\unitlength}{150.74998807bp}%
    \ifx\svgscale\undefined%
      \relax%
    \else%
      \setlength{\unitlength}{\unitlength * \real{\svgscale}}%
    \fi%
  \else%
    \setlength{\unitlength}{\svgwidth}%
  \fi%
  \global\let\svgwidth\undefined%
  \global\let\svgscale\undefined%
  \makeatother%
  \begin{picture}(1,0.99999994)%
    \lineheight{1}%
    \setlength\tabcolsep{0pt}%
    \put(0,0){\includegraphics[width=\unitlength,page=1]{ball_mapped.pdf}}%
    \put(0.11194024,0.51167126){\color[rgb]{0,0,0}\makebox(0,0)[t]{\smash{\begin{tabular}[t]{c}$\SPH$\end{tabular}}}}%
    \put(0.35074575,0.56142212){\color[rgb]{0,0,0}\makebox(0,0)[t]{\smash{\begin{tabular}[t]{c}$c(\SPH)$\end{tabular}}}}%
    \put(0.49999978,0.78530274){\color[rgb]{0,0,0}\makebox(0,0)[t]{\smash{\begin{tabular}[t]{c}$\BLLMC$\end{tabular}}}}%
  \end{picture}%
\endgroup%

%% file: figs/ball3.pdf_tex
\begingroup%
  \makeatletter%
  \providecommand\color[2][]{%
    \errmessage{(Inkscape) Color is used for the text in Inkscape, but the package 'color.sty' is not loaded}%
    \renewcommand\color[2][]{}%
  }%
  \providecommand\transparent[1]{%
    \errmessage{(Inkscape) Transparency is used (non-zero) for the text in Inkscape, but the package 'transparent.sty' is not loaded}%
    \renewcommand\transparent[1]{}%
  }%
  \providecommand\rotatebox[2]{#2}%
  \newcommand*\fsize{\dimexpr\f@size pt\relax}%
  \newcommand*\lineheight[1]{\fontsize{\fsize}{#1\fsize}\selectfont}%
  \ifx\svgwidth\undefined%
    \setlength{\unitlength}{150.75001504bp}%
    \ifx\svgscale\undefined%
      \relax%
    \else%
      \setlength{\unitlength}{\unitlength * \real{\svgscale}}%
    \fi%
  \else%
    \setlength{\unitlength}{\svgwidth}%
  \fi%
  \global\let\svgwidth\undefined%
  \global\let\svgscale\undefined%
  \makeatother%
  \begin{picture}(1,1.0000003)%
    \lineheight{1}%
    \setlength\tabcolsep{0pt}%
    \put(0,0){\includegraphics[width=\unitlength,page=1]{ball3.pdf}}%
    \put(0.0771144,0.62351365){\color[rgb]{0,0,0}\makebox(0,0)[lt]{\lineheight{1.25}\smash{\begin{tabular}[t]{l}$A$\end{tabular}}}}%
    \put(0.0771144,0.32500564){\color[rgb]{0,0,0}\makebox(0,0)[lt]{\lineheight{1.25}\smash{\begin{tabular}[t]{l}$B$\end{tabular}}}}%
    \put(0.25124375,0.59863765){\color[rgb]{0,0,0}\makebox(0,0)[lt]{\lineheight{1.25}\smash{\begin{tabular}[t]{l}$c(A)$\end{tabular}}}}%
    \put(0.25124373,0.37475765){\color[rgb]{0,0,0}\makebox(0,0)[lt]{\lineheight{1.25}\smash{\begin{tabular}[t]{l}$c(B)$\end{tabular}}}}%
  \end{picture}%
\endgroup%

%% file: figs/ball4.pdf_tex
\begingroup%
  \makeatletter%
  \providecommand\color[2][]{%
    \errmessage{(Inkscape) Color is used for the text in Inkscape, but the package 'color.sty' is not loaded}%
    \renewcommand\color[2][]{}%
  }%
  \providecommand\transparent[1]{%
    \errmessage{(Inkscape) Transparency is used (non-zero) for the text in Inkscape, but the package 'transparent.sty' is not loaded}%
    \renewcommand\transparent[1]{}%
  }%
  \providecommand\rotatebox[2]{#2}%
  \newcommand*\fsize{\dimexpr\f@size pt\relax}%
  \newcommand*\lineheight[1]{\fontsize{\fsize}{#1\fsize}\selectfont}%
  \ifx\svgwidth\undefined%
    \setlength{\unitlength}{150.74998807bp}%
    \ifx\svgscale\undefined%
      \relax%
    \else%
      \setlength{\unitlength}{\unitlength * \real{\svgscale}}%
    \fi%
  \else%
    \setlength{\unitlength}{\svgwidth}%
  \fi%
  \global\let\svgwidth\undefined%
  \global\let\svgscale\undefined%
  \makeatother%
  \begin{picture}(1,0.99999994)%
    \lineheight{1}%
    \setlength\tabcolsep{0pt}%
    \put(0,0){\includegraphics[width=\unitlength,page=1]{ball4.pdf}}%
  \end{picture}%
\endgroup%

%% file: figs/ball5.pdf_tex
\begingroup%
  \makeatletter%
  \providecommand\color[2][]{%
    \errmessage{(Inkscape) Color is used for the text in Inkscape, but the package 'color.sty' is not loaded}%
    \renewcommand\color[2][]{}%
  }%
  \providecommand\transparent[1]{%
    \errmessage{(Inkscape) Transparency is used (non-zero) for the text in Inkscape, but the package 'transparent.sty' is not loaded}%
    \renewcommand\transparent[1]{}%
  }%
  \providecommand\rotatebox[2]{#2}%
  \newcommand*\fsize{\dimexpr\f@size pt\relax}%
  \newcommand*\lineheight[1]{\fontsize{\fsize}{#1\fsize}\selectfont}%
  \ifx\svgwidth\undefined%
    \setlength{\unitlength}{228.26896242bp}%
    \ifx\svgscale\undefined%
      \relax%
    \else%
      \setlength{\unitlength}{\unitlength * \real{\svgscale}}%
    \fi%
  \else%
    \setlength{\unitlength}{\svgwidth}%
  \fi%
  \global\let\svgwidth\undefined%
  \global\let\svgscale\undefined%
  \makeatother%
  \begin{picture}(1,0.69742483)%
    \lineheight{1}%
    \setlength\tabcolsep{0pt}%
    \put(0,0){\includegraphics[width=\unitlength,page=1]{ball5.pdf}}%
    \put(0.26449048,0.34605963){\color[rgb]{0,0,0}\makebox(0,0)[lt]{\lineheight{1.25}\smash{\begin{tabular}[t]{l}$\vect{0}$\end{tabular}}}}%
    \put(0,0){\includegraphics[width=\unitlength,page=2]{ball5.pdf}}%
    \put(0.39591433,0.49391153){\color[rgb]{0,0,0}\makebox(0,0)[lt]{\lineheight{1.25}\smash{\begin{tabular}[t]{l}$L$\end{tabular}}}}%
  \end{picture}%
\endgroup%

%% file: figs/ball6.pdf_tex
\begingroup%
  \makeatletter%
  \providecommand\color[2][]{%
    \errmessage{(Inkscape) Color is used for the text in Inkscape, but the package 'color.sty' is not loaded}%
    \renewcommand\color[2][]{}%
  }%
  \providecommand\transparent[1]{%
    \errmessage{(Inkscape) Transparency is used (non-zero) for the text in Inkscape, but the package 'transparent.sty' is not loaded}%
    \renewcommand\transparent[1]{}%
  }%
  \providecommand\rotatebox[2]{#2}%
  \newcommand*\fsize{\dimexpr\f@size pt\relax}%
  \newcommand*\lineheight[1]{\fontsize{\fsize}{#1\fsize}\selectfont}%
  \ifx\svgwidth\undefined%
    \setlength{\unitlength}{150.74998807bp}%
    \ifx\svgscale\undefined%
      \relax%
    \else%
      \setlength{\unitlength}{\unitlength * \real{\svgscale}}%
    \fi%
  \else%
    \setlength{\unitlength}{\svgwidth}%
  \fi%
  \global\let\svgwidth\undefined%
  \global\let\svgscale\undefined%
  \makeatother%
  \begin{picture}(1,0.99999994)%
    \lineheight{1}%
    \setlength\tabcolsep{0pt}%
    \put(0,0){\includegraphics[width=\unitlength,page=1]{ball6.pdf}}%
  \end{picture}%
\endgroup%

%% file: discussion.tex
\section{Topics of Further Discussion}

\subsection{A case of Incorrect Attribution and its Widespread Adoption}
\label{sec:attribution}

Remarkably, even though Stefan Banach and Alfred Tarski are usually credited for the Banach-Tarski paradox, an equivalent statement by the same ideas and
construction had already been presented by Felix Hausdorff ten years earlier, in
\cite[pp.~469--473]{hausdorffbok} in 1914. It is
striking how much of their work consists of what Hausdorff had already published, including the same
ideas in the construction. This fact seems to have remained relatively unknown even to this day.

In his book, Hausdorff presented not just the Hausdorff paradox but also a result equivalent to the Banach-Tarski paradox and its proof. The reasoning was covered in \Cref{sec:proof} where we
saw how the weak Banach-Tarski paradox, \Cref{b-t2}, follows from the Hausdorff paradox,
\Cref{hausdorff}. To Hausdorff, this seems to have been a simple result,
following immediately from the far more intricate Hausdorff paradox, covered in \Cref{sec:hausdorff}. 

But beyond using the Hausdorff paradox itself, Banach and Tarski only credit one of their
constructions to
Hausdorff. Which is the usage of rotations to recover the countable subset of $\mathbb{S}^2$, as
covered in \Cref{bt_proof_1}. However, in the same book, \cite[pp.~472]{hausdorffbok}, Hausdorff
also briefly mentions the construction by \enquote{conical bodies} as subsets of the ball, of which
\Cref{bt_proof_2} is obtained. Again, this seems to be an almost obvious result to him and a detail
he only briefly explains in parenthesis while discussing the resulting existence of
non-measurable sets.

In his paper, Hausdorff does not mention an equivalence to the final proof of \Cref{sec:proof}:
adding the center of the ball. However, it is a trivial extra step using a simpler version of the same reasoning
already performed in \Cref{bt_proof_1}. It is also a result of no real significance to him since his goal
was to show the existence of non-measurable sets to which a single point of zero measure does not matter. However, Banach and Tarski seem to
have made the final step by extending the result from the ball to more general geometrical shapes,
giving it a more immediately striking result.

Although extending his result from a subset of the sphere to the entire ball is a necessary step
to show the existence of non-measurable
subsets, Hausdorff only seems to mention this construction in his book,
\cite{hausdorffbok}, and not in his article published the same year, \cite{hausdorffart}. The
latter only covers the construction of the sphere, the Hausdorff paradox. One could then suspect that Banach
and Tarski might have been unaware of Hausdorff's final construction. However, when citing his
paradox for the sphere, in \cite[pp.~261]{btart}, they cite his book and not the article, where
Hausdorff presented the 
same result Banach and Tarski would present themselves. We can conclude that there
is no reason to believe they were unaware of this fact and unfortunately did not fully credit their
results. Perhaps their intention was not to claim originality but only to present contemporary
discoveries. Of course, they did not coin the name, but in contrast, they did use the name
\enquote{Hausdorff's paradox} when referring to what we call the Hausdorff paradox. Perhaps
it is our fault that this error keeps repeating; it might be as simple as people tending to cite
\cite{hausdorffart} instead of \cite{hausdorffbok}. Or could it be that few ever take the time to
verify the often cited sources? 

It is even possible that Banach and Tarski ended up receiving more attention simply by releasing the
paper in French instead of German, which is instead what Hausdorff did for his article and book.
Also, the book by Hausdorff is almost 500 pages long, with an educative goal, making his results
less noticeable as their primary purpose is to demonstrate the existence of non-measurable sets.
Nevertheless, this is not to dismiss the work
by Banach and Tarski, such as Banach's version of Schröder-Bernstein, \Cref{thm:b-s-b}, and the
construction of the strong form of the paradox, \Cref{b-t1} by \Cref{equivalence_bt}. Their paper also covers many
other cases of
decomposition of sets in different dimensions, including the construction of Vitali. The work's
title translates to \enquote{On the decomposition of sets,} even though this specific result
certainly is the part that became well known. They also emphasized the more striking result of the
paradox, in contrast to Hausdorff, who used it only to prove the existence of
non-measurable sets. It should be noted that the latter result is more significant from a mathematical context, and
the proof of the Hausdorff paradox itself is far more revolutionary and sophisticated than the
results built on top of it. Nevertheless, the striking result that Banach and Tarski presented gave it a significant mainstream appeal by enabling colorful descriptions in layman's terms.


It is also worth noting that there seems to be much confusion regarding the various attributions in
the article by Banach and Tarski; for instance, many sources claim that they also credit Vitali
for the result of the Banach-Tarski paradox, which is incorrect. They do mention \enquote{an idea
from Vitali}, but it seems to relate to a different part of the article and not the construction of
the famous paradox and its results. In that part, they seem to be considering only one and two
dimensions. They do not clarify what this idea is, but it most likely refers to the Vitali set.

Regarding the originality of Hausdorff's result: outside of the obvious construction of Vitali,
whom Hausdorff credits among several others. In \cite[pp.~21]{wagonbok} Wagon also speculates that the results
relating to what we call \textit{free groups} might have been motivated by the work of Felix Klein, which had inspired similar reasoning in another result, the \textit{Sierpiński–Mazurkiewicz
paradox} released the same year as Hausdorff and containing many similar ideas to Hausdorff's work.
This paradox is also noteworthy as it creates a paradox by a construction that involves the same
idea of the lifting in \Cref{G_7}, yet it does never require the \nameref{ax_choice}.

To conclude, the Banach-Tarski paradox's weak form,  the one usually referred to as
the Banach-Tarski paradox, might deserve a change of name. An obvious choice could have
been to call it the Hausdorff paradox, but since that name has already been established, one might
suggest the \textit{Hausdorff-Banach-Tarski paradox} instead. However, the chance of a
new name gaining traction seems slim, considering the mainstream adoption of the current name and
the attribution it implies.

\raggedbottom
\pagebreak
\subsection{In Defence of the Axiom of Choice}
\flushbottom
\label{disc_ac}
In this thesis, we have seen several paradoxes with constructions relying on the \nameref{ax_choice}, both the so-called Vitali paradox in \Cref{ex:vitali} and the Hausdorff paradox in \Cref{hausdorff}, from which the two forms of the Banach-Tarski paradox in \Cref{b-t1,b-t2} follows. To repeat, it says:

\axchoice*

It is not uncommon to see the \nameref{ax_choice} blamed for these paradoxical results, mainly due
to its non-constructive nature. This is especially common in the case of the Banach-Tarski paradox
since the result so strikingly disobeys our expectations. And it is true that without the
\nameref{ax_choice} these results would be impossible under some conditions. If we assume the usual
Zermelo-Fraenkel axioms and some basic conditions, construction of any non-measurable set would be
impossible without the \ACC, as proven in \cite{solovayart}. This includes the results of Vitali and Hausdorff.

However, if we disagree with the paradoxical results, we need to examine all assumptions we put in. Before we declare the \ACC to be solely the cause of these results and possibly dismiss them as quirks caused by a controversial axiom, we need to widen our perspective on what the \ACC means to the larger mathematical field and consider what other assumptions enable the results.

The fundamental result from which the Banach-Tarski paradox is obtained is the Hausdorff paradox, which was covered in \Cref{sec:hausdorff}. There, we saw that the core results in its construction were not based on the \ACC but on finding rotations that enable the creation of a free group of rank $2$, in \Cref{G_6} and obtaining its paradoxicality in \Cref{G_5}. Thus paradoxicality was already obtained before the \ACC was used. The \ACC was only needed in \Cref{G_7} to turn the paradoxicality of the group into that of a set. However, as noted in \Cref{hausdorff_noac}, the \ACC is not even required for some sets.

Furthermore, there are plenty of other paradoxes not relying on the \ACC, with results still as
staggering as the Banach-Tarski paradox. Wagon lists some significant examples in
\cite[pp.~31]{wagonbok} and makes a much more comprehensive discussion on the subject of the \ACC
in \cite[Ch.~15]{wagonbok}. We also saw many other straightforward paradoxes that did not rely on
the \ACC: from the irrational rotation in \Cref{ex:irrational} to the Hyperwebster in
\Cref{ex:hyperwebster}. And many paradoxical results have come directly from the integers and their subsets throughout history.


The most important point comes from considering all the results the \ACC allows us to prove. We
have seen that removing the axiom would not prevent all paradoxes, even those of a geometrical
nature, but it would remove other significant results. Beyond not being able to show the existence of non-measurable subsets, we would also lose significant results such as \textit{Tychonoff's theorem}, the \textit{Hahn-Banach theorem}, the existence of the \textit{Stone–Čech compactification} for all spaces and the \textit{continuum hypothesis}, to name a few.

If we try to dismiss the \ACC to prevent certain paradoxes, we need to be careful: we can only dismiss the results of the \ACC if we are working in a very focused field where we do not need the \ACC. However, it seems unlikely that the paradoxes that depend on it would be significant in that situation. At least compared to the other fundamental paradoxes that would still hold. Indeed, the \ACC is typically taken as true today by most mathematicians.

A final point worth mentioning is a result presented by Kurt Gödel in 1938, in the form of his \textit{constructible universe} in \cite{godelart}. In it, he shows that for a suitable class of sets together with the Zermelo-Fraenkel system, both the \nameref{ax_choice} and the closely related Generalized Continuum Hypothesis must be true. Although not representing the entire class of sets, it does show that the \ACC is not inconsistent with our core axioms. In some ways, this is a strong argument for the validity of the axiom.

We have seen some reasons for keeping the \ACC and how rejecting it would fail to eliminate all paradoxes, even geometrical ones. Let us now broaden our discussion to other, often overlooked, assumptions that have received less scrutiny and how changing them could prevent these paradoxes. This includes paradoxes whose constructions do not depend on the \nameref{ax_choice}. There is little point in blaming the \ACC for the paradoxical results; instead, we should cherish the insight it offers us into other assumptions we put into our system and what we expect them to make possible.

We can observe that all mentioned paradoxes have something in common, not only those in focus of this thesis but also the irrational rotation of \Cref{ex:irrational}, the Hyperwebster in \Cref{ex:hyperwebster} as well as historical results such as \textit{Hilbert's hotel} and \textit{Galileo's paradox}. They all rely on \textit{infinite sets}. That is, they rely on the existence of sets of infinite cardinality, the \textit{axiom of infinity}. This property appears to lie at the heart of all veridical paradoxes, and removing it would prevent far more paradoxes than just removing the \ACC. However, doing so would also remove an even more significant part of our mathematical theory. It is a core part of Zermelo-Fraenkel set theory and a core part of modern mathematics. It would be impossible to discard the axiom of infinity without a major rework of centuries of mathematical foundations, as we would discard the existence of the set of integers.

Furthermore, removing the axiom of infinity would also affect all results that rely on the \ACC: since choices from finite collections of finite sets would not require anything outside the usual axioms. This would prevent most of the previously mentioned results that rely on the \ACC or, in many cases, make it impossible even to formulate them. Not allowing the existence of countably infinite sets not only seems to prevent all paradoxes but also many fundamental results. This is an axiom that would cause severe consequences if removed, even more so than the \ACC.

By accepting the axiom of infinity, we also have to accept some of its immediate paradoxical results. Consider the set of integers: the subset of all even integers is only half of all integers, yet it can be turned into all integers by simply dividing them by two. These kinds of results also conflict with our view of reality and are a consequence of countable infinity, but are something we learn to accept as part of the foundational assumptions that empower our mathematics. If we accept the axiom of infinity and the paradoxes it brings, the results enabled by the \ACC are not much more strange or harder to accept. The \ACC merely helps expand on, or perhaps visualize, the results caused by the axiom of infinity. Infinity gives us the core of the paradoxes, from the simple and immediate ones by Galileo and Hilbert to the Hyperwebster and the results by Vitali, Hausdorff, Banach and Tarski.

However, accepting the \ACC and axiom of infinity does not mean we have to accept the geometrical outcomes of the paradoxes. We now reach the final and most exciting observation: our usage of points in the geometrical paradoxes. This relates to another core assumption we regularly use but with minimal reflection. Compared to the \ACC and axiom of infinity, there is a construction in our usual geometrical and topological reasoning that never receives much scrutiny but still enables all mentioned paradoxes with geometrical results. Furthermore, it is an assumption that we typically use implicitly and often without realizing it.

We typically consider topology to describe properties of geometrical objects, including the study of open sets and mappings. However, we often overlook how the objects and their open sets are defined. When we consider a circle, we tend to define it as the set of all points of a given distance from its center. This way of describing objects leads to what is known as the \textit{point-set topology} or the \textit{general topology}. As the name implies, we consider geometrical objects as sets of points. This is a concept we often do not even give a name or explicitly state.

We are taught this concept early and we rarely question it. But by merely allowing constructing and
deconstructing geometrical objects by infinite sets of points, we allow something perhaps as
remarkable as the paradoxes themselves. This is the core mechanic used in all paradoxes of
geometrical nature mentioned in this thesis, including those not relying on the \ACC, such as the
irrational rotation. It is from the infinite collections of points, through the other axioms, that the paradoxes of Vitali, Hausdorff, Banach and Tarski appear, similar to what infinite collections enable in the other paradoxes, like the words in the Hyperwebster.

Describing objects by sets of points goes back to ancient times, and simply discarding this notion,
together with \textit{Euclidean geometry}, would present a monumental challenge. But what is a
\textit{point}? We rarely question this object. We consider it geometrical, yet it lacks any
geometrical property beyond representing a coordinate. The problem with the point-set construction
goes beyond merely the usage of points or the sets being infinite: the point-set approach seems to
lack aspects of the objects it tries to describe. We assume that a set of points can fully describe
a geometrical object, yet it is not the individual points of a set that describe an object to us but their broader structure together. So why should our mathematical description of objects only focus on the points?  \textit{Pointless topology}, or \textit{point-free topology}, brings us closer to this broader perspective. It is a fairly recent development, and the name pointless topology is typically used synonymously with the most popular such topology, called \textit{locales}.

Briefly, the theory of locales discards the notion of sets of points in the general topology. However, it keeps the idea of open sets, treated as the fundamental objects instead of the sets of points. The topological theory is developed from these building blocks on a purely algebraic foundation. This enables the creation of new topologies from unexpected algebraic concepts such as boolean algebra. As proposed in \cite{simpsonart}, we can apply locales to our existing notion of sets of points, basing their equivalence of open sets on those of the point-set topology. This generalizes the original topological space by building on it and gives further properties to objects beyond what the sets of points provide.

As locales can be defined without the use of points, a locale created from a traditional set
of points describe something more than the set from which it was created. Locales contain
\textit{sublocales} without points that are still nonempty, with properties such as measure, which
makes locales extend beyond the point-set definitions. In this perspective, locales from
the sets created in the construction of the Banach-Tarski paradox in \Cref{b-t2} will still overlap, unlike
the disjoint sets, and they will also be measurable. The sum of their measure will even
be greater than that of the original ball. Thus we see that the parts no longer describe a
decomposition like they originally did. Observe that this approach does not discard the use of
points of sets but rather the opposite: it expands on it and adds properties that allow more
scrutiny of the process of the paradox. It also helps increase measurability.

Locales allow us to avoid many geometrical paradoxes while keeping the \ACC. Furthermore, some
properties that require the \ACC in the point-set topology do not for locales. For example, a result
equivalent to the Tychonoff theorem can be proven without it, as shown in
\cite{preframeart}. Of course, it is still necessary for other proofs, but sometimes in a weaker
form. The possibility of avoiding the \ACC altogether also makes locales appealing in
\textit{constructive mathematics}, or simply \textit{constructivism}, in which a specific goal is
to avoid the \ACC. This should not be confused with Gödel's constructible universe. Despite the
similar names, observe the different perspectives: constructible math avoids the \ACC, whereas the
constructible universe tries to justify it. Despite their difference, this change of topology
appeals to both perspectives. Locales offer many other advantages. Some introduction to the subject
and its advantages can be found in \cite{pointlessart}.

The use of locales will not eliminate all paradoxes. As already stated, the axiom of infinity
enables paradoxes like the group in \Cref{G_5} and simpler constructions like the Hyperwebster and
those by Hilbert and Galileo. These results remain valid since none of them have a geometrical
result using points. Nevertheless, this change to our topology will prevent the more
striking geometrical results. Still, as long as we have the axiom of infinity, there will most
likely always be some paradoxical construction possible, even in a geometrical or topological
context. For example, we can obtain a somewhat paradoxical property by picking some basic shape like a
ball and repeatedly placing this shape at even distances along some axis. Moving all shapes one
step further by this distance can be seen to add and remove shapes. However, this result is less
geometrically significant, as the construction is unbounded, and its infinite volume remains
unchanged in the process.

To conclude, the \ACC can not easily be discarded and would not sufficiently prevent paradoxes.
However, we can prevent the geometrical paradoxes by Vitali, Hausdorff, Banach and Tarski by
adopting the pointless topology of locales, which also offers further advantages. One could compare
what the points are to topology right now with what infinitesimals once were to calculus: even
though the original objects can be used to some extent, the new approach reduces their limitations
and adds powerful new tools. We need to broaden our perspectives and look beyond the points.
However, just like the Riemann integral is a concept still taught in introductory calculus instead
of the Henstock–Kurzweil integral, we are taught early to think in terms of the point-set topology,
which we, in turn, pass on to others. Changing this trend would be challenging by its momentum
alone. However, locales is a relatively modern concept under research, and one can remain
optimistic about its future acceptance.

Finally, we need to consider the implication of the Banach-Tarski paradox itself. Are we too quick
to dismiss the results of the paradox? We could choose to accept it as a valid result in our
system. Does it even deviate from our reality if we just account for the change of energy by the change
of mass? We see this process all the time in radioactive materials and particle accelerators. It
does not seem more surprising than some observations made on the nature of our reality, like the
wave-particle duality.

When used to describe objects, the points fail to capture the concept of materials in reality by
not describing the underlying atoms and forces that keep matter together. Also, the decompositions in the
paradox involve parts strikingly disconnected from reality, essentially clouds of points.
However, that does not rule out the result on a smaller scale: \cite{augensteinart} is a well-known
work that speculates on the connection between the paradox, with infinite sets of points in
general, and properties in particle physics.

Furthermore, it would not be the first time math predicts some property of reality. From Maxwell predicting electromagnetic waves as a result of his own equations, to results such as $1+2+\dots=-1/12$ from the Zeta function finding practical use in string theory and quantum fields. In comparison, does it seem impossible that these paradoxes, and the infinite sets of points they are built on, could one day be used to describe aspects of our reality?

If nothing else, it still gives us insight into our own creation. Perhaps the strong reactions to
the paradox say more about our implicit assumption of mathematics and maybe an overconfidence in the
axioms we use. It is wrong to think of math as describing an absolute truth about the world.
Instead, it preserves the truth of what we put into it. Just like Gödel famously proved the
incompleteness of any sufficiently strong axiomatic system, we might have to accept that, despite
our best effort, no choice of axioms will ever truly allow us to describe reality as we know it.
Should we consider that a problem or an invitation to explore mathematics beyond the confinements of reality? That is left for the reader to decide.

%
%

I, the author, began working on this thesis without a personal opinion on the \ACC but have grown
to strongly support it, with a newfound skepticism to our trust in the points. The \ACC makes sense
in a world where infinite collections are daily occurrences. Any criticism of that notion should
instead be directed towards the axiom of infinity, one of the Zermelo-Fraenkel axioms. Likewise,
the Banach-Tarski paradox itself is not much of a paradox if we live in a world where objects
consist of infinitely small points that we can freely grasp and move as we want.

\section{Acknowledgements}

I want to thank Per Åhag for mentoring this thesis and, through one of his many interesting courses,
inspired me to look into the paradox. Happy to answer all my \enquote{tiny
questions} over the years. It was his encouragement that convinced me to write about this subject.

I also want to thank Lars-Daniel Öhman, who has been an invaluable mentor through my years of
university studies, answering both mathematical and study-related questions. And for
teaching many exciting courses, including in the area of discrete mathematics.